\documentclass{amsart}

\usepackage{amsmath, pdfsync, amssymb, changepage, fancyhdr,lipsum,etoolbox,graphicx,caption}
\usepackage{stmaryrd}
\usepackage{tikz}
\usepackage{pgfplots}
\usepackage{subfig}
\usepackage{mathtools}
\usepackage{amsthm}
\usepackage{amstext}
\usepackage{enumitem} 
\usepackage{mathrsfs} 
\usepackage{hyperref}
\usepackage{placeins}
\usepackage{float}
\usepackage[inner=22mm, outer=22mm, textheight=220mm]{geometry}

\newcommand\sfrac[2]{#1/#2}

\makeatletter
\newtheorem*{rep@theorem}{\rep@title}
\newcommand{\newreptheorem}[2]{%
	\newenvironment{rep#1}[1]{%
 		\def\rep@title{#2 \ref*{##1}}%
 		\begin{rep@theorem}}%
 	{\end{rep@theorem}}}
\makeatother

\newtheorem{theorem}{Theorem}[subsection]
\newreptheorem{theorem}{Theorem}
\theoremstyle{definition}
\newtheorem{definition}[theorem]{Definition}
\newtheorem{example}[theorem]{Example}
\theoremstyle{remark}
\newtheorem{remark}[theorem]{Remark}
\newtheorem{corollary}[theorem]{Corollary}

\numberwithin{equation}{section}

\DeclareMathOperator{\MCG}{MCG}
\DeclareMathOperator{\BMCG}{BMCG}
\DeclareMathOperator{\Homeo}{Homeo}
\DeclareMathOperator{\Isom}{Isom}
\DeclareMathOperator{\Hom}{Hom}
\DeclareMathOperator{\SL}{SL}
\DeclareMathOperator{\PSL}{PSL}
\DeclareMathOperator{\Id}{Id}
\DeclareMathOperator{\tr}{tr}
\DeclareMathOperator{\csch}{csch}

\def\C{\mathbb{C}}
\def\R{\mathbb{R}}
\def\N{\mathbb{N}}
\def\Z{\mathbb{Z}}
\def\H{\mathbb{H}}
\def\RPlus{\R_+}
\def\RPlusSpace{\R_+^3}
\def\SLZ{\SL_2(\Z)}
\def\SLR{\SL_2(\R)}
\def\SLC{\SL_2(\C)}
\def\PSLR{\PSL_2(\R)}

\def\Sym{\Sigma_3}
\def\SymE{\Sigma}
\def\L{\mathcal{L}}
\def\ML{\mathcal{ML}}

\def\SCC{\mathscr{S}}
\def\curves{\mathcal{C}}
\def\T{\mathcal{T}}

\def\surf{S_{g,n}}
\def\punc{S_{1,1}}
\def\kap{\kappa^{-1}}
\def\comp{\T_{\tr}}
\def\compL{\T_{\ell}}
\def\compF{\T_{\text{FN}}}
\def\RR{\RPlus\times\R}
\def\puncChar{\chi_{\SL_2(\C)}(\pi_1(\punc))}
\def\earthquakeTrace{E_\alpha^{\tr}}
\def\flowTrace{F_\alpha^{\tr}}

\def\earthquakeHyp{E_\alpha^{\ell}}
\def\earthquakeAFN{E_\alpha^{\text{FN}}}
\def\earthquakeBFN{E_\beta^{\text{FN}}}
\def\earthquakeCFN{E_{\alpha\beta}^{\text{FN}}}

\def\v{\mathbf{v}}
\def\w{\mathbf{w}}
\def\u{\mathbf{u}}
\def\x{x}
\def\y{y}
\def\z{z}
\def\f{F}

\def\sl{\text{sl}}

\def\ella{\ell_\alpha}
\def\taua{\tau_\alpha}
\def\ellfrac{\left(\frac{\ell}{2}\right)}
\def\taufrac{\left(\frac{\tau}{2}\right)}
\def\ssfrac{\left(\frac{s}{2}\right)}

\def\HR{H_1\left(\punc;\R\right)}

\begin{document}

\title{Earthquakes on the Once-Punctured Torus}

\author{Grace S. Garden}
%\address{School of Mathematics and Statistics, The University of Sydney, Sydney, NSW, 2006, Australia}
%\email{g.garden@maths.usyd.edu.au}
\thanks{School of Mathematics and Statistics, The University of Sydney, Sydney, NSW, 2006, Australia\\
Email address: g.garden@maths.usyd.edu.au}

\begin{abstract}
We study earthquake deformations on Teichm\"uller space associated with simple closed curves of the once-punctured torus. We describe two methods to get an explicit form of the earthquake deformation for any simple closed curve. The first method is rooted in linear recurrence relations, the second in hyperbolic geometry. The two methods align, providing both an algebraic and geometric interpretation of the earthquake deformations. We convert the expressions to other coordinate systems for Teichm\"uller space to examine earthquake deformations further. Two families of curves are used as examples. Examining the limiting behaviour of each gives insight into earthquakes about measured geodesic laminations, of which simple closed curves are a special case.
\end{abstract}

\maketitle

%% --- INTRODUCTION ------------------------------------------------------------------------------------------------------------------------------------------------------------------- 

\section{Introduction}
\label{sec:intro}

Earthquakes were first introduced by Thurston, generalising the concept of Dehn twists (see \cite{CanaryEtAl06} for an exposition). Earthquake deformations have proved to be useful in the exploration of hyperbolic manifolds and hyperbolic structures (for example \cite{Ker83,Ker85,Mir08,Thurston1,Thurston2,WatWol83,Wei89}). Most notably, a hyperbolic structure on a surface can be related to any other hyperbolic structure on a surface by a unique earthquake \cite{Thurston1,Thurston2}. This was used in the proof of the Nielsen realisation problem \cite{CanaryEtAl06,Ker83}. 

Consider a hyperbolic orientable surface $S=\surf$ with genus $g\geq0$ and $n\geq0$ punctures. Consider also a homotopy class of simple closed curves on $S$ represented by $\gamma$. Defined intuitively, a Dehn twist about $\gamma$ is achieved by cutting the surface open along $\gamma$, twisting around $\gamma$, and then gluing the surface together again. An earthquake deformation about $\gamma$ is achieved by cutting the surface open along $\gamma$, twisting around $\gamma$ by some distance $r$, and then gluing the surface together again.  In this way, earthquake deformations about simple closed curves can be viewed as ``fractional" Dehn twists.

Wolpert \cite{Wol83}, Kerckhoff \cite{Ker85}, and Weiss's \cite{Wei89} research into earthquakes reveal them to be specific integral curves. Wolpert builds on work from Fenchel-Nielsen twist deformations and associated vector fields, which are proven to be Hamiltonian using the Weil-Petersson symplectic form. Waterman and Wolpert \cite{WatWol83} continue this with a focus on the once-punctured torus, producing pictures of example earthquakes in this setting. Kerckhoff proves that earthquakes are real analytic, and that the earthquake flow is precisely the trajectories for an analytic flow defined on Teichm\"uller space. Weiss exploits parallels between the earthquake flow on Teichm\"uller space and geodesic flow on a Riemannian manifold to derive a form for the system of differential equations. 

We study left earthquake deformations on the once-punctured torus, denoted $\punc$, from an alternative perspective. We describe two methods in different coordinate systems to derive explicit formulas for the earthquake deformation for any simple closed geodesic (Section \ref{sec:results}). The first uses the character variety and linear recurrence relations to extend orbits of Dehn twists to a continuous flow. The second uses Teichm\"uller space and hyperbolic geometry to determine the earthquake deformation. The two methods give results that agree, proving the continuous flow in the first method is indeed the earthquake deformation and giving both algebraic and geometric interpretations, respectively, of earthquakes. We then focus on earthquake deformations in Fenchel-Nielsen length-twist coordinates, where we get more results for the {behaviour} and asymptotics. The formulas are later translated to other coordinate systems (Section \ref{sec:pics}). Results for simple closed geodesics can be extended to any measured geodesic lamination using a density argument from Thurston \cite{Thurston1}.

The fundamental group of $\punc$ is the free group of rank two, mapping class group of $\punc$ is $\SLZ$, and the universal cover of $\punc$ is the hyperbolic plane $\H^2$. The Teichm\"uller space of $\punc$, $\T(\punc)$, can be identified with $\H^2$, with isometry group $\PSLR$. The $\SLC$-character variety of $\punc$ is isomorphic to $\C^3$ via the following map (see \cite{Fricke1896,FrickeKlein1912,Horowitz75}),
\[
\begin{split}
	\Hom(\pi_1(\punc), \SLC)\sslash\SLC  &\to \C^3 \\
	\rho &\mapsto \begin{pmatrix} x \\ y \\ z \end{pmatrix} = \begin{pmatrix} \tr(\rho(\alpha)) \\ \tr(\rho(\beta)) \\ \tr(\rho(\alpha\beta)) \end{pmatrix},
\end{split}
\]
where $\alpha$ and $\beta$ are the chosen generators for $\pi_1(\punc)$. The choice of $\alpha$ and $\beta$ is called the framing. The $(x,y,z)$ coordinates in $\C^3$ are referred to as trace coordinates.

The action of Dehn twists on $\punc$ gives an action of polynomial automorphisms on the $\SLC$-character variety that preserve the level sets of a polynomial. This polynomial corresponds to the trace of the commutator $K=\alpha\beta\alpha^{-1}\beta^{-1}$,
\[
	\kappa(x,y,z) = \tr(\rho(K)) = x^2+y^2+z^2-xyz-2.
\]

The commutator is the generator for a peripheral subgroup. The value of the trace of the commutator represents the geometry of the surface around the puncture. In particular, the level set $\kap(-2)$ within the $\SLR$-character variety of $\punc$ is naturally identified with holonomies of complete hyperbolic structures of finite volume on $\punc$.

Let $\tilde{T}_\alpha$ be the polynomial automorphism associated with the left Dehn twist about $\alpha$. Using linear recurrence relations we can derive the following result extending the left Dehn twist about $\alpha$. Similar expressions can be derived for the left Dehn twists about $\beta$, and $\alpha\beta$.
\begin{reptheorem}{res:earthquakeTrace}
	For a fixed starting point $\v=(\x,\y,\z) \in \kap(t) \cap \RPlusSpace$ with $t \in \R$, the left Dehn twist about $\alpha$ can be extended to a continuous flow in trace coordinates given by the map $\flowTrace(\v): \R \mapsto \R^3$,
	\begin{equation*}
	\begin{split}
		\flowTrace(\v): \R 	&\to 		\R^3\\
		r 					&\mapsto 	 \begin{pmatrix} \x \\ y_{\alpha}(r) \\ y_{\alpha}(r-1) \end{pmatrix}
	\end{split}
	\end{equation*}
	where
	\begin{equation*}
	\begin{split}
		y_\alpha: \R 	&\to 		\R \\
		r 			&\mapsto 	
		\begin{cases} 
			C_+^{\alpha}(\v)S_+^{\alpha}(\v)^r+C_-^{\alpha}(\v)S_-^{\alpha}(\v)^r  \hspace{0.3cm}	\textnormal{ if } \x \neq 2 \\
			\y+\left(y-z\right)r \hspace{2.97cm}											\textnormal{ if } \x = 2  
		\end{cases}
	\end{split}
	\end{equation*}
	with 
	\begin{equation*}
	\begin{split}
		C_\pm^\alpha(\v) &= \frac{\y(\x^2-4)\pm(\x\y-2\z)\sqrt{\x^2-4}}{2(\x^2-4)}, \\
		S_\pm^\alpha(\v) &= \frac{\x\pm\sqrt{\x^2-4}}{2}.
	\end{split}
	\end{equation*}
	Moreover, $\flowTrace(\v)(n)=\tilde{T}_\alpha^n(\v)$, $\forall n \in \Z.$
\end{reptheorem}

We can also examine the action on Teichm\"uller space using a more geometric method. Cooper and Pfaff introduce coordinates called triangle lengths \cite{CooperPfaff}, which measure the half lengths of the unique geodesics associated with $\alpha, \beta, \alpha\beta$,
\[
	\begin{pmatrix} a \\ b \\ c \end{pmatrix} 
	= \begin{pmatrix} \frac{\ell(\alpha)}{2} \\[4pt] \frac{\ell(\beta)}{2} \\[4pt] \frac{\ell(\alpha\beta)}{2} \end{pmatrix} 
	= \begin{pmatrix} \cosh^{-1}\left(\frac{x}{2}\right) \\[4pt] \cosh^{-1}\left(\frac{y}{2}\right) \\[4pt] \cosh^{-1}\left(\frac{z}{2}\right) \end{pmatrix}.
\]

There is a similar relation on the coordinates here given by the collar equation,
\[
	\cosh(a)^2+\cosh(b)^2+\cosh(c)^2=-2\cosh(a)\cosh(b)\cosh(c).
\]

Let $\bar{T}_\alpha$ be the Dehn twist about $\alpha$ in triangle lengths. Using maximal collar neighbourhoods and hyperbolic geometry returns the following result for the left earthquake deformation about $\alpha$. Similar expressions can be derived for the left earthquake deformations about $\beta$ and $\alpha\beta$. 
\begin{reptheorem}{res:earthquakeLength}
	For a fixed starting point $\mathbf{w}=(a,b,c)\in\RPlusSpace$, the left earthquake deformation about $\alpha$ in triangle lengths is given by the map $\earthquakeHyp(\mathbf{w}):\R \to \RPlusSpace$,
	\begin{equation*}
	\begin{split}
		\earthquakeHyp(\mathbf{w}): \R 	&\to 		\RPlusSpace\\
		r 							&\mapsto 	\begin{pmatrix} a \\ b_\alpha(r) \\ b_\alpha(r-1) \end{pmatrix}
	\end{split}
	\end{equation*}
	where
	\begin{equation*}
	\begin{split}
		b_\alpha: \R 	&\to 		\R_+ \\
		r 			&\mapsto 	\cosh^{-1}\left(\cosh(b)\cosh(ra)-\left(\cosh(c)-\cosh(a)\cosh(b)\right)\frac{\sinh(ra)}{\sinh(a)}\right).
	\end{split}
	\end{equation*}
	Moreover, $\earthquakeHyp(\v)(n)=\bar{T}_\alpha^n(\v)$, $\forall n \in \Z.$
\end{reptheorem}

It is straightforward to show the expressions align on the level set $\kap(-2)$, proving that the continuous flow in Theorem \ref{res:earthquakeTrace} is the left earthquake deformation in trace coordinates. Then both Theorem \ref{res:earthquakeTrace} and Theorem \ref{res:earthquakeLength} provide explicit formulas for the left earthquake deformations for $\alpha$, $\beta$, and $\alpha\beta$ in the framing $\alpha$, $\beta$. The parameter $r$ represents the magnitude of the twist, and depends on a choice of orientation.

We can use an algorithm, referred to as the change of coordinates algorithm, to extend the results to the left earthquake deformation about any simple closed curve $\gamma$ in the same framing. To summarise the algorithm briefly, given $\gamma$ a simple closed curve as a word in $\alpha$ and $\beta$,
\begin{enumerate}
	\item Find $\delta$ such that $\pi_1(\punc) = \langle \gamma, \delta \rangle$;
	\item \label{item:alg} Apply Theorem \ref{res:earthquakeTrace} to $\gamma$ for framing $\gamma$, $\delta$;
	\item Convert the results to the framing $\alpha$, $\beta$.
\end{enumerate}
This same method provides explicit formulas for the left earthquake deformations for $\delta$ and $\gamma\delta$ at the same time. 

We convert the results for the left earthquake deformations to Fenchel-Nielsen length-twist coordinates for the pants decomposition given by $\alpha$, 
\[
	\begin{pmatrix} \ella \\ \taua \end{pmatrix} = \begin{pmatrix} 2\cosh^{-1}\left(\frac{x}{2}\right) \\ 2\cosh^{-1}\left(\frac{y}{2\coth\left(\cosh^{-1}\left(\frac{x}{2}\right)\right)}\right) \end{pmatrix} = \begin{pmatrix} 2a \\ 2\cosh^{-1}\left(\frac{\cosh(b)}{\coth\left(a\right)}\right) \end{pmatrix}.
\]

This is a natural coordinate system for Teichm\"uller space, and we find some results relating the asymptotics of the earthquake deformations and the initial lamination.

Denote a left earthquake deformation about $\gamma$ in Fenchel-Nielsen coordinates by
\begin{equation*}
\begin{split}
	E^\text{FN}_{\gamma}(\u,\cdot): \R &\mapsto \T(\punc) \\
	s &\mapsto \begin{pmatrix} \ell_\alpha^\gamma(\u,s) \\ \tau_\alpha^\gamma(\u,s) \end{pmatrix}.
\end{split}
\end{equation*}
\begin{reptheorem}{thm:earthquakeAsymptotics}
	Take simple closed curve $\gamma$. The slope of the associated earthquake $E^\text{FN}_{\gamma}(\u,s)$ converges to the inverse of the slope of $\gamma$ for all starting points $\u\in\T(\punc)$. We have, 
	\[
	\begin{split}
		\frac{\tau_{\alpha}^\alpha(\u,s)}{\ell_{\alpha}^\alpha(\u,s)}\,&\to\,-\frac{1}{\sl(\alpha)}=-\infty \hspace{0.45cm} \text{ as }s\to\infty, \, \forall \u\in\T(\punc) \text{ and }\\
		\frac{\tau_\alpha^{\gamma}(\u,s)}{\ell_\alpha^{\gamma}(\u,s)}\,&\to\,\frac{1}{\sl(\gamma)}\hspace{0.3cm}=\frac{i(\gamma,\beta)}{i(\gamma,\alpha)} \text{ as }s\to\infty, \, \forall \u\in\T(\punc).
	\end{split}
	\]
\end{reptheorem}
This is especially useful for studying earthquakes about measured geodesic laminations, which can be difficult to quantify.

Each of the left earthquake deformations corresponds to the action of the $1$-parameter groups of parabolic elements on Teichm\"uller space associated with the Dehn twists in $\SLZ$. We initially conjectured this would extend to an action of $\SLR$, but this does not appear to be the case. See also the work of \cite{AHW22}.

The paper proceeds as follows: Section \ref{sec:prelim} describes preliminary concepts and results, in particular regarding Dehn twists, earthquake deformations, the character variety, and Teichm\"uller space; Section \ref{sec:results} discusses the results for the earthquake deformations, first looking at the action on the character variety in Section \ref{sec:recur}, second looking at the action on Teichm\"uller space in Section \ref{sec:hypGeo}, third expanding on the change of coordinates algorithm in Section \ref{sec:algorithm}, fourth looking at the action in Fenchel-Nielsen coordinates, and fifth providing some examples in Section \ref{sec:examples}; Section \ref{sec:pics} presents and compares pictures of earthquakes in various coordinate systems. {Our examples include two families of curves and examining the limiting behaviour gives insight into earthquakes about measured geodesic laminations, of which simple closed curves are a special case.}

For the remainder of the paper we use $\alpha$ and $\beta$ the meridian and longitude respectively as the default framing for $\punc$. We often refer to an equivalence class of simple closed curves $[\gamma]$ by a representative curve $\gamma$. 

%% --- PRELIMINARIES ------------------------------------------------------------------------------------------------------------------------------------------------------------------- 

\section{Preliminaries}
\label{sec:prelim}

In this section, we provide relevant definitions and results regarding deformations, the character variety, and Teichm\"uller space with a focus on the once-punctured torus. 

%% --- DEHN TWISTS AND EARTHQUAKES -----------------------------------------------------------------------------------------------------------------------------

\subsection{Dehn twists and earthquakes}
\label{sec:twist}

We formalise the definition of Dehn twists and earthquakes, and introduce the mapping class group. Examples are provided for the once-punctured torus.

\begin{definition} [Left-handed Dehn twist, \cite{FarbMargalit12}]
	Let $\gamma$ be a simple closed curve on $S$. Find an annulus $A \subset S$ such that $\gamma$ is the core of $A$. Consider the annulus $\hat{A} = S^1\times[0,1] = \{ (\theta,r) \mid \theta \in [0, 2\pi),r \in [0,1]\}$. Identify $A$ with $\hat{A}$ by an orientation-preserving homeomorphism $\phi: A \to \hat{A}$ such that $\phi(\gamma(t)) = (2\pi t,\frac{1}{2})$.

	A \emph{left-handed Dehn twist} about $\gamma$, denoted $T_\gamma$, is then a twist of $A$ constructed by a map $\hat{T}_\gamma: \hat{A} \to \hat{A}$ defined by $(\theta, r) \mapsto (\theta + 2\pi r,r)$. Then $T_\gamma: S \to S$ is defined by 
	\begin{center}
		$T_\gamma(x)= \begin{cases} x & x \in S\setminus A, \\ \phi^{-1}\hat{T}_\gamma\phi(x) & x \in A. \end{cases}$
	\end{center} 
	This is continuous and $\hat{T}_\gamma\mid_{\partial\hat{A}} = \Id$.
\end{definition}

Adjusting the definition of $\hat{T}_\gamma$ to be $\hat{T}_\gamma(\theta,r) = (\theta - 2\pi r,r)$ gives the right-handed Dehn twist. It can be easily proven that Dehn twists are well-defined on homotopy equivalence classes of simple closed curves. For this paper we use Dehn twist to refer to a left-handed Dehn twist.

Dehn twists have a close relationship with the mapping class group. The \emph{mapping class group} of $S$ is the set of all isotopy classes of orientation-preserving homeomorphisms from $S$ to $S$,
\begin{equation*}
	\MCG(S) = \Homeo^+(S)/\sim.
\end{equation*}

\begin{theorem}[$\MCG$ is generated by Dehn twists]
	For an orientable closed surface of genus $g \geq 0$ a finite number of Dehn twists about non-separating simple closed curves generates $\MCG(S)$.
\end{theorem}

See the work of Dehn \cite{Dehn87}, Lickorish \cite{Lickorish64}, and Humphries \cite{Humphries79} for a proof and appropriate sets of generating curves.

There is a similar result for non-compact topologically finite surfaces. Take the subgroup of all orientation-preserving mapping classes that are the identity on the boundary of the surface, $\BMCG(S)$. If $S=\surf$ is an orientable surface then a finite number of Dehn twists about non-separating curves generates $\BMCG(S)$. For $\punc$ the simple closed curves to consider are precisely the meridian and longitude given by the framing, $\alpha$, $\beta$, with the mapping class group isomorphic to $\SLZ$.

Consider the Dehn twists for $\alpha$, $\beta$, $\alpha\beta$,
\begin{eqnarray}
	T_\alpha &=& 
	\begin{cases} 
		\alpha \mapsto \alpha \\ \beta \mapsto \beta\alpha^{-1}
	\end{cases}
	\hspace{.35cm}
	{T_\alpha}^{-1} = 
	\begin{cases} 
		\alpha \mapsto \alpha \\ \beta \mapsto \beta\alpha
	\end{cases} \\
	T_\beta &=& 
	\begin{cases} 
		\alpha \mapsto \alpha\beta \\ \beta \mapsto \beta
	\end{cases}
	\hspace{.7cm}
	{T_\beta}^{-1} = 
	\begin{cases} 
		\alpha \mapsto \alpha\beta^{-1} \\ \beta \mapsto \beta
	\end{cases} \\
	T_{\alpha\beta} &=& 
	\begin{cases} 
		\alpha \mapsto \alpha\beta\alpha \\ \beta \mapsto \alpha^{-1}
	\end{cases}
	\hspace{.24cm}
	{T_{\alpha\beta}}^{-1} = 
	\begin{cases} 
		\alpha \mapsto \beta^{-1} \\ \beta \mapsto \beta\alpha\beta
	\end{cases}
\end{eqnarray}

The above Dehn twists have respective representations $A, B, C \in \SLZ$, 
\begin{eqnarray}
	A	&=&	\begin{pmatrix} 1 & -1 \\ 0 & 1 \end{pmatrix}, \hspace{0.55cm} 	A^{-1}	=	\begin{pmatrix} 1 & 1 \\ 0 & 1 \end{pmatrix}, \\
	B	&=&	\begin{pmatrix} 1 & 0 \\ 1 & 1 \end{pmatrix}, \hspace{0.75cm} 	B^{-1}	=	\begin{pmatrix} 1 & 0 \\ -1 & 1 \end{pmatrix}, \\
	C	&=&	\begin{pmatrix} 2 & -1 \\ 1 & 0 \end{pmatrix}, \hspace{0.5cm} 	C^{-1}	=	\begin{pmatrix} 0 & 1 \\ -1 & 2 \end{pmatrix}.
\end{eqnarray}

The concept of Dehn twists can be extended using earthquake deformations to consider fractional twists around measured geodesic laminations. Loosely, a measured geodesic lamination on $S$ is a pair consisting of a foliation of a closed subset of $S$ by complete, simple geodesics and a transverse measure. 

We provided an intuitive definition for earthquake deformations about simple closed curves in Section \ref{sec:intro}. Simple closed curves equipped with a transverse measure, for instance the counting measure, are examples of measured geodesic laminations. A \emph{left earthquake deformation} about simple closed curve $\gamma$ with transverse measure $\mu$ is achieved by cutting the surface along $\gamma$, performing a left twist around $\gamma$ by some distance $t$, and then gluing the surface together again. Similarly a right earthquake can be obtained using a right twist. 

Typically, the measure $\mu$ determines the ``speed" of the twisting; a time $t$ earthquake deformation along $\gamma$ is a twist of distance $t\mu$. For a simple closed curve, if we ignore the measure, this is precisely the Fenchel-Nielsen twist deformation \cite{FenchelNielsen}. For rest of this paper we use earthquake deformation to refer to a left earthquake deformation. Where necessary we will sometimes distinguish between left or right earthquake deformations by forwards or backwards, respectively. 

Let $\ML$ denote the set of measured geodesic laminations on $S$ and $\SCC$ the set of isotopy classes of simple closed curves on $S$. Embed $\SCC\times\R_+$ in $\ML$ by sending $(\gamma,r)$ to the geodesic isotopic to $\gamma$ with $r$ times the counting measure.

\begin{theorem}[Thurston, \cite{Thurston1}]
\label{thm:dense}
	$\SCC\times\R_+$ is dense in $\ML$.
\end{theorem}

Consider the notion of intersection and slope of curves on $\punc$. Let $i(\gamma,\delta)$ be the algebraic intersection number between simple closed curves $\gamma$ and $\delta$. Curve $\gamma$ has \emph{slope} 
\begin{equation}
\label{eqn:slope}
	\sl(\gamma)=\frac{i(\gamma,\alpha)}{i(\gamma,\beta)}.
\end{equation}
Of course, this depends on the framing $\alpha,\beta$. The definition of slope can be extended to general measured geodesic laminations using the density result Theorem \ref{thm:dense}. Simple closed curves on $\punc$ have rational slope; measured geodesic laminations on $\punc$ that are not simple closed curves have irrational slope. 

We use Theorem \ref{thm:dense} to extend the earthquake deformation definition from simple closed curves to measured geodesic laminations. This approach is well-defined as per Kerckhoff \cite{Ker85}. That is, for any sequence of simple closed curves $\gamma_i$ with respective measures $\mu_i$, if $(\gamma_i,\mu_i)$ converges to $(\gamma,\mu)$ then the sequence of earthquakes about $(\gamma_i, \mu_i)$ converges to the earthquake about $(\gamma,\mu)$. 

For the purposes of this paper, we primarily explore simple closed curves and assume $\mu$ the counting measure, but Theorem \ref{thm:dense} means our study could be extended to all measured geodesic laminations. We discuss examples of this in Section \ref{sec:examples}.

Earthquakes deformations adjust the hyperbolic structure on the surface and can be considered in several different ways: 
\begin{enumerate}
	\item as deformations of $S$; 
	\item as maps from Teichm\"uller space to itself; 
	\item as paths in Teichm\"uller space as $t$ varies.
\end{enumerate}

We focus on the third perspective of earthquake deformations as paths in Teichm\"uller space. Let us end this section with a useful result for earthquakes that points to how important they are.

\begin{theorem}[Thurston's earthquake theorem]
	For any $\mathbf{p},\mathbf{q}\in\T(S)$, there exists a unique earthquake that connects them $\mathbf{p}$ to $\mathbf{q}$.  
\end{theorem}

%% --- CHARACTER VARIETIES AND TEICHMULLER SPACE ---------------------------------------------------------------------------------------------------------------------------------------------------------------

\subsection{Character varieties and Teichm\"uller space}
\label{sec:CharacterVarieties}

We first look at character varieties and then Teichm\"uller space. We review their definitions and then restate useful results for the once-punctured torus. In particular, we reintroduce trace coordinates and triangle lengths for the once-punctured torus.

The \emph{$G$-character variety} of $S$ is the space of conjugacy classes of homeomorphisms from $\pi_1(S)$ to an algebraic group $G$,
\begin{center}
	$\chi_G(\pi_1(S)) = \Hom(\pi_1(S), G)\sslash G$.
\end{center}
The double slash denotes the algebro-geometric quotient, which ensures the ensuing space is an algebraic set. Each element of the character variety represents a geometric structure on $S$. 

Consider the $\SLC$-character variety of $\punc$, $\puncChar$. As described in Section \ref{sec:intro}, $\puncChar$ is isomorphic to $\C^3$ via the following map (see \cite{Fricke1896,FrickeKlein1912,Horowitz75}),
\begin{equation}
\begin{split}
	\puncChar 	&\to 		\C^3\\
	\rho			&\mapsto 	\begin{pmatrix} x \\ y \\ z \end{pmatrix} = \begin{pmatrix} \tr(\rho(\alpha)) \\ \tr(\rho(\beta)) \\ \tr(\rho(\alpha\beta)) \end{pmatrix}.
\end{split}
\end{equation}
The trace coordinates $(x,y,z)$ are given in terms of the chosen framing, $\alpha$, $\beta$. 

\begin{figure}[hbt]
	\centering
	\includegraphics[width=0.5\textwidth]{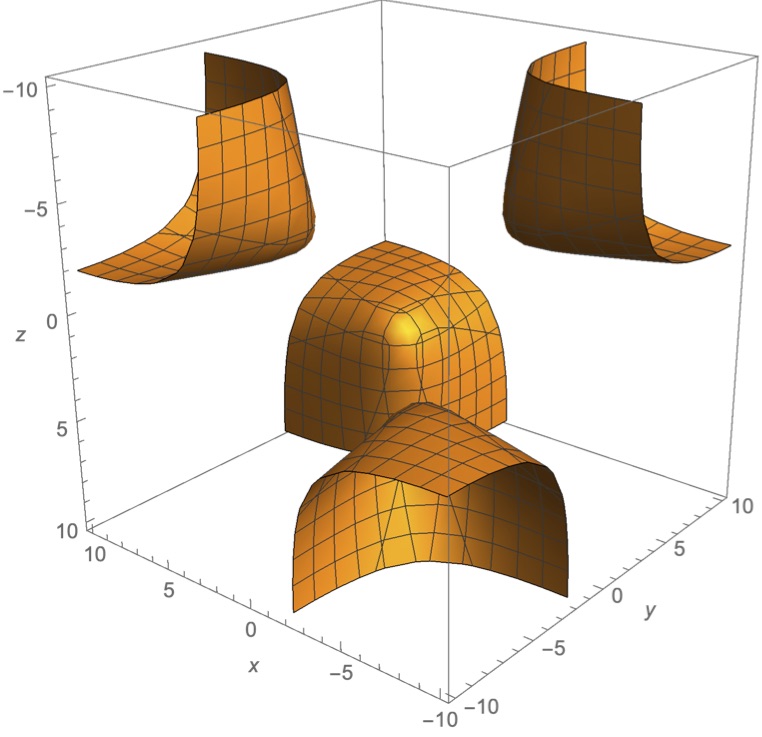}
	\caption{Level set $\kappa(x,y,z)=-2$ in $\R^3$.}
	\label{fig:levelSets}
\end{figure}

The commutator is $K= \alpha\beta\alpha^{-1}\beta^{-1}$. This is precisely the peripheral element around the puncture on $\punc$. For some element $\rho$ in the character variety, well-known trace identities can be used to get a polynomial expression for $\tr(\rho(K))$ in terms of $x,y,z$,
\begin{equation}
\label{eqn:commutator}
	\kappa(x,y,z) = \tr(\rho(K)) = x^2 + y^2 + z^2 - xyz - 2.
\end{equation}
An example level set of $\kappa(x,y,z)$ is given in Figure \ref{fig:levelSets}, see also \cite{Goldman03}. Level sets will be denoted $\kappa(x,y,z)=t$ or $\kap(t)$, with $t\in \R$. 

We first discuss the properties of this polynomial. If $\rho: \pi_1(\punc) \to \SLC$ corresponds to a geometric structure on $\punc$ then the value of $\tr(\rho(K))$ will say something about that structure. The different structures on $\punc$ are determined by the behaviour around the removed point, which can form a geodesic (with infinite volume), a cusp (with finite volume), or a cone point. The connection between $\tr(\rho(K))$ and the geometry can be understood via the classification of isometries on $\H^2$, see \cite{FarbMargalit12}.

Finally, there are some obvious symmetries of $\kappa(x,y,z)$.
\begin{remark}
\label{rem:kappaSym}
	The polynomial $\kappa(x,y,z)$ is invariant under the maps 
	\begin{equation*}
		I_{xy}, I_{xz}, I_{yz}: \C^3 \to \C^3
	\end{equation*}
	\begin{equation}
		I_{xy} = \begin{cases} x \mapsto -x \\ y \mapsto -y \\ z \mapsto z \end{cases}, \hspace{0.5cm}
		I_{xz} = \begin{cases} x \mapsto -x \\ y \mapsto y \\ z \mapsto -z \end{cases}, \hspace{0.5cm}
		I_{yz} = \begin{cases} x \mapsto x \\ y \mapsto -y \\ z \mapsto -z \end{cases}.
	\end{equation}
	Note the relationship with the Klein four-group
	\begin{equation*}
	 K_4 = \langle I_{xy}, I_{xz}, I_{yz}\rangle \cong \Z_2 \oplus \Z_2
	\end{equation*}
	where the group operation is composition of maps.
\end{remark}

We turn from the character variety to Teichm\"uller space. Teichm\"uller space can be defined loosely as the space that parameterises hyperbolic structures of a surface $S$. Consider the specific level set $\kap(-2)$. This has four connected components as pictured in Figure \ref{fig:levelSets}. The action of $K_4$ on these components corresponds to the action of $\Hom(\pi_1(\punc))$ on $\puncChar$ and hence we may identify the component in the positive octant with the Hitchin-Teichm\"uller component in the $\PSLR$-character variety of $\punc$. Label this component $\comp$. The Teichm\"uller space of $\punc$ is itself identified with $\H^2$.

We use the Poincar\'e half-plane model of $\H^2$. The orientation-preserving isometry group of $\H^2$ is $\Isom^+(\H^2)=\PSLR$, given via M\"obius transformations. Without loss of generality, assume that $\gamma$ is the unique geodesic associated with the equivalence class $[\gamma]$. Let $M$ be the M\"obius transformation associated with $\gamma$. We can derive the equation for the length of ${\gamma}$, $\ell(\gamma)$, in terms of the trace of $M$,
\begin{equation}
\label{eqn:lengthTrace}
	\ell(\gamma) = 2\cosh^{-1}\left(\frac{\tr(M)}{2}\right). 
\end{equation}

This enables new coordinates that provide a geometric way to parametrise Teichm\"uller space.
\begin{definition}
	The \emph{triangle lengths} \cite{CooperPfaff} for the once-punctured torus are the half-lengths of the geodesics associated with generators $\alpha$, $\beta$ and their product $\alpha\beta$.
	\begin{equation}
	\label{eqn:lengthCoord}
		\begin{pmatrix} a \\ b \\ c \end{pmatrix} 
		\coloneqq \begin{pmatrix} \frac{\ell(\alpha)}{2} \\[4pt] \frac{\ell(\beta)}{2} \\[4pt] \frac{\ell(\alpha\beta)}{2} \end{pmatrix} 
		= \begin{pmatrix} \cosh^{-1}\left(\frac{x}{2}\right) \\[4pt] \cosh^{-1}\left(\frac{y}{2}\right) \\[4pt] \cosh^{-1}\left(\frac{z}{2}\right) \end{pmatrix}
	\end{equation}
where $\left(x,y,z\right)$ are the trace coordinates.
\end{definition}

The triangle lengths $(a,b,c)$ satisfy the collar equation (see \cite{CooperPfaff}),
\begin{equation}
\label{eqn:collar}
	\cosh^2{a}+\cosh^2{b}+\cosh^2{c}=2\cosh{a}\cosh{b}\cosh{c}.
\end{equation}
This is the equivalent to the polynomial $\kappa(x,y,z)=-2$ in triangle length coordinates. Label the surface corresponding to the collar equation in $\RPlus^3$ as $\comp \cong \T(\punc)$.

%% --- RESULTS ---------------------------------------------------------------------------------------------------------------------------------------------------------------

\section{Results}
\label{sec:results}

In this section, we examine the action of Dehn twists on the character variety using trace coordinates and linear recurrence relations and on Teichm\"uller space using triangle lengths and hyperbolic geometry. Through this, we obtain a formula for the earthquake deformation for $\alpha$. Expressions for $\beta$ and $\alpha\beta$ can be similarly derived. We then provide an algorithm to determine a similar formula for any simple closed curve and discuss two examples of families of curves.

%% --- LINEAR RECURRENCE RELATIONS ---------------------------------------------------------------------------------------------------------------------------------------------------

\subsection{Studying deformations using linear recurrence relations}
\label{sec:recur}

We explore the action of Dehn twists on the character variety first. The expression $\kappa$ from Equation \ref{eqn:commutator} is invariant under the action of the mapping class group of the once-punctured torus. In fact, the action of $\MCG(\punc)$ on $\puncChar$ is equivalent to the action of the group $\Gamma$ on $\C^3$ \cite{Horowitz75}, where 
\begin{equation*} 
	\Gamma = \{\text{polynomial automorphisms }p\text{ on }\C^3 \mid \kappa(p(\mathbf{v}))=\kappa(\mathbf{v}),\mathbf{v}\in\C^3\}.
\end{equation*}

Write the polynomial automorphism associated with mapping class $f$ as $\tilde{f}$. Using trace identities, we can write the expressions for $T_\alpha$, $T_\beta$, $T_{\alpha\beta}$ and their inverses in terms of the trace coordinates.
\begin{equation} 
\label{eqn:polyAlpha}
	\hspace{.225cm} 
	\tilde{T}_\alpha(x,y,z) = (x,xy-z,y),
	\hspace{.5cm}
	\tilde{T}_\alpha^{-1}(x,y,z) = (x,z,xz-y),
\end{equation}
\vspace{-.6cm}
\begin{equation} 
\label{eqn:polyBeta}
	\hspace{.225cm} 
	\tilde{T}_\beta(x,y,z) = (z,y,yz-x),
	\hspace{.5cm}
	\tilde{T}_\beta^{-1}(x,y,z) = (xy-z,y,x),
\end{equation}
\vspace{-.6cm}
\begin{equation} 
\label{eqn:polyAlphaBeta}
	\tilde{T}_{\alpha\beta}(x,y,z) = (xz-y,x,z),
	\hspace{.4775cm}
	\tilde{T}_{\alpha\beta}^{-1}(x,y,z) = (y,yz-x,z).
\end{equation}

\begin{remark}
\label{rem:symAuto}
	The polynomial automorphisms are related by conjugation with elements from the symmetric group on three elements, $\Sym$. Specifically, consider rotation and reflections given by
	\begin{eqnarray*}
		\SymE_\text{Rot}(x,y,z) 	&=& (z,x,y),\\
		\SymE_\text{Ref}^1(x,y,z) &=& (x,z,y),\\
		\SymE_\text{Ref}^2(x,y,z) &=& (z,y,x),\\
		\SymE_\text{Ref}^3(x,y,z) &=& (y,x,z),
	\end{eqnarray*}
	where $\SymE_\text{Ref}^i$ is the reflection that fixes the $i$th coordinate and switches the other two coordinates. Note that these all preserve the {chosen component $\comp$.}

	The polynomial automorphisms $\tilde{T}_{\alpha}$, $\tilde{T}_{\beta}$, $\tilde{T}_{\alpha\beta}$ are related by rotations,
	\begin{equation*}
		\tilde{T}_{\alpha}(x,y,z) 
		= \SymE_\text{Rot}^{-1} \tilde{T}_{\beta} \SymE_\text{Rot} (x,y,z),
		= \SymE_\text{Rot} \tilde{T}_{\alpha\beta} \SymE_\text{Rot}^{-1} (x,y,z),
	\end{equation*}
and the polynomial automorphisms and their inverses are related by reflections,
	\begin{eqnarray*}
		\tilde{T}_{\alpha}^{-1}(x,y,z) 		&=& 	\SymE_\text{Ref}^1 	\tilde{T}_{\alpha} 		\SymE_\text{Ref}^1 (x,y,z), \\
		\tilde{T}_{\beta}^{-1}(x,y,z) 		&=& 	\SymE_\text{Ref}^2 	\tilde{T}_{\beta} 		\SymE_\text{Ref}^2 (x,y,z), \\
		\tilde{T}_{\alpha\beta}^{-1}(x,y,z) 	&=& \SymE_\text{Ref}^3 	\tilde{T}_{\alpha\beta} 	\SymE_\text{Ref}^3 (x,y,z). 
	\end{eqnarray*}
\end{remark}
As a consequence, any formula for the earthquake deformation or a flow associated with one polynomial automorphism for $\alpha$, $\beta$, or $\alpha\beta$ can be used to derive formulas for all three polynomial automorphisms using the action of $\Sym$. 

Without loss of generality, consider the polynomial automorphism $\tilde{T}_\alpha$ associated with the Dehn twist about $\alpha$. For a point $\v\in\kap(t) \cap \R^3$ we get an orbit of points, given by $\tilde{T}_\alpha^n(\v)$ with $n \in\Z$. Investigating the action of $\tilde{T}_\alpha^n$ and using results from linear recurrence relations, we get the following theorem. 

\begin{theorem}
\label{res:earthquakeTrace}
	For a fixed starting point $\v=(\x,\y,\z) \in \kap(t) \cap \RPlusSpace$ with $t \in \R$, the Dehn twist about $\alpha$ can be extended to a continuous flow in trace coordinates given by the map $\flowTrace(\v): \R \mapsto \R^3$,
	\begin{equation*}
	\begin{split}
		\flowTrace(\v): \R 	&\to 		\R^3\\
		r 				&\mapsto 	 \begin{pmatrix} \x \\ y_{\alpha}(r) \\ y_{\alpha}(r-1) \end{pmatrix}
	\end{split}
	\end{equation*}
	where
	\begin{equation*}
	\begin{split}
		y_\alpha: \R 	&\to 		\R\\
		r 			&\mapsto 	
		\begin{cases} 
			C_+^{\alpha}(\v)S_+^{\alpha}(\v)^r+C_-^{\alpha}(\v)S_-^{\alpha}(\v)^r  \hspace{0.3cm}	\textnormal{ if } \x \neq 2 \\
			\y+\left(y-z\right)r \hspace{2.97cm}											\textnormal{ if } \x = 2  
		\end{cases}
	\end{split}
	\end{equation*}
	with 
	\begin{equation*}
	\begin{split}
		C_\pm^\alpha(\v) &= \frac{\y(\x^2-4)\pm(\x\y-2\z)\sqrt{\x^2-4}}{2(\x^2-4)}, \\
		S_\pm^\alpha(\v) &= \frac{\x\pm\sqrt{\x^2-4}}{2}.
	\end{split}
	\end{equation*}
	Moreover, $\flowTrace(\v)(n)=\tilde{T}_\alpha^n(\v)$, $\forall n \in \Z.$
\end{theorem}

\begin{proof}
	Take $n\in\Z$. We introduce notation for the three components of the map $\tilde{T}_\alpha^n$ acting on points $\v=(\x,\y,\z) \in\kap(t) \cap \R^3$,
	\[
		\tilde{T}_\alpha^n(\v) = \begin{pmatrix} x_\alpha(n) \\ y_\alpha(n) \\ z_\alpha(n) \end{pmatrix}.
	\]
	Note $x_\alpha$, $y_\alpha$, $z_\alpha$ are also dependant on $\v$, although this is omitted in the notation for brevity. 
	The polynomial automorphism for $\tilde{T_\alpha}$ (Equation \ref{eqn:polyAlpha}) gives a system of recurrence relations,
	\[
		\begin{pmatrix} x_{\alpha}(n) \\ y_{\alpha}(n) \\ z_{\alpha}(n) \end{pmatrix} 
		= \begin{pmatrix} \x \\ \x y_{\alpha}(n-1) - y_{\alpha}(n-2)  \\ y_{\alpha}(n-1) \end{pmatrix}.
	\]
	This gives a linear recurrence relation for $y_\alpha$,
	\[
		y_\alpha(n)=\x y_\alpha(n-1) - y_\alpha(n-2),
	\]
with initial conditions $y_\alpha(0)=\y$ and $y_\alpha(1)=\x\y-\z$.

	The recurrence relation can be solved using the characteristic polynomial method for linear recurrence relations. There will be two forms of solution depending on if there are two distinct roots (which occurs if $\x \neq 2$) or one repeated root (which occurs if $\x = 2$) to the characteristic polynomial. Extending the solutions from $n\in\Z$ to any real number $r\in\R$ and considering $\v=(\x,\y,\z) \in \kap(t) \cap \RPlusSpace$ with $t \in \R$ gives the expression for $\flowTrace$.
\end{proof}

We are particularly interested in the component $\comp = \kap(-2)\cap \R_+^3$ mentioned in Section \ref{sec:CharacterVarieties}, which can be identified with the Teichm\"uller space of $\punc$. Note the special case $x=2$ does not exist on this component. A similar expression can be derived to extend the Dehn twists about $\beta$ and $\alpha\beta$ to flows $F_\beta^{\tr}$ and $F_{\alpha\beta}^{\tr}$ respectively. This is done either by following the same derivation or by applying the symmetry relationships from Remark \ref{rem:symAuto} to $\flowTrace$.

There are a few remarks we can make about Theorem \ref{res:earthquakeTrace}. 

\begin{remark}
	The results in Theorem \ref{res:earthquakeTrace} are given for the positive octant $\v=(\x,\y,\z) \in \kap(t) \cap \RPlusSpace$, $t \in \R$. The results could be extended for all starting points $\v=(\x,\y,\z) \in \kap(t) \cap \R^3$ using the symmetries given in Equation \ref{rem:kappaSym}. 
\end{remark}

\begin{remark}
	The right-handed Dehn twist about $\alpha$, $T_\alpha^{-1}$, can similarly be extended to a continuous flow. This is given by the formula in Equation \ref{res:earthquakeTrace} under the transformation $r \mapsto -r$, with
	\[
		\flowTrace\left(\flowTrace(\v)(r)\right)(-r)=\v=\flowTrace\left(\flowTrace(\v)(-r)\right)(r),
	\]
	as expected.
\end{remark}

\begin{remark}
	Note where imaginary numbers could appear in $y_\alpha$, $C_\pm^\alpha$, $S_\pm^\alpha$.
	\begin{enumerate}[label=(\roman*)]
		\item If $r\in\Z$ then $y_\alpha(r) \in \R$.
	\end{enumerate}
	Assume $r\in\R$, $\v=(x,y,z) \in \kap(t) \cap \RPlusSpace$ with $t \in \R$,
	\begin{enumerate}[label=(\roman*),resume]
		\item If $x\geq2$, then $C_\pm^\alpha(\v), S_\pm^\alpha(\v) \in\R \text{ } \forall y,z\in\RPlus$ and $y_\alpha \in \R$, thus $\flowTrace(\v)(r)\in\R^3$;
		\item 				
		\label{rem:squareRoot1}
			If $0\leq \x < 2$, then $C_\pm^\alpha(\v), S_\pm^\alpha(\v), \notin\R \text{ } \forall y,z\in\RPlus$, but $y_\alpha \in \R$, and thus $\flowTrace(\v)(r)\in\R^3$;
		\item If $\v=(x,y,z) \in \comp$, then $\flowTrace(\v)(r)\in\comp$.
	\end{enumerate}
\end{remark}

\begin{remark}
	Consider the limit of the flow as $r\to\pm\infty$,
	\[
		\left(x_{\pm\infty}, y_{\pm\infty}, z_{\pm\infty}\right) \coloneqq \lim_{r\to\pm\infty}\flowTrace(x,y,z)(r) = \left(x, \infty, \infty\right),
	\]
	and define 
	\[
		N_\pm(u,v,w)=v\pm\frac{2w-uv}{\sqrt{u^2-4}}.
	\]

	Then projectively,
	\[
	\begin{split}
		\left[x_{\pm\infty}:y_{\pm\infty}:z_{\pm\infty}\right] = \left[0:\frac{N_\mp(x,y,z)}{\sqrt{N_\mp(x,y,z)^2+N_\pm(x,z,y)^2}}:\frac{N_\pm(x,z,y)}{\sqrt{N_\mp(x,y,z)^2+N_\pm(x,z,y)^2}}\right].
	\end{split}
	\]

	The projective classes provide an idea of the rate of blow up at the limit. In particular, the projective classes as $r\to\infty$ and $r\to-\infty$ are distinct,
	\[
		\left[x_{\infty}: y_{\infty}: z_{\infty}\right] \neq \left[x_{-\infty}: y_{-\infty}: z_{-\infty}\right], {\forall (\x,\y,\z)\in\RPlusSpace.}
	\]
\end{remark}

\begin{figure}[t!]
	\centering
	\includegraphics[width=0.5\textwidth]{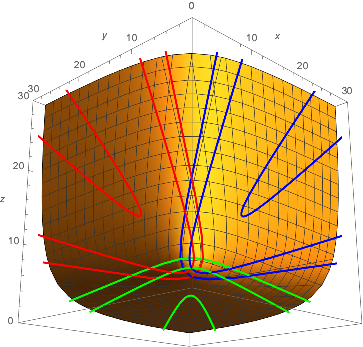}
	\caption{Example forward and backward earthquake deformations about $\alpha$ (red), $\beta$ (blue), and $\alpha\beta$ (green) on $\comp$. Starting points are given by the respective sets $\mathcal{S}_\alpha$, $\mathcal{S}_\beta$, and $\mathcal{S}_{\alpha\beta}$ from Equation \ref{eqn:startPoint}.}
	\label{fig:orbitsTrace}
\end{figure}

The approach used in the proof of Theorem \ref{res:earthquakeTrace} could potentially be applied to other simple closed curves as long as the Dehn twist is known. Whether the recurrence relation can be solved for all simple closed curves is unclear.

See Figure \ref{fig:orbitsTrace} for some example flows with starting points,
\begin{equation}
\label{eqn:startPoint}
\begin{split}
	\mathcal{S}_\alpha = 		&\left\{\left(3,3,3\right),\left(2\sqrt{2},2\sqrt{2},4\right),\left(10,10,-10\left(-5+\sqrt{23}\right)\right)\right\}, \\
	\mathcal{S}_\beta = 			&\SymE_\text{Rot}^{-1} \mathcal{S}_\alpha, \\
	\mathcal{S}_{\alpha\beta} = 	&\SymE_\text{Rot}\mathcal{S}_\alpha.
\end{split}
\end{equation}

The results presented in this section have extended the orbits of the Dehn twists about $\alpha$, $\beta$, and $\alpha\beta$ to a continuous flow in trace coordinates.

%% --- HYPERBOLIC GEOMETRY ---------------------------------------------------------------------------------------------------------------------------------------------------

\subsection{Studying deformations using hyperbolic geometry}
\label{sec:hypGeo}

We now explore the action of Dehn twists on Teichm\"uller space using triangle lengths. Unless otherwise stated, use $\alpha$, $\beta$, and $\alpha\beta$ to denote the unique geodesic associated with each equivalence class and let $\bar{T}_{\alpha}$, $\bar{T}_{\beta}$, $\bar{T}_{\alpha\beta}$ be the respective Dehn twists in triangle lengths. To analyse the Dehn twists, we turn to maximal collar neighbourhoods as the neighbourhood to twist. 

A collar neighbourhood of the geodesic $\gamma$ is an $\epsilon$-neighbourhood of $\gamma$ that is homeomorphic to an open annulus. A \emph{maximal collar neighbourhood} of $\gamma$ is a collar neighbourhood of $\gamma$ that is maximal with respect to inclusion of sets. We can generate maximal collar neighbourhoods around $\alpha$, $\beta$, and $\alpha\beta$ on $\punc$. 

\begin{figure}[htb]
	\centering
	\begin{minipage}{0.44\textwidth}
		\flushright
		\includegraphics[width=0.68\textwidth]{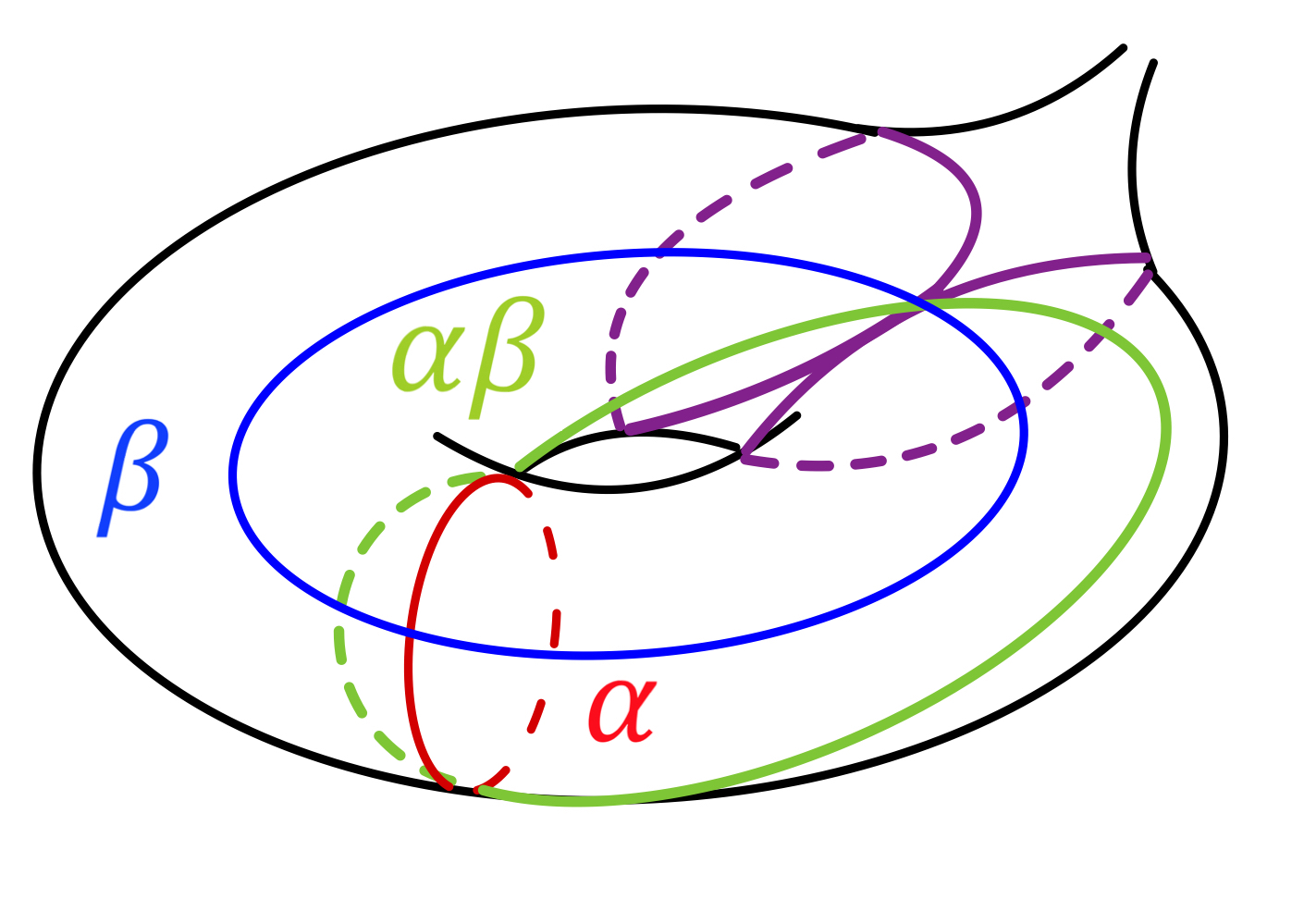}
	\end{minipage}
	\begin{minipage}{0.1\textwidth}
		\centering
		\center{$\implies$}
	\end{minipage}
	\begin{minipage}{0.44\textwidth}
		\flushleft
		\includegraphics[width=0.8\textwidth]{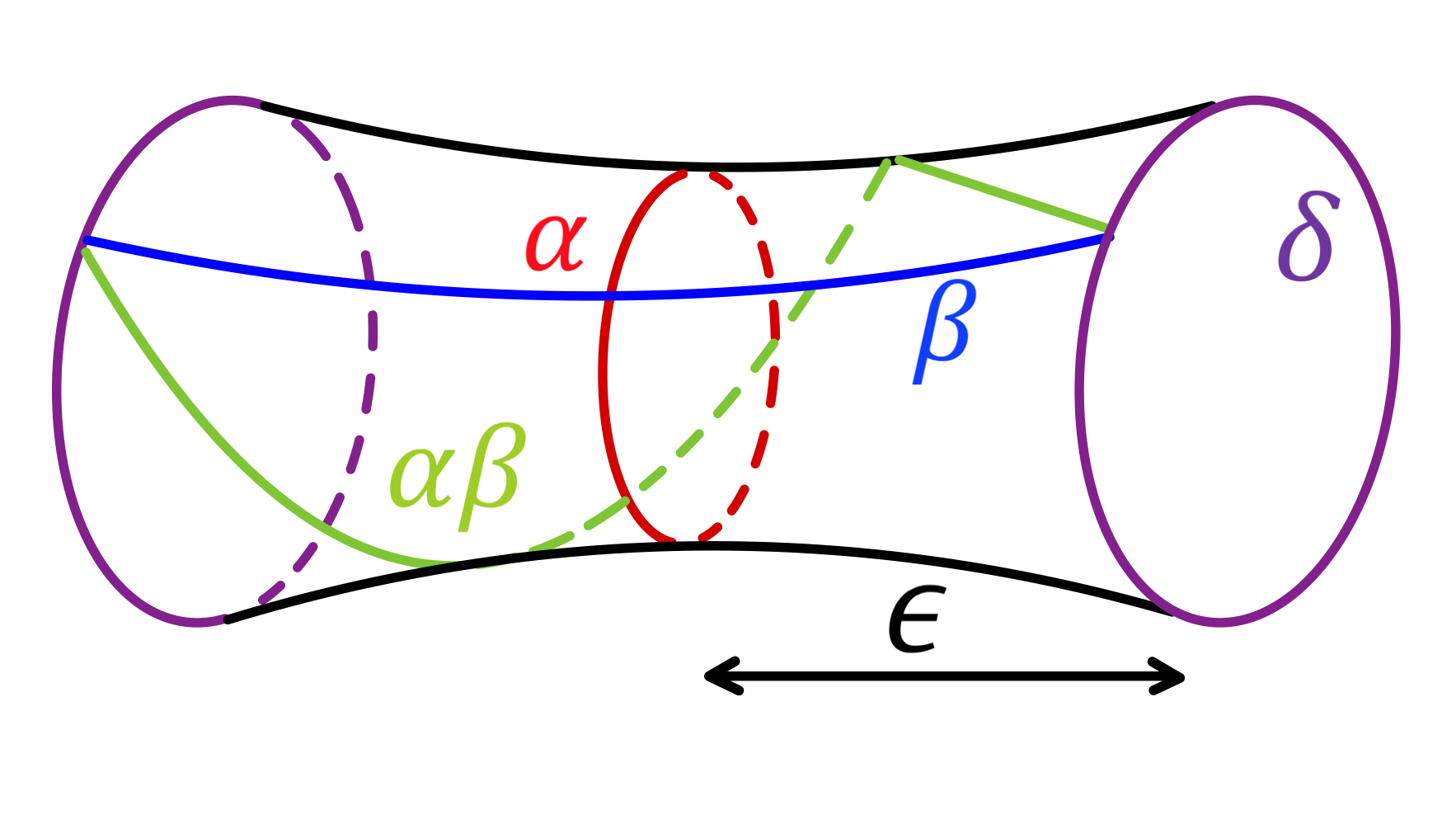}
	\end{minipage}
	\caption{Maximal collar neighbourhood for $\alpha$ on $\punc$ (left) and unfolded into a tube (right).}
	\label{fig:maxCollar}
\end{figure}

For example, the maximal collar neighbourhood for $\alpha$ is shown in Figure \ref{fig:maxCollar}.

\begin{theorem}[Buser, \cite{Buser78}]
\label{res:Buser}
	For a maximal collar neighbourhood of geodesic $\gamma$ on the once-punctured torus, the length of the boundary curve $d=\ell(\delta)$ and the length of half of the neighbourhood $\epsilon$ are given by
	\begin{eqnarray*}
		d(\gamma) 		&=& 2\coth{\frac{\ell(\gamma)}{2}},\\
		\epsilon(\gamma) 	&=& \sinh^{-1}\left(\frac{1}{\sinh{\frac{\ell(\gamma)}{2}}}\right).
	\end{eqnarray*}
\end{theorem}

\begin{figure}[thb]
	\centering
	\begin{minipage}{0.44\textwidth}
		\flushright
		\includegraphics[width=0.78\textwidth]{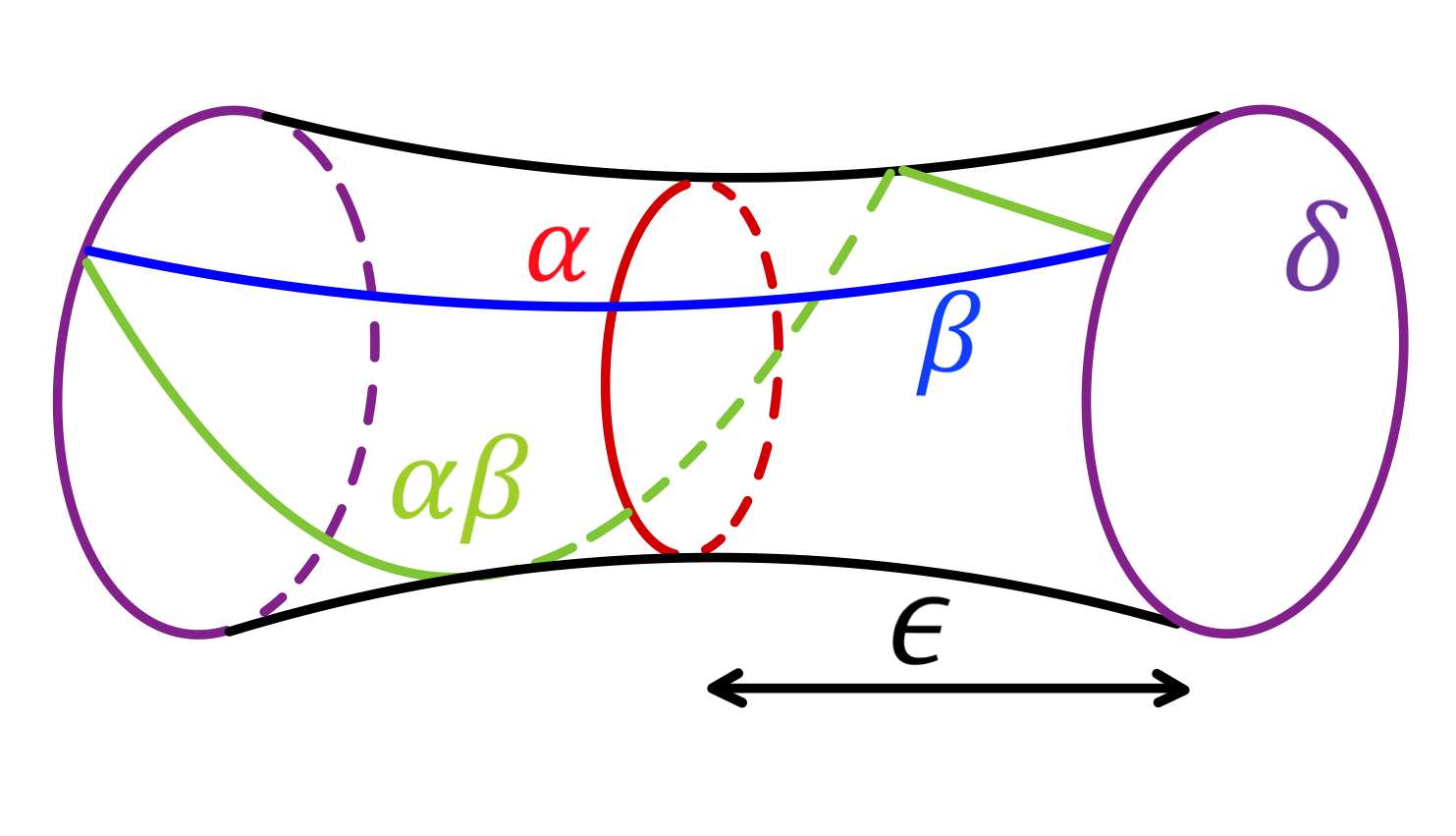}
	\end{minipage}
	\begin{minipage}{0.1\textwidth}
		\centering
		\center{$\implies$}
	\end{minipage}
	\begin{minipage}{0.44\textwidth}
		\flushleft
		\includegraphics[width=0.72\textwidth]{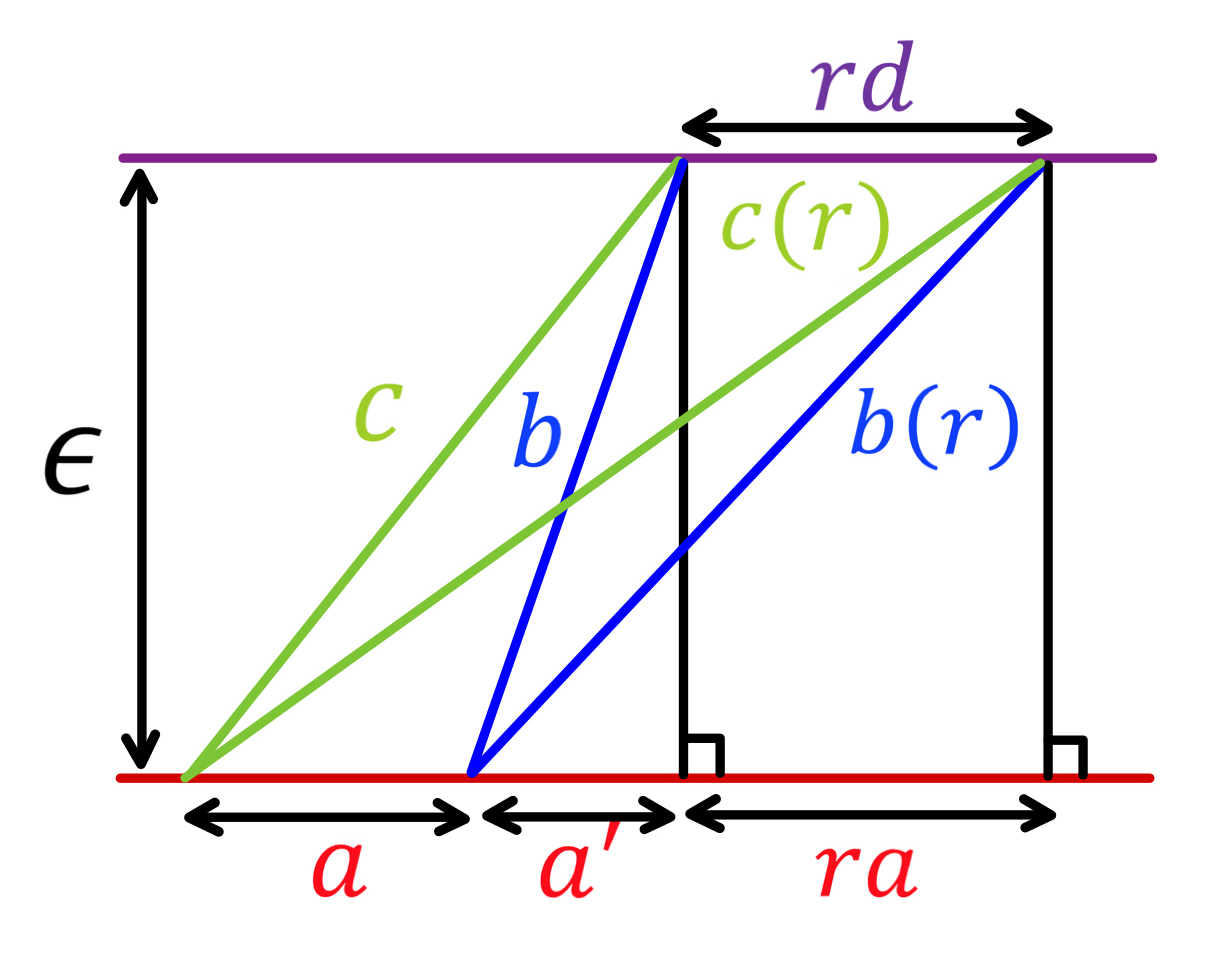}
	\end{minipage}
	\caption{The maximal collar neighbourhood for $\alpha$ on $\punc$ unfolded into a tube (left) and half of the neighbourhood unfolded in the universal cover $\H^2$ (right).}
	\label{fig:maxCollarCover}
\end{figure}

Consider the image of half of the maximal collar neighbourhood in the universal cover $\H^2$. A representation of this is shown in Figure \ref{fig:maxCollarCover}. This can be used to compute how the lengths $a, b, c$ change after taking a fractional Dehn twist about $\alpha$ of some magnitude $t \in \R$, $a(t), b(t), c(t)$.

\begin{theorem}
\label{res:earthquakeLength}
	For a fixed starting point $\mathbf{w}=(a,b,c)\in\compL\subseteq\RPlusSpace$, the earthquake deformation about $\alpha$ in triangle lengths is given by the map $\earthquakeHyp(\mathbf{w}):\R \to \RPlusSpace$,	
	\[
	\begin{split}
		\earthquakeHyp(\mathbf{w}): \R 	&\to 		\RPlusSpace\\
		r 							&\mapsto 	\begin{pmatrix} a \\ b_\alpha(r) \\ b_\alpha(r-1)) \end{pmatrix}
	\end{split}
	\]
	where
	\begin{equation*}
	\begin{split}
		b_\alpha: \R 	&\to 		\RPlus \\
		r 			&\mapsto 	\cosh^{-1}\left(\cosh(b)\cosh(ra)-\left(\cosh(c)-\cosh(a)\cosh(b)\right)\frac{\sinh(ra)}{\sinh(a)}\right).
	\end{split}
	\end{equation*}
	Moreover, $\earthquakeHyp(\v)(n)=\bar{T}_\alpha^n(\v)$, $\forall n \in \Z.$
\end{theorem}

\begin{proof}
	The arrangement in Figure \ref{fig:maxCollarCover} provides a series of four hyperbolic right-angled triangles. Two triangles represent the initial relative formation of the curves $\alpha$, $\beta$, and $\alpha\beta$ with hypotenuse lengths $b$, $c$. Two triangles represent the formation of the curves $\alpha$, $\beta$, and $\alpha\beta$ after taking a fractional Dehn twist of value $t$, with hypotenuse lengths $b(t)$, $c(t)$. 
	Given a hyperbolic right-angled triangle with side lengths $p$, $q$ and hypotenuse length $h$, the lengths are related by 
	\[
		\cosh(h)=\cosh(p)\cosh(q).
	\]
	
Using this relation and those provided in Theorem \ref{res:Buser}, the series of hyperbolic triangles can be solved to find $b(t)$ and $c(t)$ as functions in $t$ and the starting lengths $a$, $b$, $c$. Applying the collar equation from Equation \ref{eqn:collar} gives the formula for the earthquake deformation about $\alpha$.
\end{proof}

\begin{figure}[bht]
	\centering
	\includegraphics[width=0.5\textwidth]{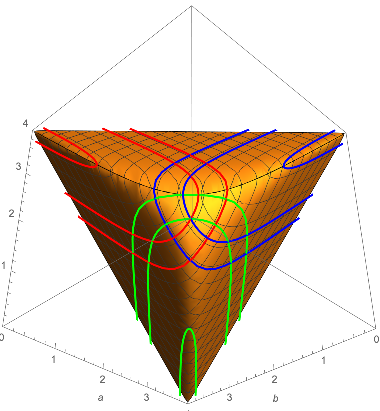}
	\caption{Example forward and backward earthquake deformations about $\alpha$ (red), $\beta$ (blue), and $\alpha\beta$ (green) on the image of $\comp$ in $\RPlusSpace$ using triangle lengths. Starting points are given by applying $\nu$ from Equation \ref{eqn:coordChange} to the respective sets $\mathcal{S}_\alpha$, $\mathcal{S}_\beta$, and $\mathcal{S}_{\alpha\beta}$ from Equation \ref{eqn:startPoint}.}
	\label{fig:orbitsHyp}
\end{figure}

A similar expression can be derived for the earthquake deformations about $\beta$ and $\alpha\beta$, $E_\beta^{\ell}$ and $E_{\alpha\beta}^{\ell}$ respectively, in triangle lengths. This is done either by following the same derivation process or by applying the symmetry relationships from Remark \ref{rem:symAuto} to $\earthquakeHyp$, adjusting for the change of coordinates. See Figure \ref{fig:orbitsHyp} for some examples.

Compare the formula for $\earthquakeHyp$ as given in Theorem \ref{res:earthquakeLength} to $\flowTrace$ from Theorem \ref{res:earthquakeTrace}. There is a natural homeomorphism between the two coordinate systems,
\begin{equation}
\begin{split}
\label{eqn:coordChange}
	\nu: \comp 						&\to \compL \\
	\begin{pmatrix} x \\ y \\ z \end{pmatrix} 	&\mapsto  \begin{pmatrix} \cosh^{-1}\left(\frac{x}{2}\right) \\ \cosh^{-1}\left(\frac{y}{2}\right) \\ \cosh^{-1}\left(\frac{z}{2}\right) \end{pmatrix}
\end{split}
\end{equation}
using Equation \ref{eqn:lengthTrace}. We can use this to relate $\earthquakeHyp$ and $\flowTrace$.

\begin{theorem}
	Take starting point $\w=(a,b,c) \in \compL\subseteq \RPlusSpace$, then
	\[
		\earthquakeHyp(\w)=\nu\left(\flowTrace\left(\nu^{-1}\left(\w\right)\right)\right).
	\]
\end{theorem}

\begin{proof}
	The proof uses the relationship between the trace coordinates and triangle lengths given in Equation \ref{eqn:lengthCoord}, the fact that we are operating on the level set $\kappa(x,y,z)=-2$ (corresponding to component $\comp$), and that 
	\[
		\left(\cosh(\theta)\pm\sinh(\theta)\right)^n=\cosh(n\theta)\pm\sinh(n\theta).
	\]
	Note the special case $x=2$ in Theorem \ref{res:earthquakeTrace} does not lie on Teichm\"uller space. This corresponds to $a=0$.
\end{proof}

\begin{corollary}
\label{cor:flowEarthquakeTrace}
	The earthquake deformation about $\alpha$ in trace coordinates is precisely the continuous flow extension of the Dehn twist about $\alpha$ on the component $\comp$, 
	\[
		\earthquakeTrace(\v)=\flowTrace(\v)\text{, } \forall \v\in\comp.
	\]
\end{corollary}

\begin{remark}
\label{rem:reparametriseHyp}
	The parameter $r\in\R$ represents the relative length twisted along the curve $\alpha$, with $r=1$ indicating a twist of magnitude the full length of $\alpha$. It can be useful to express the formula for earthquakes in terms of the objective length of the curve. This is often used as the more common convention when studying earthquakes.

	In particular, when using the counting measure, the length is precisely the measure of the curve. For the earthquake about $\alpha$, the coordinate change is $s=r\ell(\alpha)=2ra$.

	After applying $s=2ra$ the formula from Theorem \ref{res:earthquakeLength} becomes,
	\[
	\begin{split}
		\earthquakeHyp(\w): \R 	&\to 		\RPlusSpace\\
		s 							&\mapsto 	\begin{pmatrix} a \\ b_\alpha(s) \\ b_\alpha(s-2a)) \end{pmatrix}
	\end{split}
	\]
	where
	\begin{equation*}
	\begin{split}
		b_\alpha: \R 	&\to 		\RPlus \\
		s 			&\mapsto 	\cosh^{-1}\left(\cosh(b)\cosh(s/2)-\left(\cosh(c)-\cosh(a)\cosh(b)\right)\frac{\sinh(s/2)}{\sinh(a)}\right).
	\end{split}
	\end{equation*}
	
A similar reparametrisation can be done for the result in Theorem \ref{res:earthquakeTrace}.
\end{remark}

The results presented in this section have extended the orbits of the Dehn twists about $\alpha$, $\beta$, and $\alpha\beta$ to get an explicit formula for the earthquake deformation in triangle length coordinates. We have proven this is equivalent to the flow from Section \ref{sec:recur}, which extended the orbits of the Dehn twists about $\alpha$, $\beta$, and $\alpha\beta$ to get a continuous flow in trace coordinates. Thus, this continuous flow is the earthquake deformation in trace coordinates. The formula in trace coordinates gives an algebraic interpretation of earthquake deformations; the formula in triangle lengths gives a geometric interpretation of earthquake deformations.

%% --- ALGORITHM ---------------------------------------------------------------------------------------------------------------------------------------------------

\subsection{Change of coordinates algorithm}
\label{sec:algorithm}

We present a change of coordinates algorithm that outlines how to find similar earthquake deformation formula for any simple closed curve. Given the framing $\alpha$ and $\beta$, Sections \ref{sec:recur} and \ref{sec:hypGeo} provide a description the earthquake deformation about $\alpha$, $\beta$, $\alpha\beta$ in trace coordinates and triangle lengths, respectively. The algorithm uses these results to construct the earthquake deformation for any simple closed geodesic $\gamma$ in the chosen framing. For simplicity we consider only trace coordinates in the description of the algorithm. The definition of triangle lengths in Equation \ref{eqn:lengthCoord} can be used to convert the expressions.

The algorithm is broken into three steps. Given an input of a simple closed geodesic $\gamma$ as a word in $\alpha$ and $\beta$, $\gamma = w_\gamma(\alpha,\beta)$,
\begin{enumerate}
	\item Compute a simple closed curve $\delta$ such that $i(\gamma,\delta)=1$;
	\item Calculate the earthquake deformation for $\gamma$ in the framing $\gamma$ and $\delta$ using Theorem \ref{res:earthquakeTrace};
	\item Transform the expression for the earthquake deformation expression to the original framing $\alpha$ and $\beta$.
\end{enumerate}

We elaborate on each of these steps in order.

\subsubsection{Step (1)}
First compute a simple closed curve $\delta$ such that $i(\gamma,\delta)=1$. 

Such combinations of curves $\gamma$, $\delta$ can be found via a mapping class applied to $\alpha$, $\beta$. All mapping classes on $\punc$ are represented by an element 
\begin{equation*}
	M = \begin{pmatrix} m_1 & m_2 \\ n_1 & n_2 \end{pmatrix} \in \SLZ.
\end{equation*}

Consider $\gamma$ an $(m_1,n_1)$-type curve for some $m_1,n_1 \in \Z$ with $(m_1,n_1) \neq (0,0)$. Then to 
prove existence of $\delta$ we need to find an $(m_2,n_2)$-type curve with $m_2,n_2 \in\Z$ that satisfies $m_1n_2-m_2n_1=1$. As $(m_1,n_1)$ is non-trivial, such a vector will always exist. Then there is a mapping class $M$ that takes $\alpha$ to $\gamma$ and $\beta$ to $\delta$ with $i(\gamma,\delta)=1$. We can write $\gamma$ and $\delta$ as words in $\alpha$, $\beta$,
\begin{equation*}
	\gamma = w_\gamma(\alpha,\beta), \delta = w_\delta(\alpha,\beta).
\end{equation*}

\subsubsection{Step (2)}
Second, calculate the earthquake deformation for $\gamma$ in the framing given by $\gamma$ and $\delta$. 

Define new ``local" trace coordinates 
\begin{equation}
\label{eqn:localTrace}
	\begin{pmatrix} x' \\ y' \\ z' \end{pmatrix} = \begin{pmatrix} \tr(\rho(\gamma)) \\ \tr(\rho(\delta)) \\ \tr(\rho(\gamma\delta)) \end{pmatrix}.
\end{equation}

Substituting $(x',y',z')$ for $(x,y,z)$ in Theorem \ref{res:earthquakeTrace} and using Corollary \ref{cor:flowEarthquakeTrace} gives results for the earthquake deformation about $\gamma$ in local trace coordinates starting at some point $\v'=(x',y',z') \in \comp$. 
\begin{equation}
\label{eqn:algEarthquake}
\begin{split}
	E_{\gamma}^{'\text{tr}}(\v')=\earthquakeTrace(\v'): \R 	&\to 		\comp\\
	r 											&\mapsto 	\begin{pmatrix} x' \\ y_{\gamma}(r) \\ y_{\gamma}(r-1) \end{pmatrix}
\end{split}
\end{equation}
where
\begin{equation*}
\begin{split}
	y_\gamma: \R 	&\to 		\R \\
	r 			&\mapsto 	
	\begin{cases} 
		C_+^{\gamma}(\v')S_+^{\gamma}(\v')^r+C_-^{\gamma}(\v')S_-^{\gamma}(\v')^r 	\hspace{0.2cm}		\textnormal{ if } x' \neq 2 \\
		y' +(y'-z')r \hspace{3.03cm}														\textnormal{ if } x' = 2  
	\end{cases}
\end{split}
\end{equation*}
with 
\begin{equation*}
\begin{split}
	C_\pm^\gamma(\v') &= \frac{y'(x'^2-4)\pm(x'y'-2z')\sqrt{x'^2-4}}{2(x'^2-4)}, \\
	S_\pm^\gamma(\v') &= \frac{x'\pm\sqrt{x'^2-4}}{2}.
\end{split}
\end{equation*}

\subsubsection{Step (3)}
Third, transform the expression for the earthquake deformation for $\gamma$ in the framing $\gamma$ and $\delta$ back to the original framing given by $\alpha$ and $\beta$. We do this by relating the local trace coordinates to the original trace coordinates. 

We can determine a change of coordinates map $\phi$ from original to local trace coordinates by expanding $\gamma$ and $\delta$ as words in $\alpha$ and $\beta$, 
\[
	\begin{split}
		\phi: \comp 						&\to 		\comp\\
		\begin{pmatrix} x \\ y \\ z \end{pmatrix} 	&\mapsto \begin{pmatrix} x' \\ y' \\ z' \end{pmatrix} = \begin{pmatrix} \tr(\rho(w_\gamma(\alpha,\beta))) \\ \tr(\rho(w_\delta(\alpha,\beta))) \\ \tr(\rho(w_\gamma(\alpha,\beta)w_\delta(\alpha,\beta))) \end{pmatrix}.
	\end{split}
\]
The map $\phi$ is the polynomial automorphism equivalent to the mapping class $M$ from Step (1). Use trace identities to express the right-hand side in terms of $x=\tr(\rho(\alpha))$, $y=\tr(\rho(\beta))$, and $z=\tr(\rho(\alpha\beta))$. 

Similarly, as $\gamma$ and $\delta$ provide a framing of $\punc$, we can write $\alpha$ and $\beta$ as words in $\gamma$, $\delta$,
\begin{equation*}
	\alpha = w_\alpha(\gamma,\delta), \beta = w_\beta(\gamma,\delta).
\end{equation*}
This can be constructed from the inverse of the mapping class $M$.

Then we can determine an inverse change of coordinates map $\psi$ from local to original trace coordinates by expanding $\alpha$ and $\beta$ as words in $\gamma$ and $\delta$,
\[
\begin{split}
	\psi: \comp						&\to 		\comp \\
	\begin{pmatrix} x' \\ y' \\ z' \end{pmatrix} 	&\mapsto \begin{pmatrix} x \\ y \\ z \end{pmatrix} 
	= \begin{pmatrix} \tr(\rho(w_\alpha(\gamma,\delta))) \\ \tr(\rho(w_\beta(\gamma,\delta))) \\ \tr(\rho(w_\alpha(\gamma,\delta)w_\beta(\gamma,\delta))) \end{pmatrix}.
\end{split}
\]
The map $\psi$ is the polynomial automorphism equivalent to the inverse of the mapping class $M$ from Step (1). Use trace identities to expand the right-hand side in terms of $x'=\tr(\rho(\gamma))$, $y'=\tr(\rho(\delta))$, $z'=\tr(\rho(\gamma\delta))$. The maps $\phi$ and $\psi$ are inverses. 

Then use $\phi$, $\psi$, and $\earthquakeTrace$ to get the earthquake deformation about $\gamma$ in original trace coordinates starting at $\v = (x,y,z) \in \comp$,
\[
	E_\gamma^{\tr}(\v) = \psi \left(\earthquakeTrace(\phi(\v))\right).
\]

The solution for examples of $\gamma$ is given in Section \ref{sec:examples}.  

%% --- FENCHEL-NIELSEN LENGTH TWIST COORDINATES ---------------------------------------------------------------------------------------------
\subsection{Studying deformations using Fenchel-Nielsen length-twist coordinates}
\label{sec:FN}

It is useful to study the earthquakes in other coordinate systems. Of particular interest is Fenchel-Nielsen length-twist coordinates. 

Fenchel-Nielsen coordinates are based on a pants decomposition of the surface. See, for example, \cite{FarbMargalit12}. Consider the pants decomposition given by $\alpha$ on $\punc$. The Fenchel-Nielsen coordinates using $\alpha$, $(\ella,\taua)$, are related to trace coordinates and triangle lengths by
\begin{equation}
\label{eqn:FNdefn}
	\begin{pmatrix} \ella \\ \taua \end{pmatrix} = \begin{pmatrix} 2\cosh^{-1}\left(\frac{x}{2}\right) \\ 2\cosh^{-1}\left(\frac{y}{2\coth\left(\cosh^{-1}\left(\frac{x}{2}\right)\right)}\right) \end{pmatrix} = \begin{pmatrix} 2a \\ 2\cosh^{-1}\left(\frac{\cosh(b)}{\coth\left(a\right)}\right) \end{pmatrix}.
\end{equation}
We could write other coordinate changes for Fenchel-Nielsen coordinates using any other simple closed curve $\gamma$, $(\ell_\gamma,\tau_\gamma)$. Label the space $\compF = \RR$. Each choice of pants decomposition defines an approach to identify $\compF \cong \T(\punc)$.

Label the coordinate change from trace coordinates
\begin{equation}
\begin{split}
\label{eqn:FNcoordChange}
	\zeta: \comp&\to\compF \\
	\begin{pmatrix} x \\ y \\ z \end{pmatrix} &\mapsto \begin{pmatrix} 2\cosh^{-1}\left(\frac{x}{2}\right) \\ 2\cosh^{-1}\left(\frac{y}{2\coth\left(\cosh^{-1}\left(\frac{x}{2}\right)\right)}\right) \end{pmatrix}.
\end{split}
\end{equation}

These coordinates are no longer symmetric like the trace coordinates or triangle lengths, but they are projected to the plane in a meaningful way. For a fixed starting point $\u=(\ell,\tau)\in\compF$, the earthquake about $\alpha$ in Fenchel-Nielsen coordinates for $\alpha$ is given by
\begin{equation}
\begin{split}
	\label{eqn:alphaFN}
	\earthquakeAFN(\u):& \, \R \to \RR  \\
	s &\mapsto \begin{pmatrix} \ell \\ \tau-s \end{pmatrix}.
\end{split}
\end{equation}
Note this earthquake is parametrised by hyperbolic length, as per Remark \ref{rem:reparametriseHyp}.

\begin{figure}[bht]
	\centering
	\includegraphics[width=0.475\textwidth]{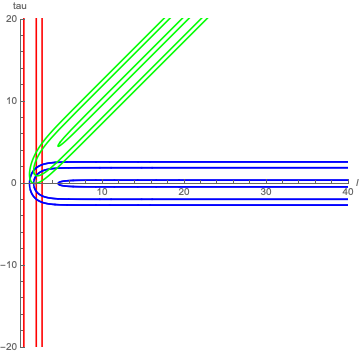}
	\caption{Example forward and backward earthquake deformations about $\alpha$ (red), $\beta$ (blue), and $\alpha\beta$ (green) using Fenchel-Nielsen length-twist coordinates. Starting points are given by applying $\zeta$ from Equation \ref{eqn:FNcoordChange} to the respective sets $\mathcal{S}_\alpha$, $\mathcal{S}_\beta$, and $\mathcal{S}_{\alpha\beta}$ from Equation \ref{eqn:startPoint}.}
	\label{fig:orbitsFN}
\end{figure}

The expressions for earthquake deformations about $\beta$ and $\alpha\beta$, $\earthquakeBFN(\u)$ and $\earthquakeCFN(\u)$ respectively, are increasingly complex. See Figure \ref{fig:orbitsFN} for examples. 

The wider set of measured geodesic laminations, $\ML$, of $\punc$ is naturally identified with $\HR \cong \R^2$. A line through the origin represents a lamination and a vector along this line represents an associated measured lamination. Taking the antipodal map represents the projective laminations. Label the set of earthquakes in Fenchel-Nielsen coordinates for Teichm\"uller space as $\mathcal{E}\subseteq\compF$.

For any fixed point $\u\in\compF$ study the map that takes each lamination to their associated earthquake in Fenchel-Nielsen coordinates for $\alpha$ starting at $\u$,
\begin{equation}
\begin{split}
	\label{eqn:earthquake}
	{E_{\u}}: \ML\subseteq \HR \cong \R^2 &\to \mathcal{E} \subseteq \compF \\
	\L &\mapsto E^\text{FN}_{\L}(\u),
\end{split}
\end{equation}
with 
\begin{equation*}
\begin{split}
	E^\text{FN}_{\L}(\u,\cdot): \R &\mapsto \compF\\
	s &\mapsto \begin{pmatrix} \ell_\alpha^\L(\u,s) \\ \tau_\alpha^\L(\u,s) \end{pmatrix}.
\end{split}
\end{equation*}
For any simple closed curve $\delta$, the functions $\ell_\delta^\L(\u,s)$ and $\tau_\delta^\L(\u,s)$ represent the length and twist parameter of $\delta$ at the point on the earthquake $E^\text{FN}_{\L}(\u,s)$. Notation for the measure of the lamination is omitted for brevity. 

\begin{figure}[bht]
	\centering
	\begin{minipage}{0.44\textwidth}
		\centering
		\includegraphics[width=0.9\textwidth]{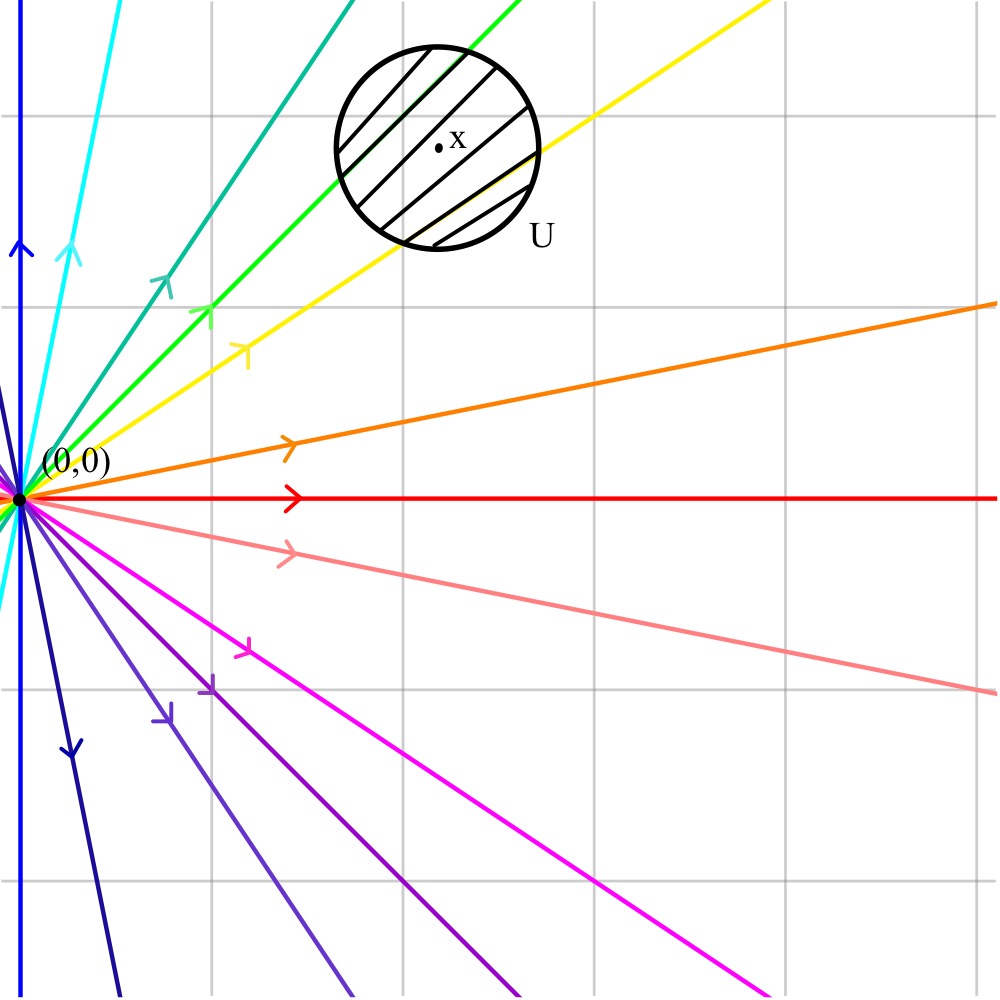}
	\end{minipage}
	\begin{minipage}{0.1\textwidth}
		\centering
		\center{$\mapsto$}
	\end{minipage}
	\begin{minipage}{0.44\textwidth}
		\centering
		\includegraphics[width=0.92\textwidth]{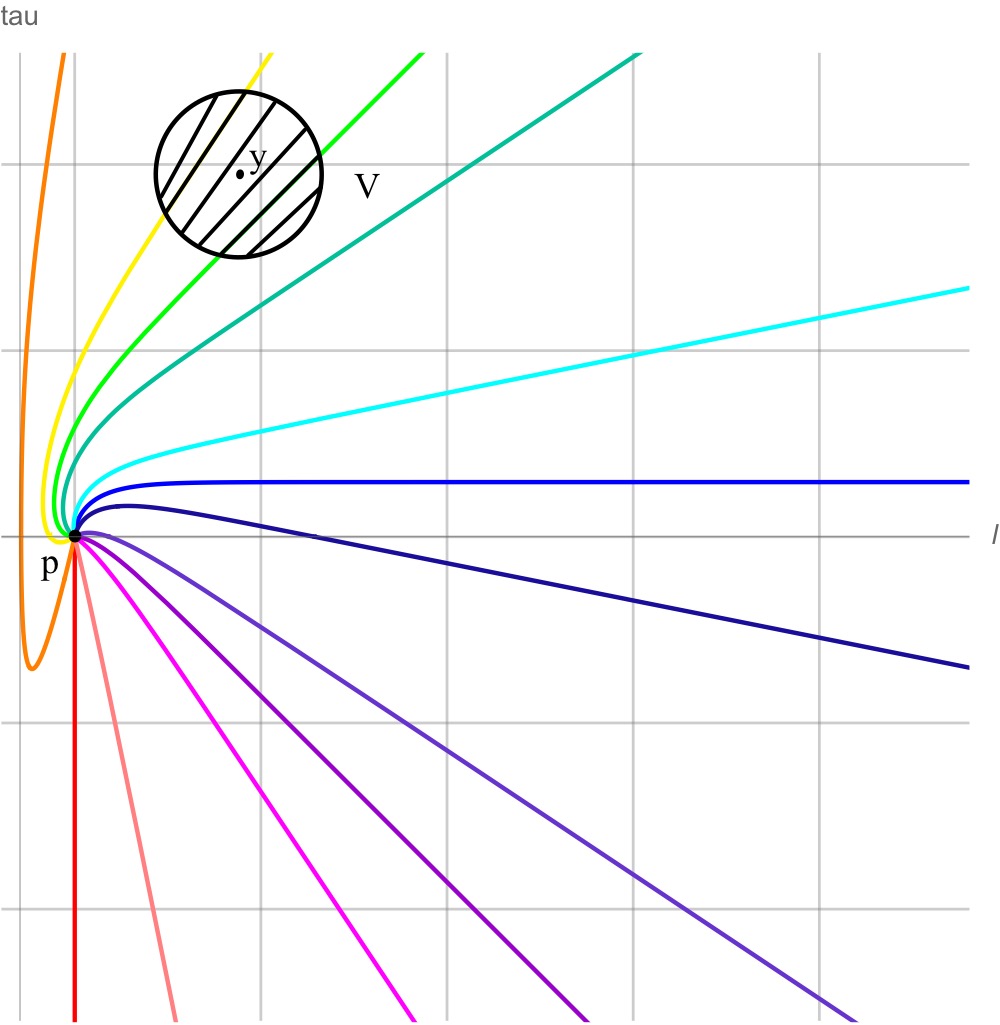}
	\end{minipage}
	
	\vspace{1cm}
	\rotatebox[origin=c]{315}{$\mapsto$}
	\hspace{8cm}
	\rotatebox[origin=c]{45}{$\mapsto$}
	
	\vspace{1cm}
	\centering
	\includegraphics[width=0.42\textwidth]{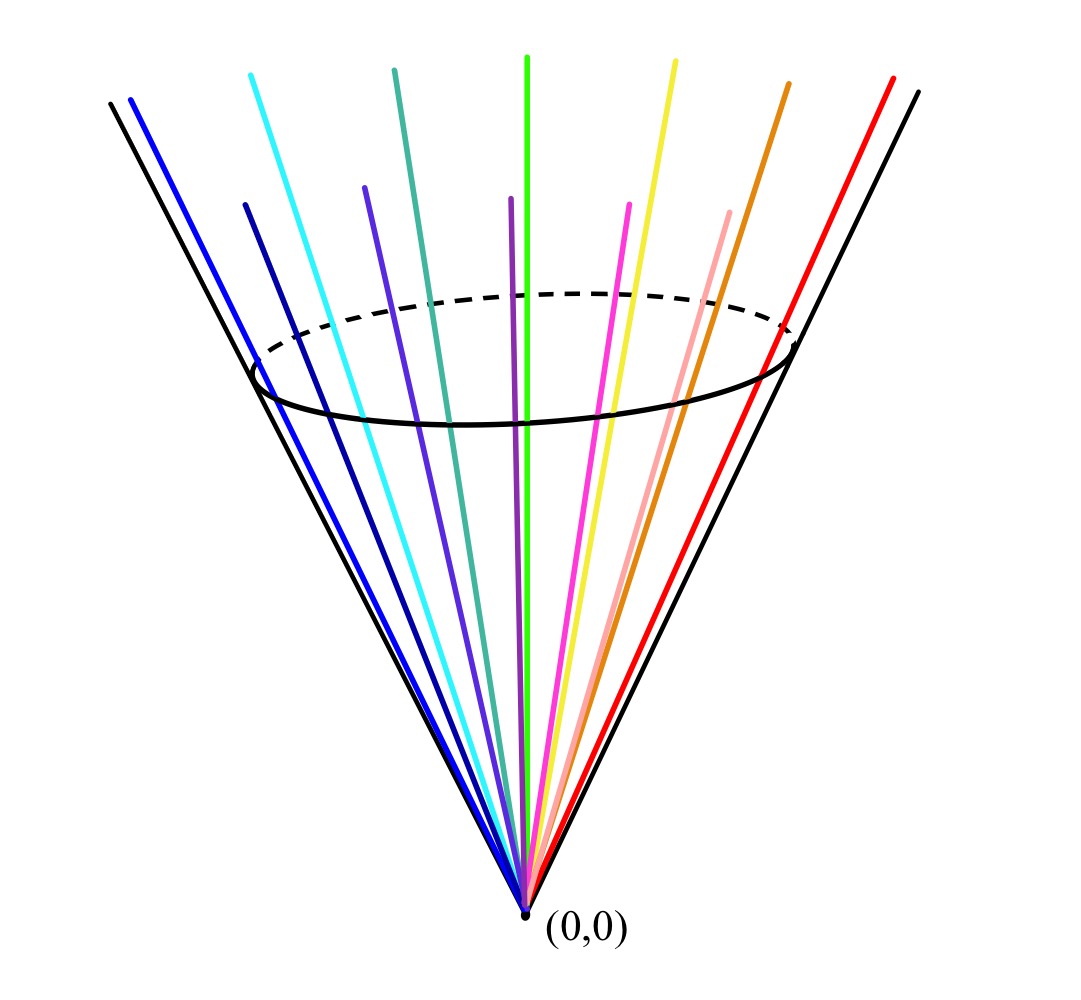}
	\caption{Illustration of map $E_\textbf{p}$ from Equation \ref{eqn:earthquake} taking laminations with basepoint $(0,0)$ in $\R^2$ (top left) to forward earthquakes with starting point $\textbf{p}$ in $\compF$ (top right) going via the space of projective laminations (bottom). Clockwise, example laminations are given by simple closed curves $\beta$ (blue), $\alpha\beta^5$, $\alpha^2\beta^3$, $\alpha\beta$, $\alpha^3\beta^2$, $\alpha^5\beta$, $\alpha$ (red), $\alpha^5\beta^{-1}$, $\alpha^3\beta^{-2}$, $\alpha\beta^{-1}$, $\alpha^2\beta^{-3}$, $\alpha\beta^{-5}$. This illustrates the behaviour from Remark \ref{rem:Emap}.}
	\label{fig:lamToQuake}
\end{figure}

\begin{remark}
\label{rem:Emap}
We make some statements about the map $E_{\textbf{p}}$.
\begin{enumerate}
	\item The map is continuous and bijective for each choice of $\textbf{p}$. Continuity is given through work of Kerckhoff (see Lemma 1.2, \cite{Ker83}) and bijectivity is given by this and Thurston's earthquake theorem. 
	\item We know $\ML$ foliates $\HR$. In particular, any neighbourhood $U$ of a point $\textbf{x}\in\R^2$ is foliated by a subset of $\ML$ and is mapped to a neighbourhood $V$ of a point $\textbf{y}\in\compF$ that is foliated by a subset of $\mathcal{E}.$ 
\end{enumerate}
	Recall the slope of a simple closed curve given in Equation \ref{eqn:slope}.
\begin{enumerate}[resume]
	\item Consider a sequence of distinct measured geodesic laminations ordered counter-clockwise by their slope,
		\[
			\L_1, \L_2, \ldots, \L_n, \hspace{0.8cm} \sl(\L_1)<\sl(\L_2)<\ldots<\sl(\L_n).
		\]
		The sequence of associated earthquakes will follow the same order clockwise.
\end{enumerate}
\end{remark}
The map is depicted in Figure \ref{fig:lamToQuake}. 

We can make further observations about the relationship between the slope of laminations and the slope of the associated earthquakes. 
\begin{theorem}
\label{thm:earthquakeAsymptotics}
	Take simple closed curve $\gamma$. The slope of the associated earthquake $E^\text{FN}_{\gamma}(\u,s)$ converges to the inverse of the slope of $\gamma$ for all starting points $\u\in\compF$. We have, 
	\[
	\begin{split}
		\frac{\tau_{\alpha}^\alpha(\u,s)}{\ell_{\alpha}^\alpha(\u,s)}\,&\to\,-\frac{1}{\sl(\alpha)}=-\infty \hspace{0.45cm} \text{ as }s\to\infty, \, \forall \u\in\compF\text{ and }\\
		\frac{\tau_\alpha^{\gamma}(\u,s)}{\ell_\alpha^{\gamma}(\u,s)}\,&\to\,\frac{1}{\sl(\gamma)}\hspace{0.3cm}=\frac{i(\gamma,\beta)}{i(\gamma,\alpha)} \text{ as }s\to\infty, \, \forall \u\in\compF.
	\end{split}
	\]
\end{theorem}

The following corollary is apparent after applying the density result in Theorem \ref{thm:dense}.
\begin{corollary}
	Take a measured geodesic lamination $\L$ that is not $\alpha$. The slope of the associated earthquake $E^\text{FN}_{\L}(\u)$ converges to the inverse of the slope of $\L$ for all starting points $\u\in\compF$. That is,
	\[
		\frac{\tau_\alpha^\L(\u,s)}{\ell_\alpha^\L(\u,s)}\,\to\,\frac{1}{\sl(\L)} \text{ as }s\to\infty, \, \forall \u\in\compF.
	\]
\end{corollary}

\begin{proof}[Proof of Theorem \ref{thm:earthquakeAsymptotics}]
	The theorem is proved in three main steps.
	\begin{enumerate}[label=(\roman*)]
		\item \label{pf:stepOne} Prove the result directly for $\alpha$ and $\beta$ using explicit formula for the earthquakes.
	\end{enumerate}
	Consider now a simple closed curve $\gamma$ that intersects both $\alpha$ and $\beta$ at least once.
	\begin{enumerate}[label=(\roman*),resume]
		\item \label{pf:stepTwo} Prove that 
			\[
				\frac{\ell_\beta^{\gamma}(\u,s)}{\ell_\alpha^{\gamma}(\u,s)}\,\to\,\left\lvert\frac{i(\gamma,\beta)}{i(\gamma,\alpha)}\right\rvert \text{ as }s\to\infty, \, \forall \u\in\compF.
			\]
		\item \label{pf:stepThree} Prove that
			\[
					\ell_\beta^\gamma(\u,s) \sim	\begin{cases} 
									\tau_\alpha^\gamma(\u,s) & \text{ if $\sl(\gamma)>0$}\\
									-\tau_\alpha^\gamma(\u,s) & \text{ if $\sl(\gamma)<0$}
								\end{cases}
					\text{ as }s\to\infty, \, \forall \u\in\compF.
			\]
	\end{enumerate}
	Bringing this all together proves the theorem. 
	
	Part \ref{pf:stepOne}.
	The result for $\alpha$ is clear from the formula from Equation \ref{eqn:alphaFN}. Fix $\u\in\compF$. We have
	\[
		\frac{\tau_\alpha^{\gamma}(\u,s)}{\ell_\alpha^{\gamma}(\u,s)}\,\to -\infty = -\frac{1}{\sl(\alpha)} \text{ as }s\to\infty.
	\]
	
	Using the change of coordinates formulas we find an expression for the earthquake about $\beta$. 
	\[
	\begin{split}
		\earthquakeBFN(\u,s)	&= \zeta(\SymE_\text{Rot}(\earthquakeTrace(\SymE_\text{Rot}^{-1}(\zeta^{-1}(\u)),s))\\
							&= \begin{pmatrix} \ell_\alpha^\beta(\u,s) \\ \tau_\alpha^\beta(\u,s) \end{pmatrix}
	\end{split}
	\]
	with 
	\[
	\begin{split}
		\f(\ell,\tau) &= \cosh^{-1}\left(\cosh\taufrac\coth\ellfrac\right),\\
		\cosh\left(\frac{\ell_\alpha^\beta(\u,s)}{2}\right) &= \cosh\ellfrac \left( \cosh\ssfrac + \frac{\sinh\ssfrac \sinh\taufrac}{\sinh\left(\f(\ell,\tau)\right)}\right),\\
		\cosh\left(\frac{\tau_\alpha^\beta(\u,s)}{2}\right) &= \cosh\taufrac \cosh\ellfrac \tanh\left(\frac{\ell_\beta(\u,s)}{2}\right).
	\end{split}
	\]
	Then the slope is
	\[
	\begin{split}
		\frac{\tau_\alpha^\beta(\u,s)}{\ell_\alpha^\beta(\u,s)}=\frac{2\cosh^{-1}\left(\cosh\taufrac \cosh\ellfrac \tanh\left(\frac{\ell_\beta(\u,s)}{2}\right)\right)}{\ell_\beta(\u,s)}.
	\end{split}
	\]
	We find 
	\[
		\ell_\alpha^\beta(\u,s) \to \infty \text{ and } \tanh\left(\ell_\beta^\alpha(\u,s)\right) \to 1 \text{ as } s\to\infty
	\]
	and
	\[
		\frac{\tau_\alpha^\beta(\u,s)}{\ell_\alpha^\beta(\u,s)} \to \frac{2\cosh^{-1}\left(\cosh\taufrac \cosh\ellfrac\right)}{\infty} = 0 = \frac{1}{\sl(\beta)} \text{ as } s\to\infty.
	\]
	
	Part \ref{pf:stepTwo}.
	We exploit the relationship between earthquakes and Dehn twists. Let $T_\gamma$ be the Dehn twist about $\gamma$ and $\delta$ another simple closed curve. Consider the length of $\ell_{T_\gamma^n(\delta)}(\u)$. 
	We claim
	\begin{equation}
	\label{eqn:claim}
		\ell_{T_\gamma^n(\delta)}(\u) \sim n \lvert i(\gamma,\delta) \rvert \ell_\gamma(\u) \text{ as }n\to\infty 
	\end{equation}
	If $\lvert i(\gamma,\delta) \rvert=k$, we can split $\delta$ into $k$ segments, $\epsilon_i$, $1\leq i \leq k$ such that  
	\[
	\begin{split}
		\delta &=\epsilon_1\epsilon_2\ldots\epsilon_k, \\
		\lvert i(\gamma,\epsilon_i)\rvert &=1 \hspace{0.25cm} \forall 1\leq i \leq k,\text{ and }\\
		T_\gamma^n(\delta) &=\gamma^n\epsilon_1\gamma^n\epsilon_2\ldots\gamma^n\epsilon_k.
	\end{split}
	\]
	Note that the second condition implies that the pair $\gamma,\epsilon_i$ generate $\pi_1(\punc)$ for all $1\leq i \leq k$. In particular, that $\tr(\rho(\gamma\epsilon_i\gamma^{-1}\epsilon_i^{-1})=-2$ for all $1\leq i \leq k$ and for any representation $\rho$ that represents a hyperbolic structure in $\T(\punc)$.
	
	Take a such a representation $\rho$ such that,
	\[
		\rho(\gamma) = \begin{pmatrix} \lambda & 0 \\ 0 & \lambda^{-1} \end{pmatrix} \text{ with } \lambda>0,\hspace{1cm} \rho(\epsilon_i) = \begin{pmatrix} a_i & b_i \\ c_i & d_i \end{pmatrix}.
	\]
	Note, if $a_i=0$, $\rho(\epsilon_i)$ sends the attracting fixed point of $\rho(\gamma)$ to the repelling fixed point and their invariant axes are distinct. This contradicts work by Goldman (Lemma 3.4.5, \cite{Goldman03}). Thus, $a_i\neq0$ for all $1\leq i \leq k$.
	
	Then 
	\[
	\begin{split}
		\ell_{T_\gamma^n(\delta)}(\u) &=2\cosh^{-1}\left(\frac{\tr\left(\rho\left(\gamma^n\epsilon_1\gamma^n\epsilon_2\ldots\gamma^n\epsilon_k\right)\right)}{2}\right) \\
		& = 2\cosh^{-1}\left(\frac{\lambda^{kn}a_1a_2\ldots a_k+\ldots+\lambda^{-kn}d_1d_2\ldots d_k}{2}\right) \\
		& \sim 2\cosh^{-1}\left(\frac{\lambda^{kn}a_1a_2\ldots a_k}{2}\right) \text{ as }n\to\infty \\
		& \hspace{0.5cm} = 2 \log\left(\frac{\lambda^{kn}a_1a_2\ldots a_k}{2}+\sqrt{\left(\frac{\lambda^{kn}a_1a_2\ldots a_k}{2}\right)^2-1}\right)\\
		& \hspace{0.5cm} \sim 2 \log\left(\lambda^{kn}a_1a_2\ldots a_k\right) \text{ as }n\to\infty \\ 
		& \hspace{1cm} = nk 2 \log(\lambda) + \log(a_1a_2\ldots a_k) \\
		& \hspace{1cm} = n\lvert i(\gamma,\delta) \rvert \ell_\gamma(\u) + \log(a_1a_2\ldots a_k) \\
		& \hspace{1cm} \sim n \lvert i(\gamma,\delta) \rvert \ell_\gamma(\u) \text{ as }n\to\infty. 
	\end{split}
	\]
	The same asymptotic expression can be deduced taking the inverse Dehn twist,
	\[
		\ell_{T_\gamma^{-n}(\delta)}(\u) \sim n \lvert i(\gamma,\delta) \rvert \ell_\gamma(\u) \text{ as }n\to\infty.
	\]
	
	We get the following consequences to Equation \ref{eqn:claim},
	\[
		\ell_{\delta}\left(T_\gamma^n(\u)\right)=\ell_{T_\gamma^{-n}(\delta)}(\u) \sim n \lvert i(\gamma,\delta) \rvert \ell_\gamma(\u) \text{ as }n\to\infty,
	\]
	Then,
	\[
		\ell_{\delta}\left(E_\gamma^{\text{FN}}(\u,2\pi n)\right) = \ell_{\delta}\left(T_\gamma^n(\u)\right) \sim n \lvert i(\gamma,\delta) \rvert \ell_\gamma(\u) \text{ as }n\to\infty.
	\]
	The length function is convex on Teichm\"uller space (see \cite{Ker85}) and we can deduce
	\[
		\ell_{\delta}^\L\left(\u,s\right)=\ell_{\delta}\left(E_\gamma^{\text{FN}}(\u,s)\right)  \sim \frac{s}{2\pi} \lvert i(\gamma,\delta) \rvert \ell_\gamma(\u) \text{ as }s\to\infty.
	\]
	Then compare the ratio
	\[
	\begin{split}
		\frac{\ell_{\beta}^\L\left(\u,s\right)}{\ell_{\alpha}^\L\left(\u,s\right)}	&\sim \frac{\frac{s}{2\pi} \lvert i(\gamma,\beta) \rvert \ell_\gamma(\u)}{\frac{s}{2\pi} \lvert i(\gamma,\alpha) \rvert \ell_\gamma(\u)} \text{ as }s\to\infty\\
		&\to \left\lvert \frac{i(\gamma,\beta)}{i(\gamma,\alpha)} \right\rvert \text{ as }s\to\infty.
	\end{split}
	\]
	This concludes the proof of part \ref{pf:stepTwo}.
	
	Part \ref{pf:stepThree}.
	Consider the following change of pants decomposition formula from \cite{Papa12,Okai93} in the case of a cusp,
	\[
		\cosh\left(\frac{\ell_\beta^\gamma(\u,s)}{2}\right) = \sqrt{\frac{\cosh(\ell_\alpha^\gamma(\u,s))+1}{2\sinh\left(\frac{\ell_\alpha^\gamma(\u,s)}{2}\right)^2}}\cosh\left(\frac{\tau_\alpha^\gamma(\u,s)}{2}\right). 
	\]
	We have $\ell_\alpha^\gamma(\u,s)\to\infty$ as $s\to\infty$. Then, 
	\[
	\begin{split}
		\frac{\cosh(\ell_\alpha^\gamma(\u,s))+1}{2\sinh\left(\frac{\ell_\alpha^\gamma(\u,s)}{2}\right)^2} & =\frac{\cosh\left(\frac{\ell_\alpha^\gamma(\u,s)}{2}\right)^2+\sinh\left(\frac{\ell_\alpha^\gamma(\u,s)}{2}\right)^2+1}{2\sinh\left(\frac{\ell_\alpha^\gamma(\u,s)}{2}\right)^2}\\	
		&=\frac{1}{2}\left(1+\coth\left(\frac{\ell_\alpha^\gamma(\u,s)}{2}\right)^2+\csch\left(\frac{\ell_\alpha^\gamma(\u,s)}{2}\right)^2 \right)\\
		&\to 1 \text{ as } s\to\infty.
	\end{split}
	\]
	Using this gives 
	\[
	\begin{split}
		\cosh\left(\frac{\ell_\beta^\gamma(\u,s)}{2}\right) &\sim \cosh\left(\frac{\tau_\alpha^\gamma(\u,s)}{2}\right)\text{ as } s\to\infty \text{ and }\\
		\ell_\beta^\gamma(\u,s) &\sim \pm \tau_\alpha^\gamma(\u,s) \text{ as } s\to\infty.
	\end{split}
	\]
	Where the sign is determined by if $\gamma$ intersects $\alpha$ positively or negatively and $\beta$ positively or negatively . This is equivalent to,
	\[
	\ell_\beta^\gamma(\u,s) \sim	\begin{cases} 
									\tau_\alpha^\gamma(\u,s) & \text{ if $\sl(\gamma)>0$},\\
									-\tau_\alpha^\gamma(\u,s) & \text{ if $\sl(\gamma)<0$}.
								\end{cases}
	\]
	This finishes the proof.
\end{proof}

\begin{remark}
The result about slope does not rely on the orientation of the earthquake. That is, for each lamination $\L$, forwards and backwards earthquakes will tend to the same slope, 
	\[
	\begin{split}
		\frac{\tau_\alpha^\L(\u,-s)}{\ell_\alpha^\L(\u,-s)}\,\to\,\frac{1}{\sl(\L)} \text{ as }s\to\infty, \, \forall \u\in\compF.
	\end{split}
	\]
This is clear from the examples in Figure \ref{fig:orbitsFN} and is proven using the same approach as Theorem \ref{thm:earthquakeAsymptotics}. The only difference arises for $\L=\alpha$, where
	\[
		\frac{\tau_{\alpha}^\alpha(\u,-s)}{\ell_{\alpha}^\alpha(\u,-s)}\, \to\,\frac{1}{\sl(\alpha)}=\infty \hspace{0.45cm} \text{ as }s\to\infty, \, \forall \u\in\compF.
	\]
\end{remark}

We make a final remark regarding the relationship between forwards and backwards earthquakes.
\begin{remark}
It appears that, for a fixed starting point $\u\in\compF$, any backwards earthquake $E_\L^\text{FN}(\u,-s)$ about a lamination $\L$, will intersect all other forwards earthquakes $E_\mathcal{K}^\text{FN}(\u,s)$ for all laminations $\mathcal{K}$ in exactly one point (excluding $\u$). 

That is, for each $\L$ and each $\mathcal{K}$, there exists $s_{\L,\mathcal{K}}\neq0$ such that
	\[
		E_\L^\text{FN}(\u,-s_{\L,\mathcal{K}})=E_\mathcal{K}^\text{FN}(\u,s_{\L,\mathcal{K}}).
	\]
\end{remark}

The results from this section help to study general earthquake behaviour, as we will see in examples in the next section. 

%% --- EXAMPLES ---------------------------------------------------------------------------------------------------------------------------------------------------

\subsection{Examples}
\label{sec:examples}

We provide examples of the change of coordinates algorithm for two families of curves
\begin{enumerate}
	\item $\gamma_n = \alpha\beta\alpha^{n-1}$; and 
	\item $\gamma_n = T^n(\alpha)$ for $T=T_\beta T_\alpha^{-1}$.
\end{enumerate}

The first family of curves is an example of a sequence of simple closed curves that converges to a simple closed curve as $n\to\infty$. The second family of curves is an example of a sequence of simple closed curves that converges to a measured geodesic lamination that is not a simple closed curve as $n\to\infty$. 

For both examples we first describe the simple closed curve $\delta$ and the maps $\phi$ and $\psi$ required for the algorithm in both trace and triangle length coordinates. To this end, define new ``local" triangle length coordinates using the local trace coordinates from Equation \ref{eqn:localTrace},
\[
	\begin{pmatrix} a' \\ b' \\ c' \end{pmatrix} = \begin{pmatrix} \cosh^{-1}\left(\frac{x'}{2}\right) \\ \cosh^{-1}\left(\frac{y'}{2}\right) \\ \cosh^{-1}\left(\frac{z'}{2}\right) \end{pmatrix}.
\]
We then examine behaviour at the limit as $n\to\infty$ and investigate other ways to study the earthquakes.  

\begin{example}[$\gamma_n = \alpha\beta\alpha^{n-1}$]
\label{ex:one}

First, compute a sequence of simple closed curves $\delta_n$ such that $i(\gamma_n,\delta_n)=1$ for each $n\in\N$. One such example is $\delta_n = \alpha^{-1}$. The mapping class that takes $\alpha$, $\beta$ to $\gamma_n$, $\delta_n$ is $T_\alpha^{-n+1}T_\beta T_\alpha$. Given $\delta_n$ has no reliance on $n$ for this family of curves we exclude the subscript. Then
\[
	\begin{split}
		\gamma_n &= \alpha\beta\alpha^{n-1},\hspace{.2cm}		\delta 	= \alpha^{-1},\\
		\alpha 	&= \delta^{-1}, \hspace{0.8cm}	\beta 	= \delta\gamma_n\delta^{n-1}.
	\end{split}
\]

Second, calculate the earthquake deformation for $\gamma_n$ in the framing $\gamma_n$ and $\delta$ as given in Equation \ref{eqn:algEarthquake}.
	
Third, construct the series of maps $\phi_n$ and $\psi_n$ to go between the original framing $\alpha$, $\beta$ and the new framing $\gamma_n$, $\delta$ for each $n$. 

Going to local trace coordinates, 
\[
	\begin{split}
		\phi_n: \comp 						&\to 		\comp \\
		\begin{pmatrix} x \\ y \\ z \end{pmatrix} 	&\mapsto \begin{pmatrix} x'(n) \\ y'(n) \\ z'(n) \end{pmatrix} 
		= \begin{pmatrix} \tr(\rho(\alpha\beta\alpha^{n-1})) \\ \tr(\rho(\alpha^{-1})) \\ \tr(\rho(\alpha\beta\alpha^{n-2})) \end{pmatrix} 
		= \begin{pmatrix} x'(n) \\ x \\ x'(n-1) \end{pmatrix},
	\end{split}
\]
where $x'(n)=xx'(n-1)-x'(n-2)$, with $x'(0)=y$, $x'(1)=z$. The characteristic equation method for solving linear recurrence relations can be used to solve for $x'(n)$.

Going to original trace coordinates,
\[
	\begin{split}
		\psi_n: \comp 						&\to 		\comp \\
		\begin{pmatrix} x' \\ y' \\ z' \end{pmatrix} 	&\mapsto \begin{pmatrix} x(n) \\ y(n) \\ z(n) \end{pmatrix} 
		= \begin{pmatrix} \tr(\rho(\delta^{-1})) \\ \tr(\rho(\delta\gamma_n\delta^{n-1})) \\ \tr(\rho(\gamma_n\delta^{n-1})) \end{pmatrix} 
		= \begin{pmatrix} y' \\ y(n) \\ y(n-1) \end{pmatrix},
	\end{split}
\]
where $y(n)=y'y(n-1)-y(n-2)$, with $y(0)=x'$, $y(1)=z'$. Again, the characteristic equation method for solving linear recurrence relations can be used to solve for $y(n)$.

Then given some starting point $\v = (x,y,z) \in \comp$ the earthquake about the family $\gamma_n$ is given by the map $E_{n}^{\tr}(\v):\R \to \comp$, 
\[
	E_n^{\tr}(\v) = {\psi}_n \left(\earthquakeTrace({\phi}_n(\v))\right).
\]

In triangle length coordinates the family of maps become
\[
\begin{split}
	\bar{\phi}_n: \compL &\to \compL \\
	\begin{pmatrix} a \\ b \\ c \end{pmatrix} &\mapsto \begin{pmatrix} a'(n) \\ b'(n) \\ c'(n) \end{pmatrix} = \begin{pmatrix} 	\cosh^{-1}\left(A(n)\right) \\ a \\ \cosh^{-1}\left(A(n-1)\right)\end{pmatrix}
\end{split}
\]
with 
\[
	A(n)=\cosh{b}\cosh{\left(na\right)}+\left(\cosh{c}-\cosh{a}\cosh{b}\right)\frac{\sinh{\left(na\right)}}{\sinh{a}}
\]
and
\[
\begin{split}
	\bar{\psi}_n: \compL &\to \compL  \\
	\begin{pmatrix} a' \\ b' \\ c' \end{pmatrix} &\mapsto \begin{pmatrix} a(n) \\ b(n) \\ c(n) \end{pmatrix} = \begin{pmatrix} b' \\ \cosh^{-1}\left(B(n)\right) \\ \cosh^{-1}\left(B(n-1)\right) \end{pmatrix}
\end{split}
\]
where 
\[
	{B}(n)=\cosh{a'}\cosh{\left(nb'\right)}+\left(\cosh{c'}-\cosh{a'}\cosh{b'}\right)\frac{\sinh{\left(nb'\right)}}{\sinh{b'}}.
\]

Then given some starting point $\w = (a,b,c) \in \RPlusSpace$ the earthquake about the family $\gamma_n$ is given by the map $E_{n}^\ell(\w):\R \to \RPlusSpace$,
\[
	E_n^{\ell}(\w) = \bar{\psi}_n \left(\earthquakeHyp(\bar{\phi}_n(\w))\right).
\]

We can calculate the limit as $n\to\infty$ by expanding the right-hand side expression. The sequence $\gamma_n=\alpha\beta\alpha^{n-1}$ asymptotically approaches $\alpha$ and thus we expect the earthquake about $\gamma_n$ to converge to the earthquake about $\alpha$. We use hyperbolic trigonometry identities in addition to the collar equation to simplify expressions. 

Rescale using $r\mapsto \sfrac{r}{n^2}$ to account for the growing length of the curve and the appropriate sequence of transverse measures in $n$. Using Taylor series expansions to evaluate certain terms in the limit we find,
\[
	\lim_{n\to\infty}E_n^{\ell}(\w)\left(\frac{r}{n^2}\right) = E_{\alpha}^{\ell}(\w)(r).
\]
 
 \begin{figure}[hbt]
	\centering
	\includegraphics[width=0.5\textwidth]{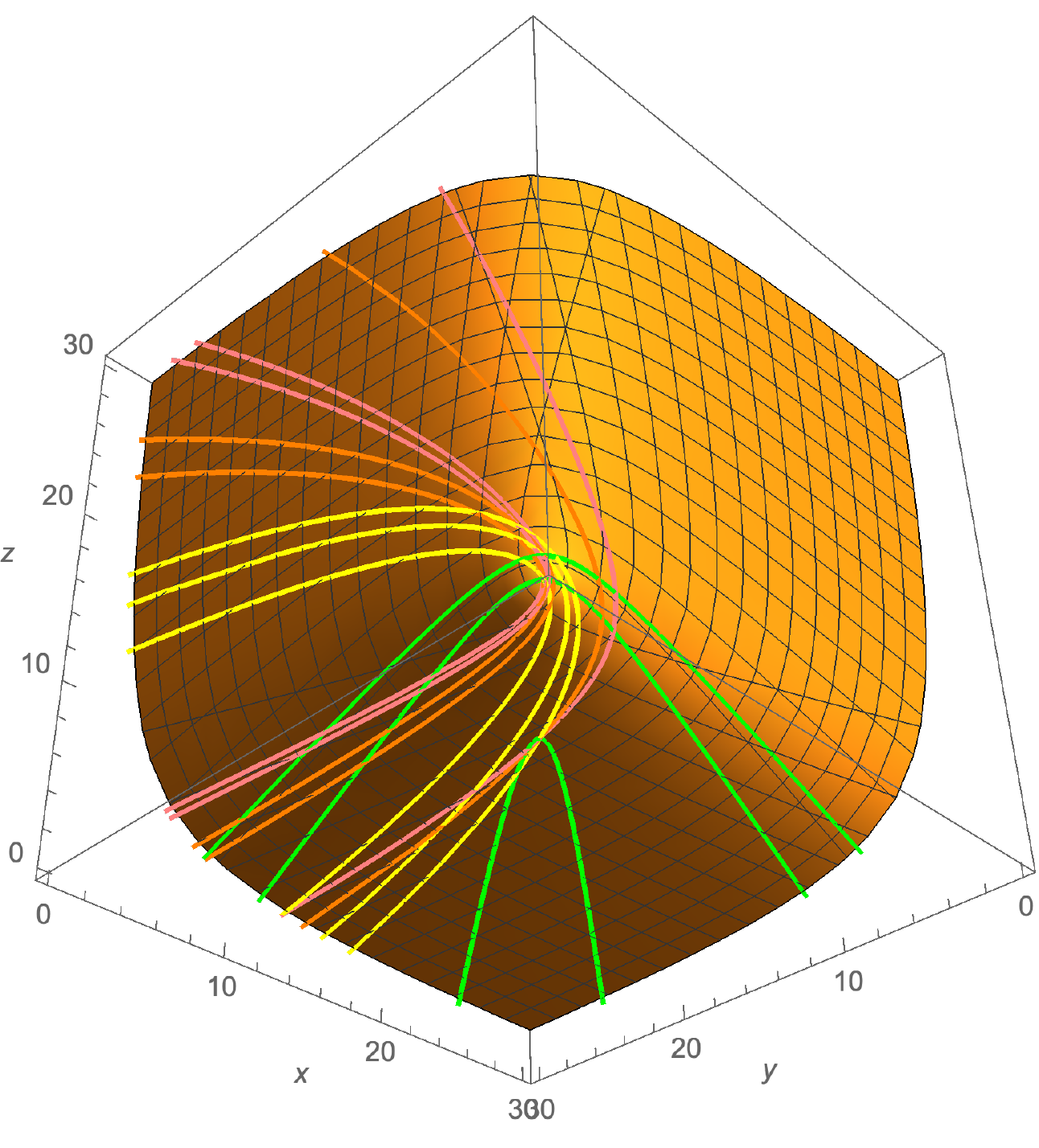}
	\includegraphics[width=0.5\textwidth]{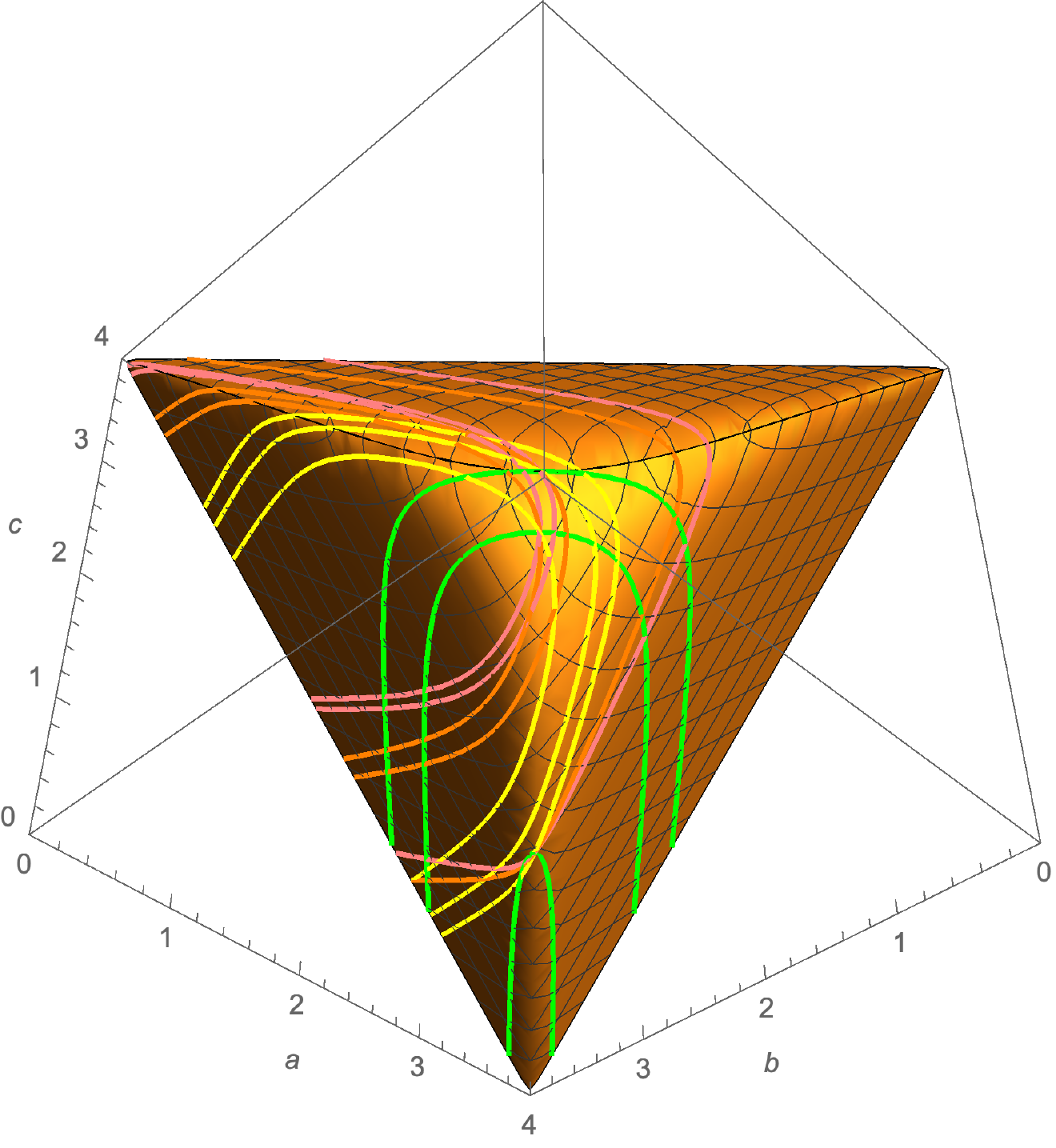}
	\caption{Some example earthquake deformations about $\gamma_n = \alpha\beta\alpha^{n-1}$ for in trace coordinates (top) and triangle length coordinates (bottom). Examples shown are for $n=1,2,3,4$ in green, yellow, orange, pink, respectively. Starting points are given by $\mathcal{S}_\alpha$ (top) and $\nu\mathcal{S}_\alpha$ (bottom) from Equations \ref{eqn:startPoint} and  \ref{eqn:coordChange}.}
	\label{fig:example1}
\end{figure}

\FloatBarrier
If we used the hyperbolic length parametrisation as in Remark \ref{rem:reparametriseHyp} the appropriate rescaling is $s\mapsto \sfrac{s}{n}$, 
\begin{equation}
	\lim_{n\to\infty}E_n^{\ell}(\w)\left(\frac{s}{n}\right) = E_{\alpha}^{\ell}(\w)(s).
\end{equation}

Note the rescaling is important, otherwise the limit would tend to $\infty$ in multiple coordinates. This result is consistent with Kerckhoff's result \cite{Ker83} that defining earthquakes for general measured geodesic laminations through the limits of sequences is well-defined.

See Figure \ref{fig:example1} for examples of the orbits of $\gamma_n$ for $n=1,2,3,4$.
\end{example}

\begin{example}[$\gamma_n = T^n(\alpha)$]

Take the following mapping class,
\[
	T=T_\beta T_\alpha^{-1} = 
	\begin{cases} 
		\alpha \mapsto \alpha\beta \\ \beta \mapsto \beta\alpha\beta
	\end{cases},
\]
with matrix representation
\[
	N= \begin{pmatrix} 1&1\\1&2\end{pmatrix}.
\]

This is a pseudo-Anosov mapping class with irrational eigenvalues $\lambda_\pm=\sfrac{\left(3\pm\sqrt{5}\right)}{2}$. Repeated applications of $T$ to $\alpha$ will converge to a lamination $\ell$ with irrational slope $\sfrac{(1+\sqrt{5})}{2}$. That is, a measured geodesic lamination that is not a simple closed curve.

First, compute a sequence of simple closed curves $\delta_n$ such that $i(\gamma_n,\delta_n)=1$ for each $n\in\N$. One such example is $\delta_n = T^n(\beta)$. Then
\[
	\begin{split}
		\gamma_n 	&= T^n(\alpha), \hspace{0.5cm} \delta_n = T^n(\beta),\\
		\alpha 		&= T^{-n}(\gamma_n), \hspace{0.275cm} \beta = T^{-n}(\delta_n).
	\end{split}
\]

Examples for $n=1,2$ are shown below,
\[
	\begin{split}
		\gamma_1 	&= \alpha\beta, \hspace{0.65cm} \delta_1 	= \beta\alpha\beta,\\
		\alpha 		&= \gamma_1^2\delta_1^{-1}, \hspace{0.275cm} \beta 	= \delta_1\gamma_1^{-1},
	\end{split}
\]
\[
	\begin{split}
		\gamma_2 	&= \alpha\beta^2\alpha\beta, \hspace{1.75cm} \delta_2 = \beta\alpha\beta\alpha\beta^2\alpha\beta,\\
		\alpha 		&= \gamma_2^2\delta_2^{-1}\gamma_2^2\delta_2^{-1}\gamma_2\delta_2^{-1}, \hspace{0.2cm} \beta = \delta_2\gamma_2^{-1}\delta_2\gamma_2^{-2}.
	\end{split}
\]

Note that as $n\to\infty$, $\gamma_n$ and $\delta_n$ will converge to the same measured geodesic lamination with the same irrational slope, which we can use to our advantage.

Second, calculate the earthquake deformation for $\gamma_n$ and $\delta_n$ in the framing $\gamma_n$ and $\delta_n$ as given in Equation \ref{eqn:algEarthquake}. 
	
Third, construct the series of maps $\phi_n=\tilde{T}^n$ and $\psi_n=\tilde{T}^{-n}$ to go between the original framing $\alpha$, $\beta$ and the new framing $\gamma_n$, $\delta_n$ for each $n$. Both $\tilde{T}^{n}$ and $\tilde{T}^{-n}$ involve a system of recurrence relations that are not linear and we have found no typical method for solving them. We can still understand the underlying maps $\tilde{T}$ and $\tilde{T}^{-1}$,
\[
\begin{split}
	\tilde{T}: \comp 						&\to 		\comp \\
	\begin{pmatrix} x \\ y \\ z \end{pmatrix} 	&\mapsto  \begin{pmatrix} x' \\ y' \\ z' \end{pmatrix} 
	= \begin{pmatrix} \tr(\rho(\alpha\beta)) \\ \tr(\rho(\beta\alpha\beta)) \\ \tr(\rho(\alpha\beta^2\alpha\beta)) \end{pmatrix} = \begin{pmatrix} z \\ yz-x \\ z\left(yz-x\right)-y\end{pmatrix}
\end{split}
\]
and
\[
\begin{split}
	\tilde{T}^{-1}: \comp 					&\to 		\comp \\
	\begin{pmatrix} x' \\ y' \\ z' \end{pmatrix} 	&\mapsto \begin{pmatrix} x \\ y \\ z \end{pmatrix} 	
	= \begin{pmatrix} \tr(\rho(\gamma_n^2\delta_n^{-1})) \\ \tr(\rho(\delta_n\gamma_n^{-1})) \\ \tr(\rho(\gamma_n)) \end{pmatrix} = \begin{pmatrix} x'\left(x'y'-z'\right)-y' \\ x'y'-z' \\ x' \end{pmatrix}.
\end{split}
\]

Given some starting point $\v = (x,y,z) \in \comp$ we can define the earthquake about the family $\gamma_n$ recursively by $E_{n}^{\text{tr}}(\v):\R \to \comp$,
\[
\begin{split}
	E_n^{\text{tr}}(\v) 	&= \tilde{T}^{-n} \left(\earthquakeTrace(\tilde{T}^{n}(\v))\right),\\
					&= \tilde{T}^{-1}\left(E_{n-1}^{\text{tr}}(\tilde{T}(\v))\right).
\end{split}
\]

\begin{figure}[hbt]
	\centering
	\includegraphics[width=0.5\textwidth]{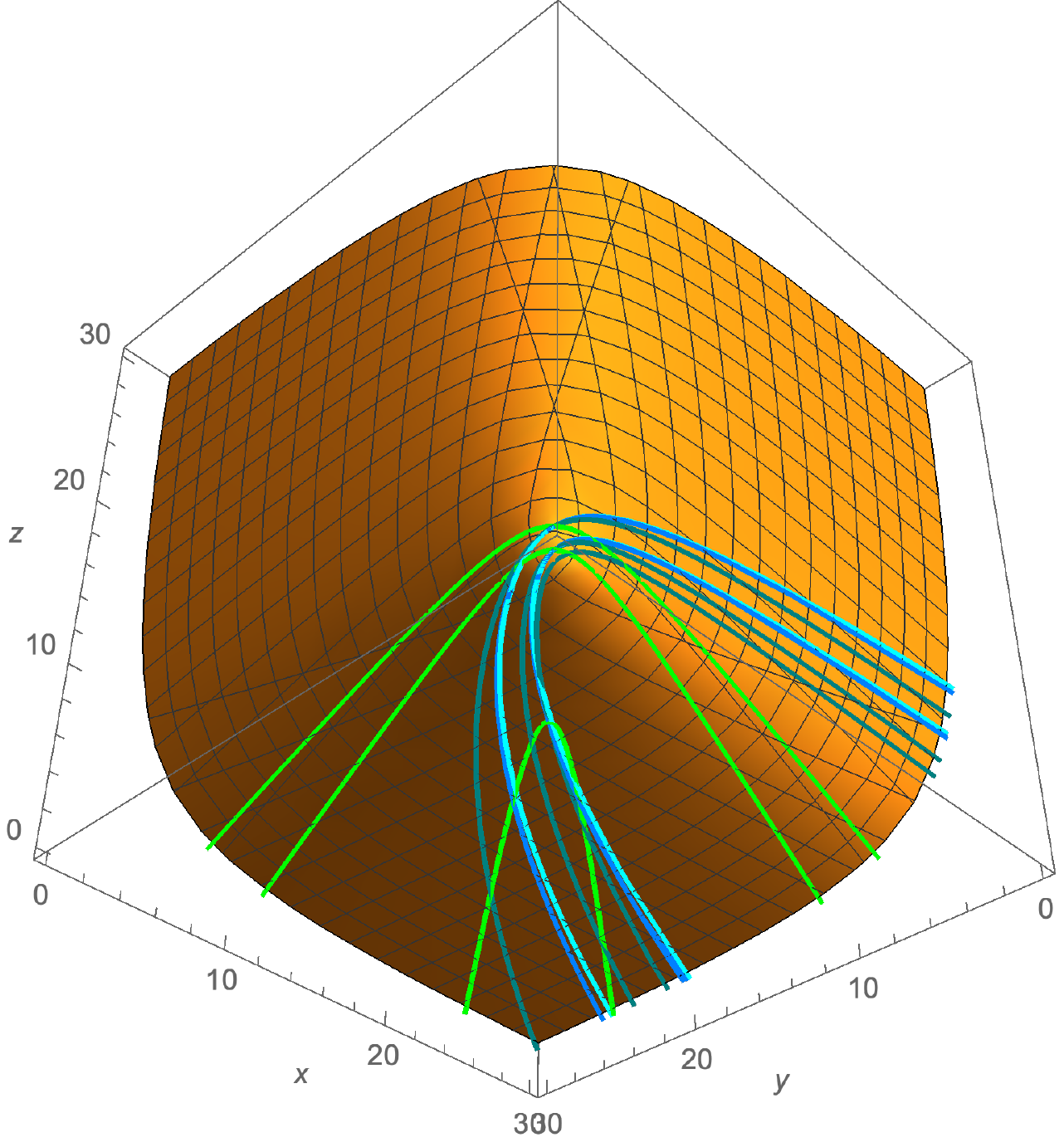}
	\includegraphics[width=0.5\textwidth]{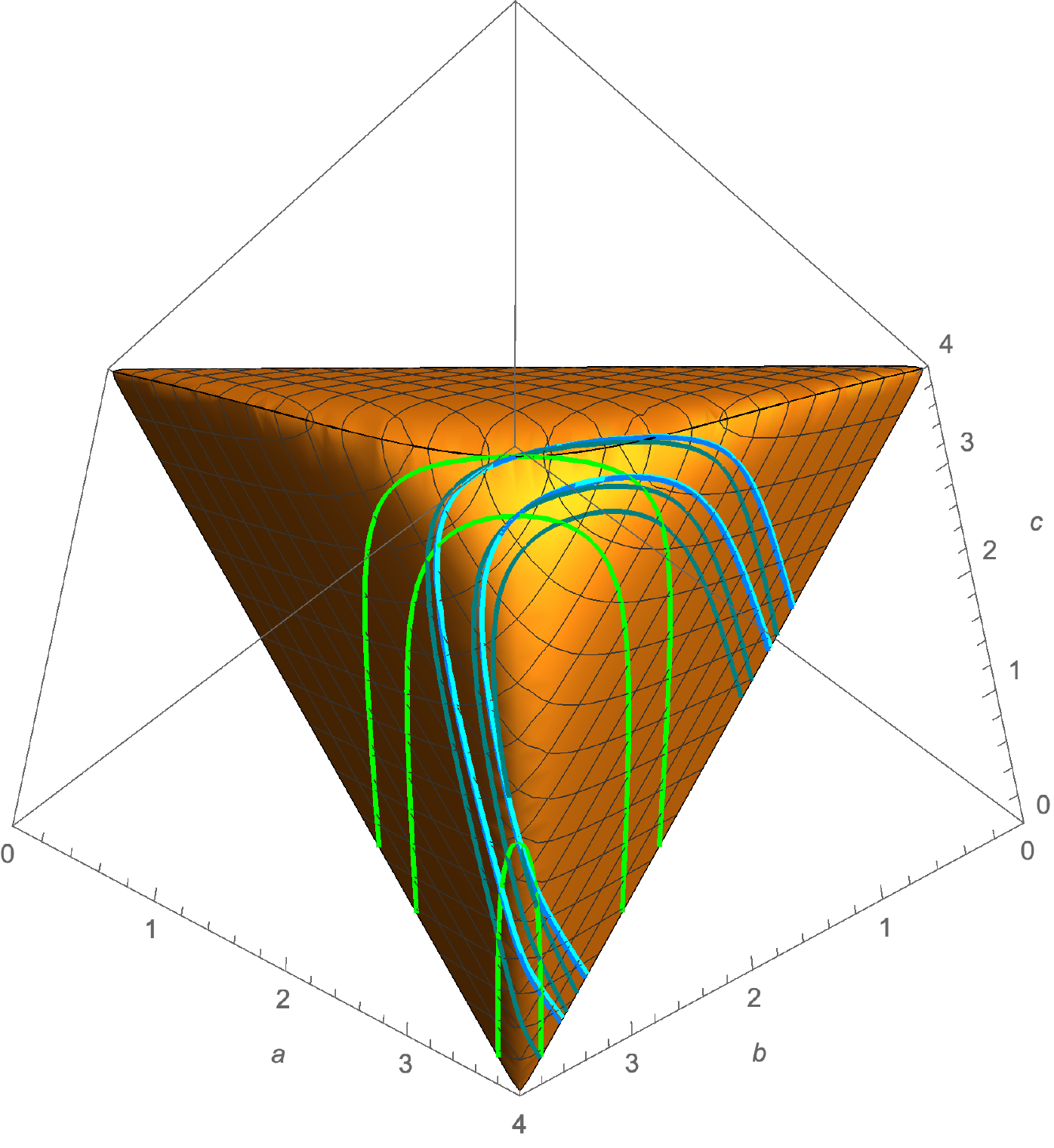}
	\caption{Some example earthquake deformations about $\gamma_n = T^n(\alpha)$ for $n=1,2,3,4$ in trace coordinates (top) and triangle length coordinates (bottom). Examples shown are for $n=1,2,3,4$ in green, turquoise, light blue, and cyan, respectively. Starting points are given by $\mathcal{S}_\alpha$ (top) and $\nu\mathcal{S}_\alpha$ (bottom) from Equations \ref{eqn:startPoint} and  \ref{eqn:coordChange}.}
	\label{fig:example2}
\end{figure}

In triangle length coordinates the maps $\tilde{T}$ and $\tilde{T}^{-1}$ become,
\[
\begin{split}
	\bar{T}: \RPlusSpace 				&\to 		\RPlusSpace \\
	\begin{pmatrix} a \\ b \\ c \end{pmatrix} 	&\mapsto \begin{pmatrix} a' \\ b' \\ c' \end{pmatrix} 
	= \begin{pmatrix} c \\ \cosh^{-1}\left(2\cosh{b}\cosh{c}-\cosh{a}\right) \\ \cosh^{-1}\left(2\cosh{c}\left(2\cosh{b}\cosh{c}-\cosh{a}\right)-\cosh{b}\right) \end{pmatrix}
\end{split}
\]
and
\[
\begin{split}
	\bar{T}^{-1}: &\RPlusSpace 			\to 		\RPlusSpace \\
	\begin{pmatrix} a' \\ b' \\ c' \end{pmatrix} 	&\mapsto \begin{pmatrix} a \\ b \\ c \end{pmatrix} 	
	= \begin{pmatrix} \cosh^{-1}\left(2\cosh{a'}\left(2\cosh{a'}\cosh{b'}-\cosh{c'}\right)-\cosh{b'}\right) \\ \cosh^{-1}\left(2\cosh{a'}\cosh{b'}-\cosh{c'}\right) \\ a' \end{pmatrix}
\end{split}	
\]

Given some starting point $\w = (a,b,c) \in \comp$ we can define the earthquake about the family $\gamma_n$ recursively by $E_{n}^{\ell}(\w):\R \to \RPlusSpace$,
\[
	\begin{split}
	E_n^{\ell}(\w) 	&= \bar{T}^{-n} \left(\earthquakeHyp(\bar{T}^{n}(\w))\right),\\
				&= \bar{T}^{-1}\left(E_{n-1}^{\ell}(\bar{T}(\w))\right).
	\end{split}
\]
We can similarly define the earthquake about the family $\delta_n$ in trace coordinates and triangle lengths.

Again, to assess the limit as $n\to\infty$ a rescaling factor is required, otherwise the length and the transverse measure would become unbounded. Recall $\lambda_+$ is the largest eigenvalue of the matrix representation $N$, $\lambda_+=\sfrac{\left(3+\sqrt{5}\right)}{2}$. The appropriate rescaling to take for our original parametrisation is $r\mapsto\sfrac{r}{\lambda_+^{2n}}$ and for our hyperbolic length reparametrisation is $s\mapsto\sfrac{r}{\lambda_+^{n}}$. 

See Figure \ref{fig:example2} for examples of the orbits of $\gamma_n$ for {$n=1,2,3,4$}. Note the variation in the earthquake from $n=3$ to $n=4$ is already getting small. The recursive formula makes it difficult to evaluate any terms in the limit directly. 

Compare the difference between $\gamma_n$ and $\delta_n$ as $n$ increases. Given Thurston's earthquake theorem, we know this narrows the region of $\T(\punc)$ where the limiting earthquake can appear. Given the relationship between slope in Theorem \ref{thm:earthquakeAsymptotics}, we know the slope of the limiting earthquake must be $\sfrac{2}{(1+\sqrt{5}}$. We can estimate the limiting earthquake from these two pieces of information. See Figure \ref{fig:example2FN} for examples of the forwards and backwards orbits of $\gamma_n$ and $\delta_n$ with an estimation of the limiting earthquake about $\ell$.

\begin{figure}[hbt]
	\centering
	\includegraphics[width=0.44\textwidth]{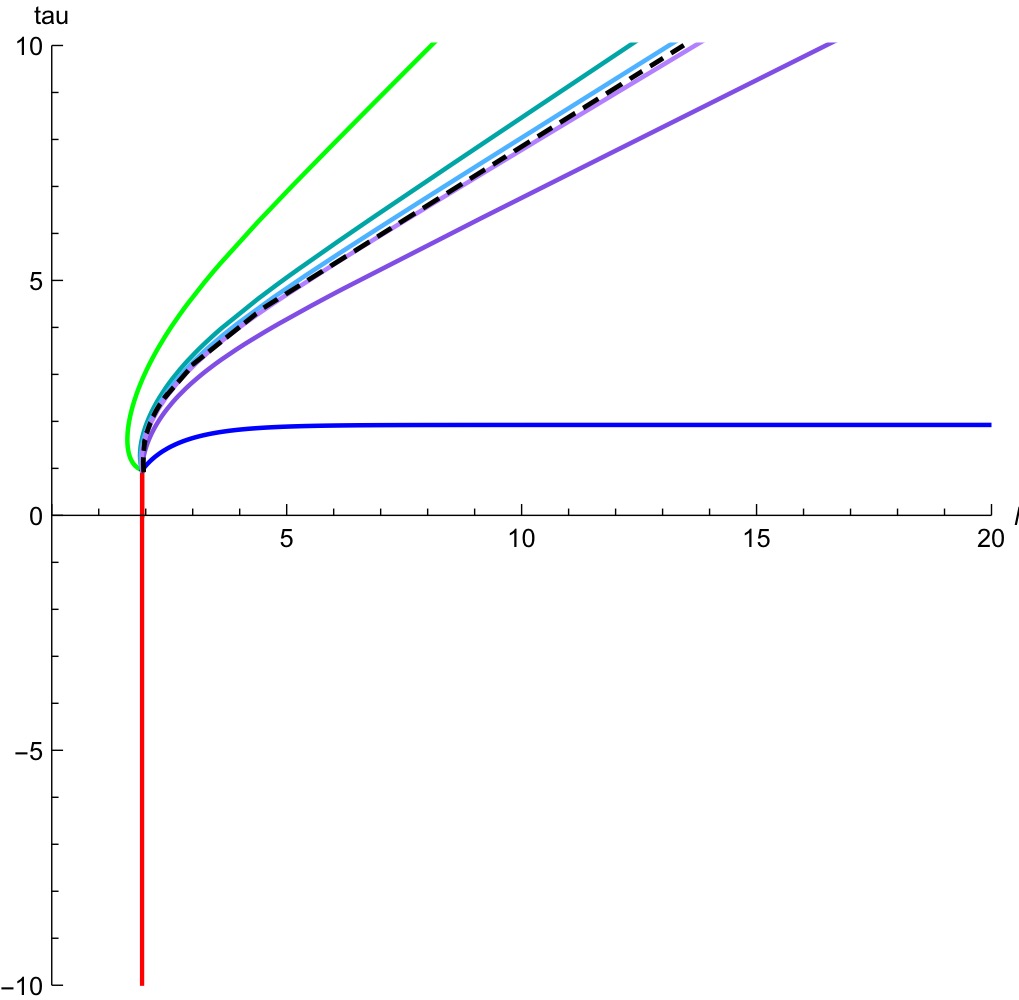}
	\includegraphics[width=0.44\textwidth]{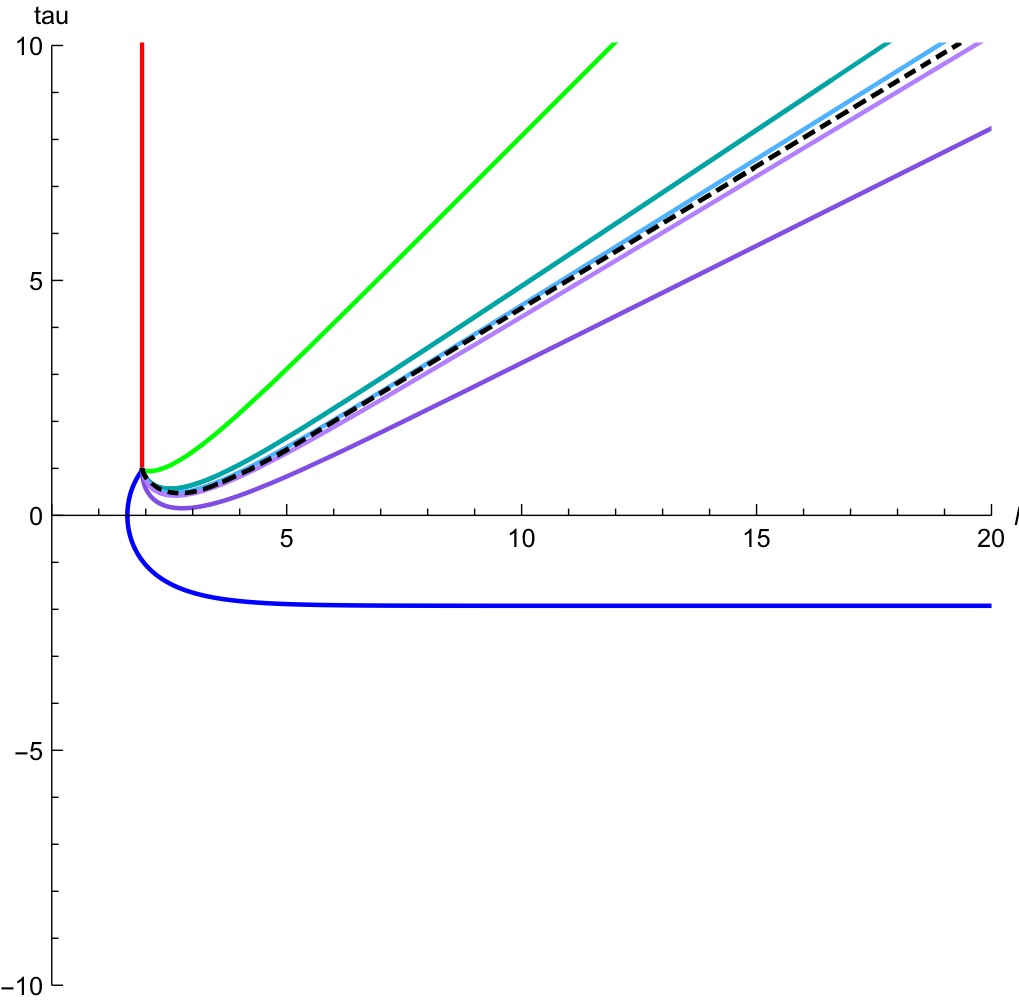}
	\caption{Example forwards (left) and backwards (right) earthquake deformations about $\gamma_n = T^n(\alpha)$ for $n=0,1,2,3$ and $\delta_n = T^n(\beta)$ for $n=0,1,2$ with starting point $\left(2\cosh^{-1}\left(\frac{3}{2}\right), 2\cosh^{-1}\left(\frac{\sqrt{5}}{2}\right)\right)$. The dashed black line denotes an estimation of the limiting earthquake between $\gamma_3$ and $\delta_2$.}
	\label{fig:example2FN}
\end{figure}

\end{example}

\FloatBarrier

\clearpage

\section{Pictures of earthquakes}
\label{sec:pics}

We use the results presented in Section \ref{sec:results} to generate pictures of earthquakes in multiple different coordinate systems. We focus on earthquakes for the framing $\alpha$, $\beta$, as well as those for the families of simple closed curves in the two examples. Denote the set of all of these curves to be $\curves$,
\[
	\curves=\{\alpha, \beta, \alpha\beta, \alpha\beta\alpha, \alpha\beta\alpha^2, \alpha\beta\alpha^3,T^2(\alpha),T^3(\alpha),T^4(\alpha)\}.
\]

We specifically look at trace coordinate starting points $\left(3,3,3\right)$, $\left(2\sqrt{2},2\sqrt{2},4\right)$, and $\left(10, 10, -10(-5 + \sqrt{23})\right)$. Note $(3,3,3)$ represents the hexagonal torus and $(2\sqrt{2},2\sqrt{2},4)$ represents the square torus. 

We consider the earthquakes in six coordinate systems: trace coordinates on the component $\comp$, two modified trace coordinates referred to as ``spherical trace coordinates" and ``inverted trace coordinates", triangle length coordinates, Fenchel-Nielsen length-twist coordinates, and the simplex. It would be useful to study the earthquakes directly on $\T(\punc)\cong\H^2$, but it is difficult to calculate change of coordinates between $\H^2$ and the other representations of $\T(\punc)$, see Abe \cite{Abe}.

\begin{figure}[H]
	\centering
	\includegraphics[width=0.275\textwidth]{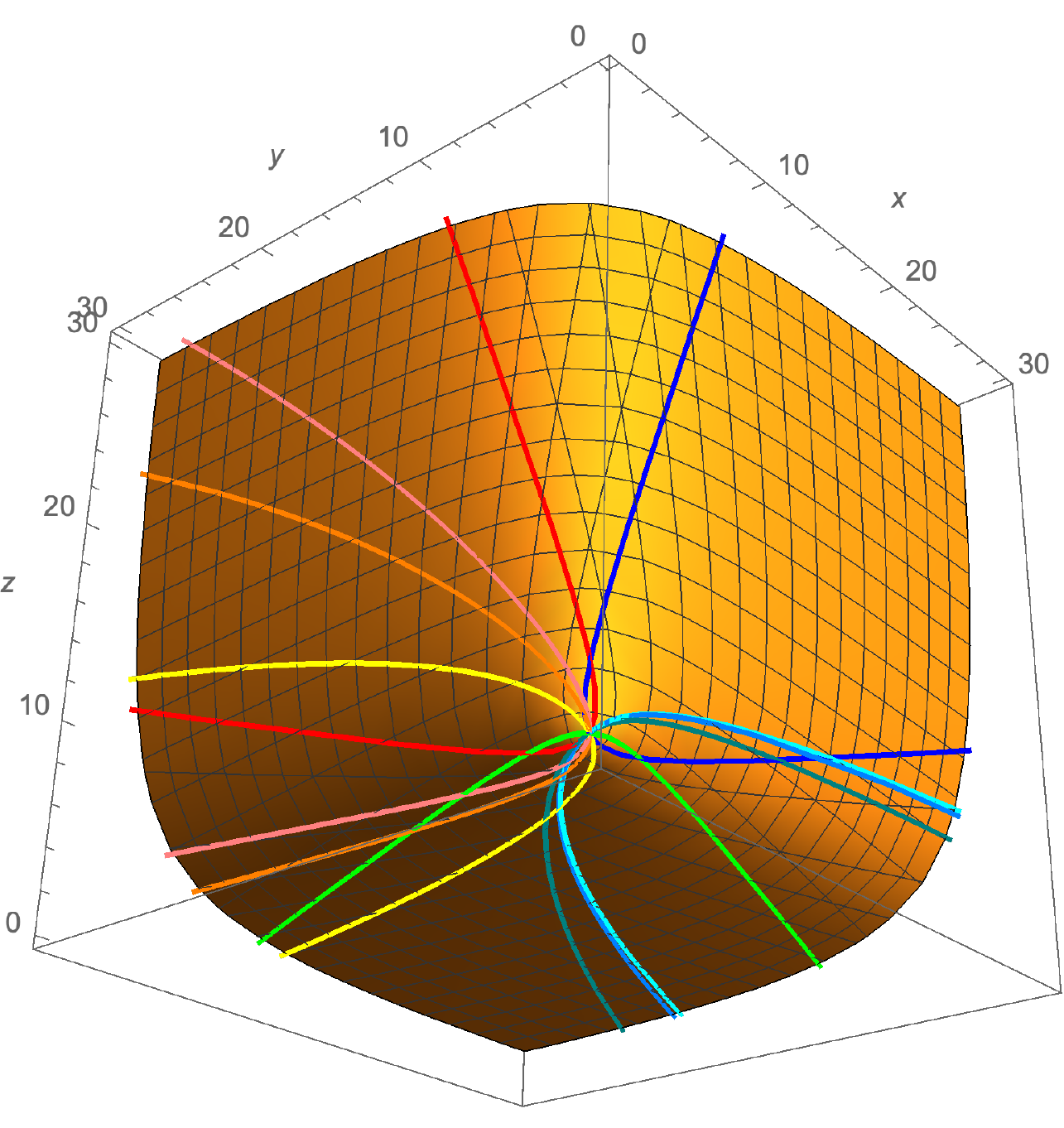}
	\includegraphics[width=0.275\textwidth]{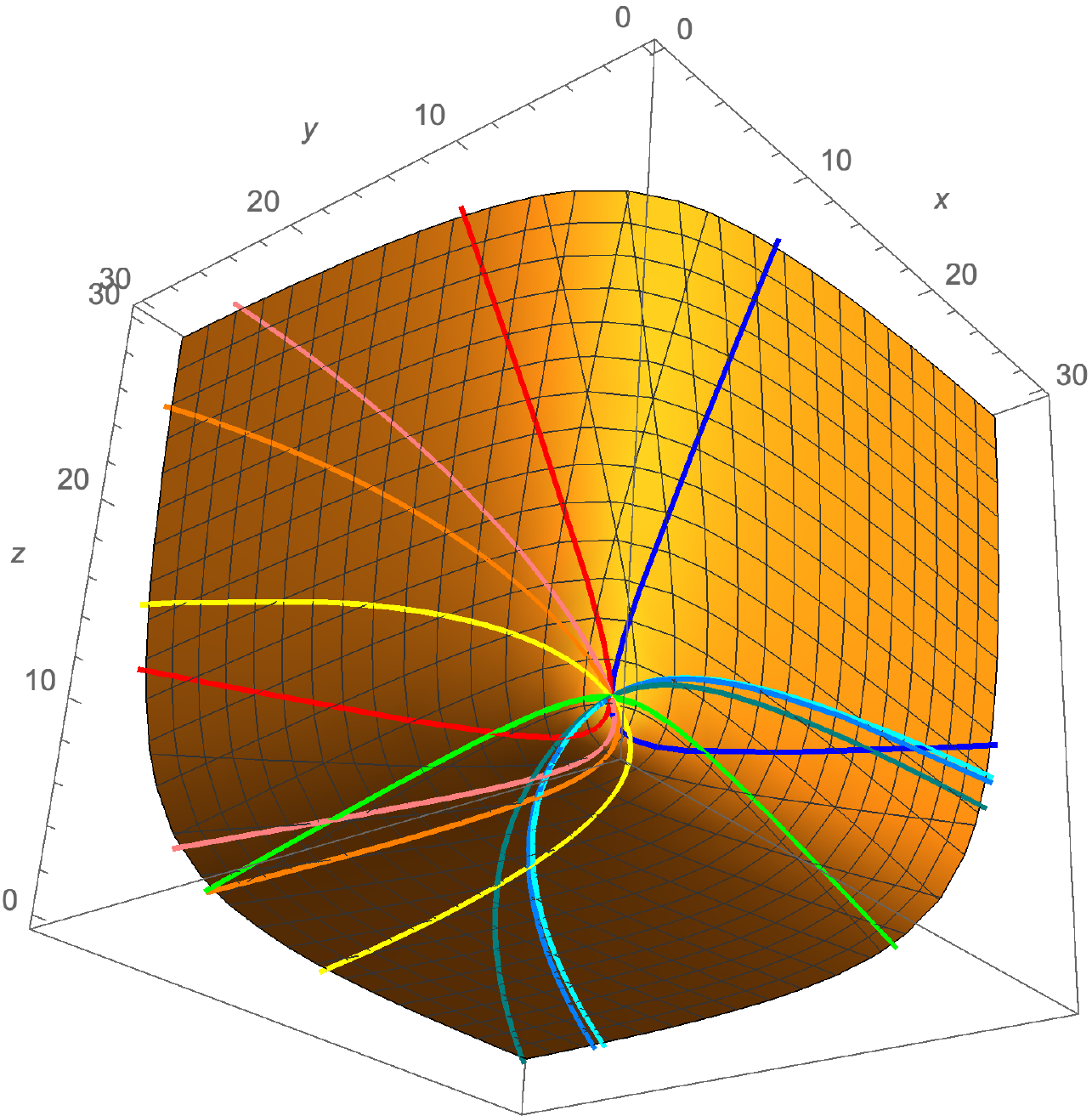}
	\includegraphics[width=0.275\textwidth]{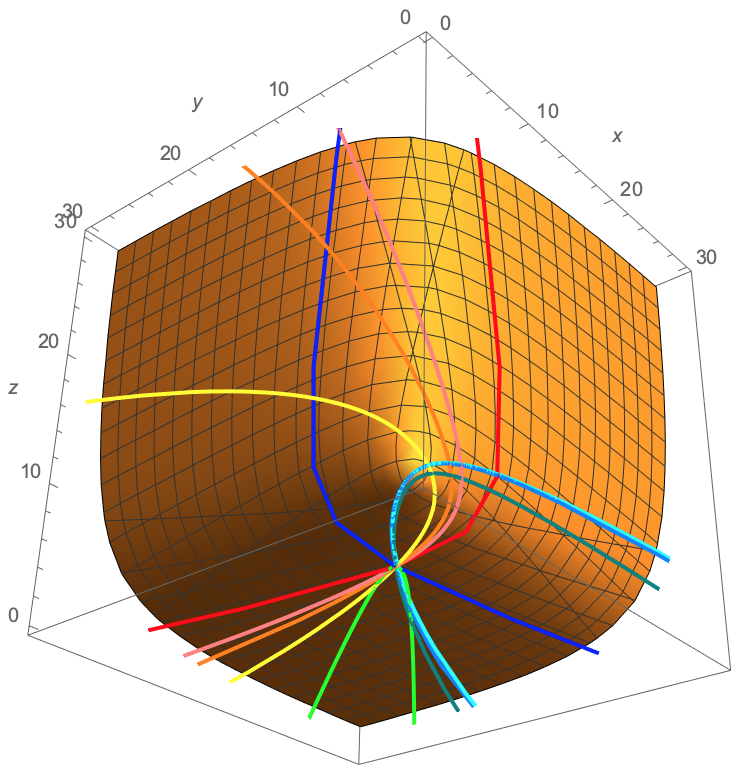}
	\caption{Forward and backward earthquake deformations about curves in $\curves$ on $\comp$ in trace coordinates given starting points $(3,3,3)$ (left), $\left(2\sqrt{2},2\sqrt{2},4\right)$ (centre), and $\left(10, 10, -10(-5 + \sqrt{23})\right)$ (right).}
	\label{fig:coordTrace}
\end{figure}

Recall trace coordinates $(x,y,z)$ on the component $\comp$ with $x^2+y^2+z^2-xyz=0$. Pictures of earthquakes about the curves in $\curves$ in trace coordinates are given in Figure \ref{fig:coordTrace}.

\begin{figure}[H]
	\centering	
	\includegraphics[width=0.275\textwidth]{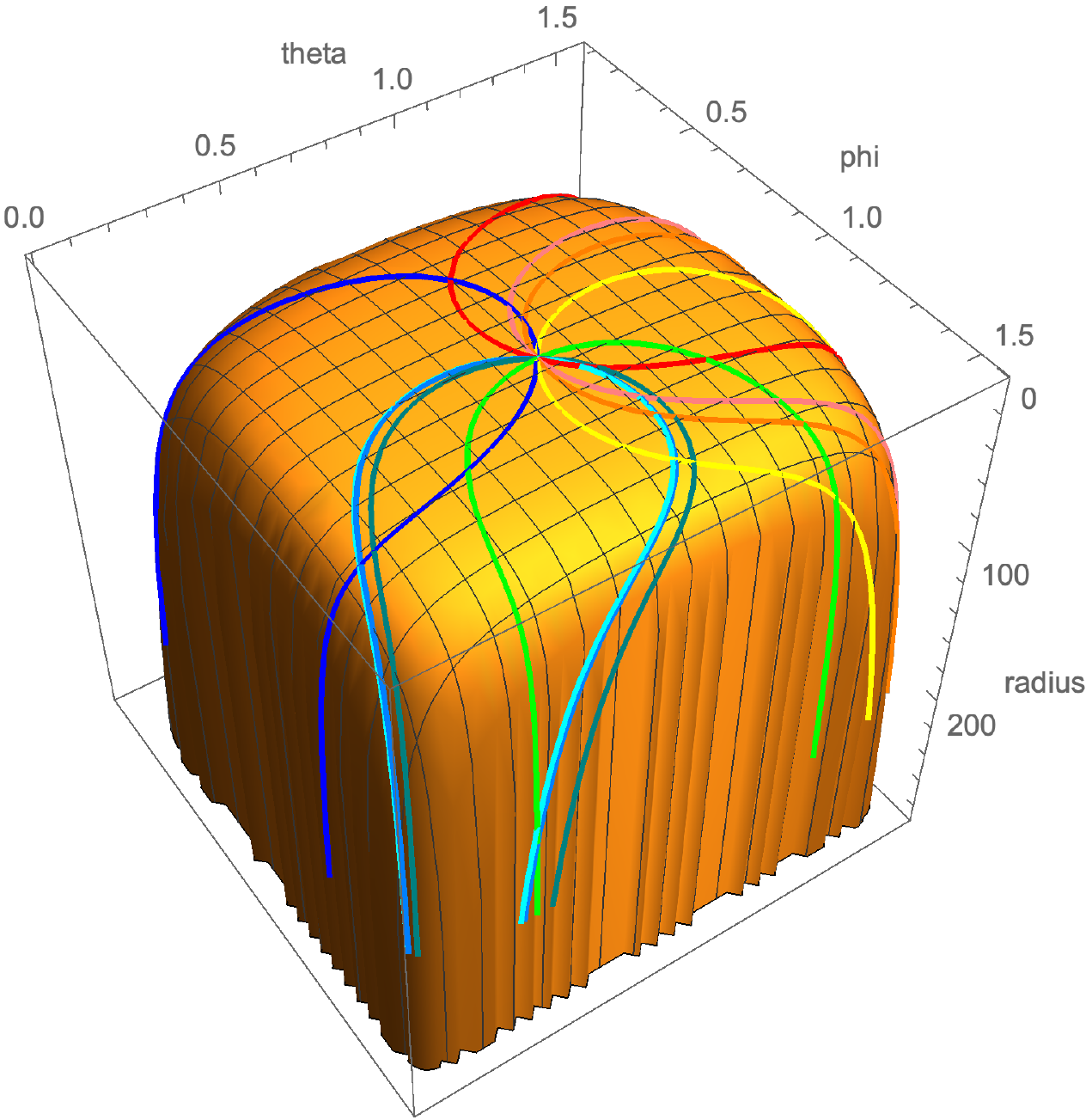}
	\includegraphics[width=0.275\textwidth]{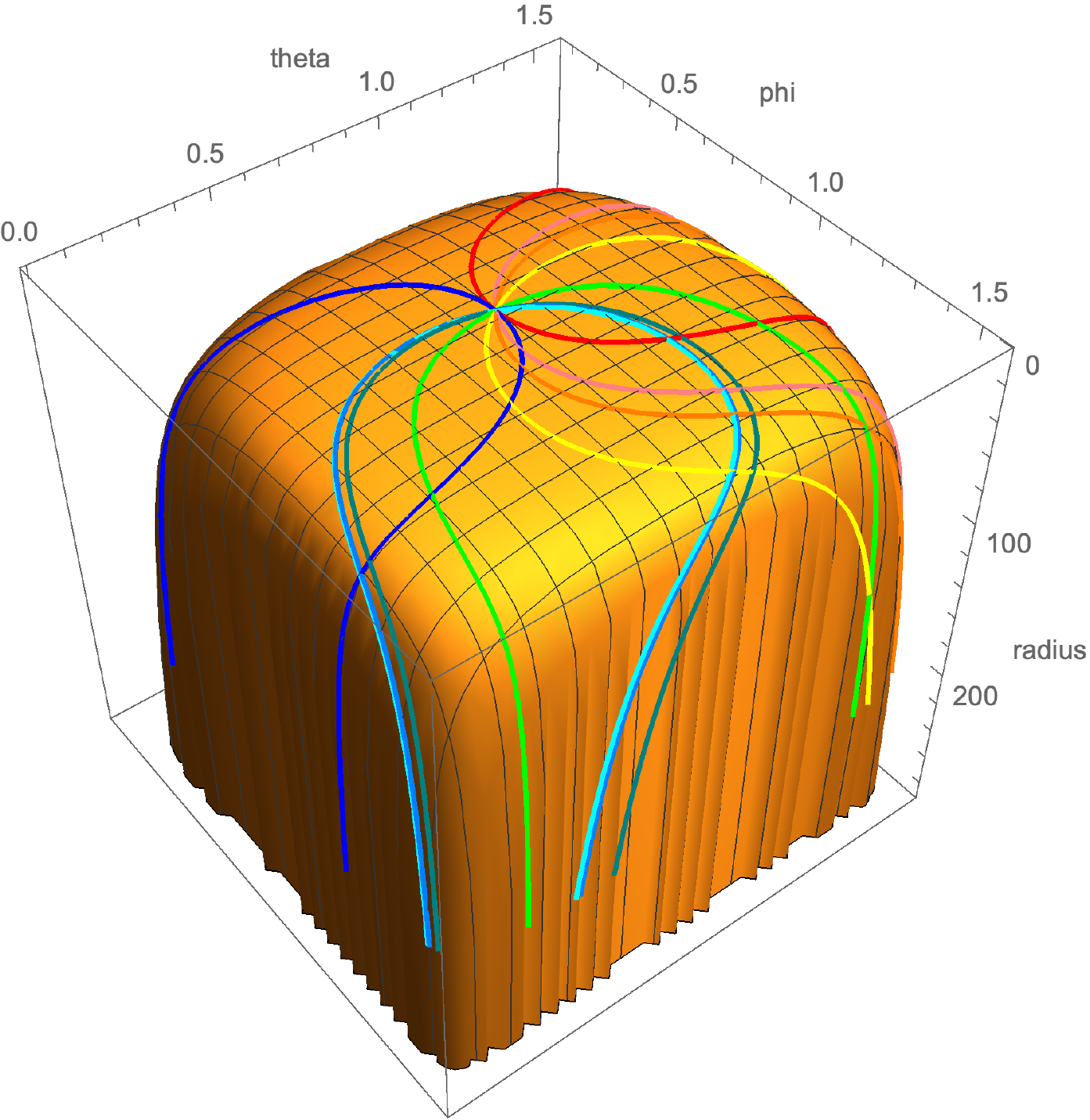}
	\includegraphics[width=0.275\textwidth]{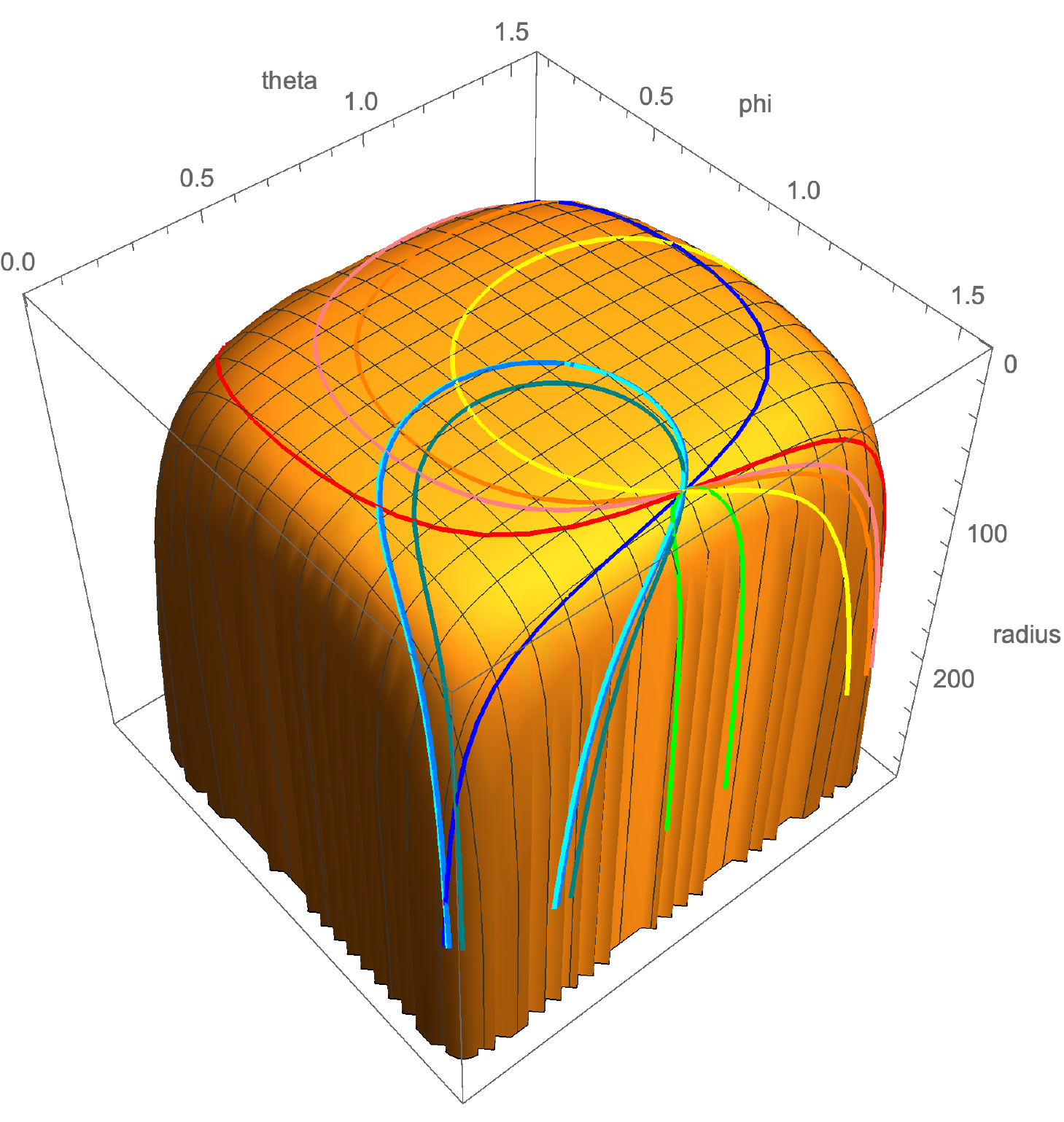}
	\caption{Forward and backward earthquake deformations about curves in $\curves$ on $\comp$ in spherical trace coordinates given starting points corresponding to trace coordinates $(3,3,3)$ (left), $\left(2\sqrt{2},2\sqrt{2},4\right)$ (centre), and $\left(10, 10, -10(-5 + \sqrt{23})\right)$ (right).}
	\label{fig:coordSph}
\end{figure}

\emph{Spherical trace coordinates} are given by applying a spherical coordinate change to $(x,y,z)$,
\[
	\begin{pmatrix} \theta \\ \phi \\ r \end{pmatrix} =  \begin{pmatrix} \tan^{-1}\left(\frac{y}{x}\right) \\ \tan^{-1}\left(\frac{\sqrt{x^2+y^2}}{z}\right) \\ \sqrt{x^2+y^2} \end{pmatrix}
\]
with $1-r\cos(\theta)\sin(\theta)\cos(\phi)\sin(\phi)^2=0$. Pictures of earthquakes about the curves in $\curves$ in spherical trace coordinates are given in Figure \ref{fig:coordSph}.

\begin{figure}[thb]
	\centering
	\includegraphics[width=0.275\textwidth]{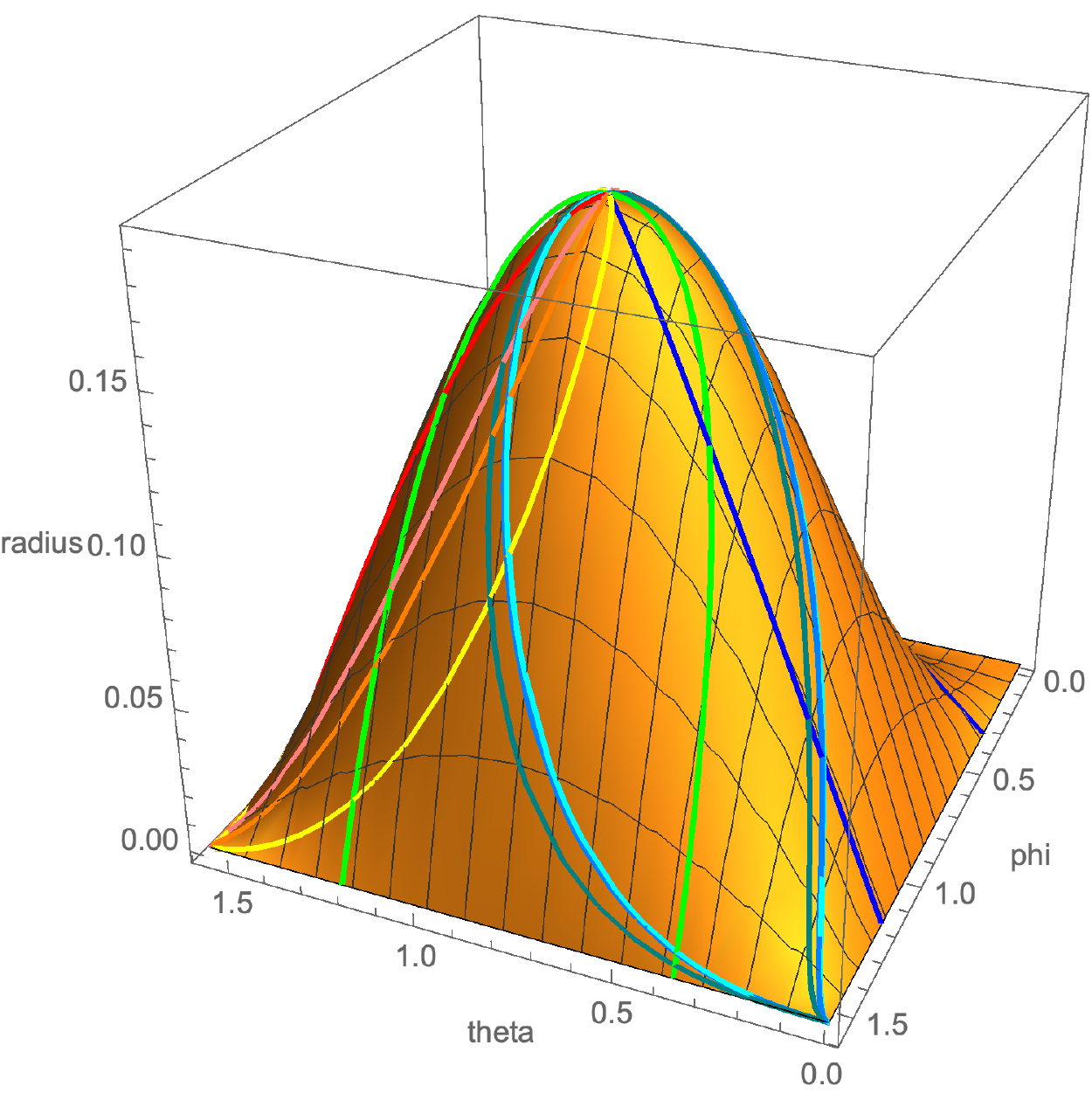}
	\includegraphics[width=0.275\textwidth]{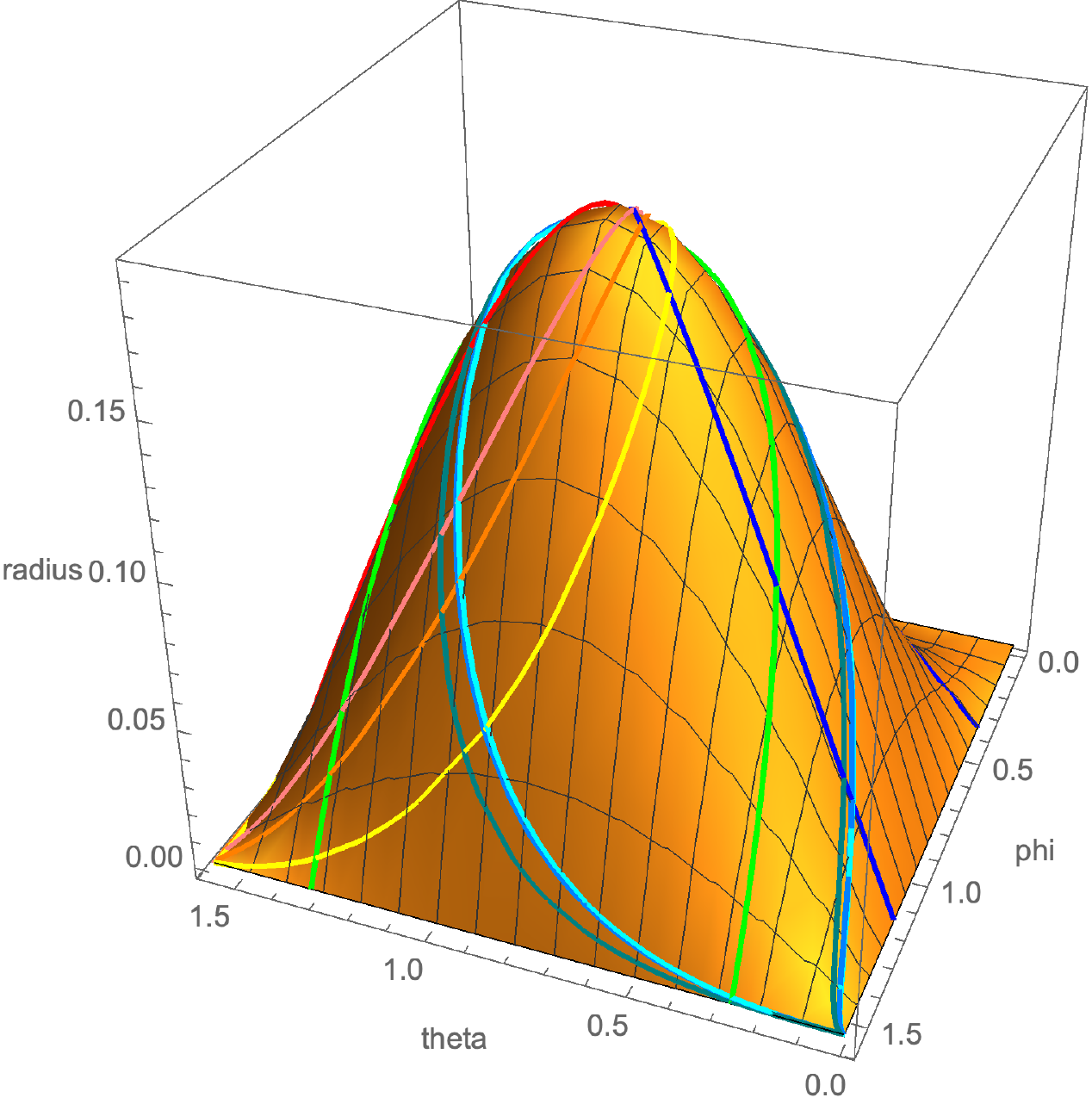}
	\includegraphics[width=0.275\textwidth]{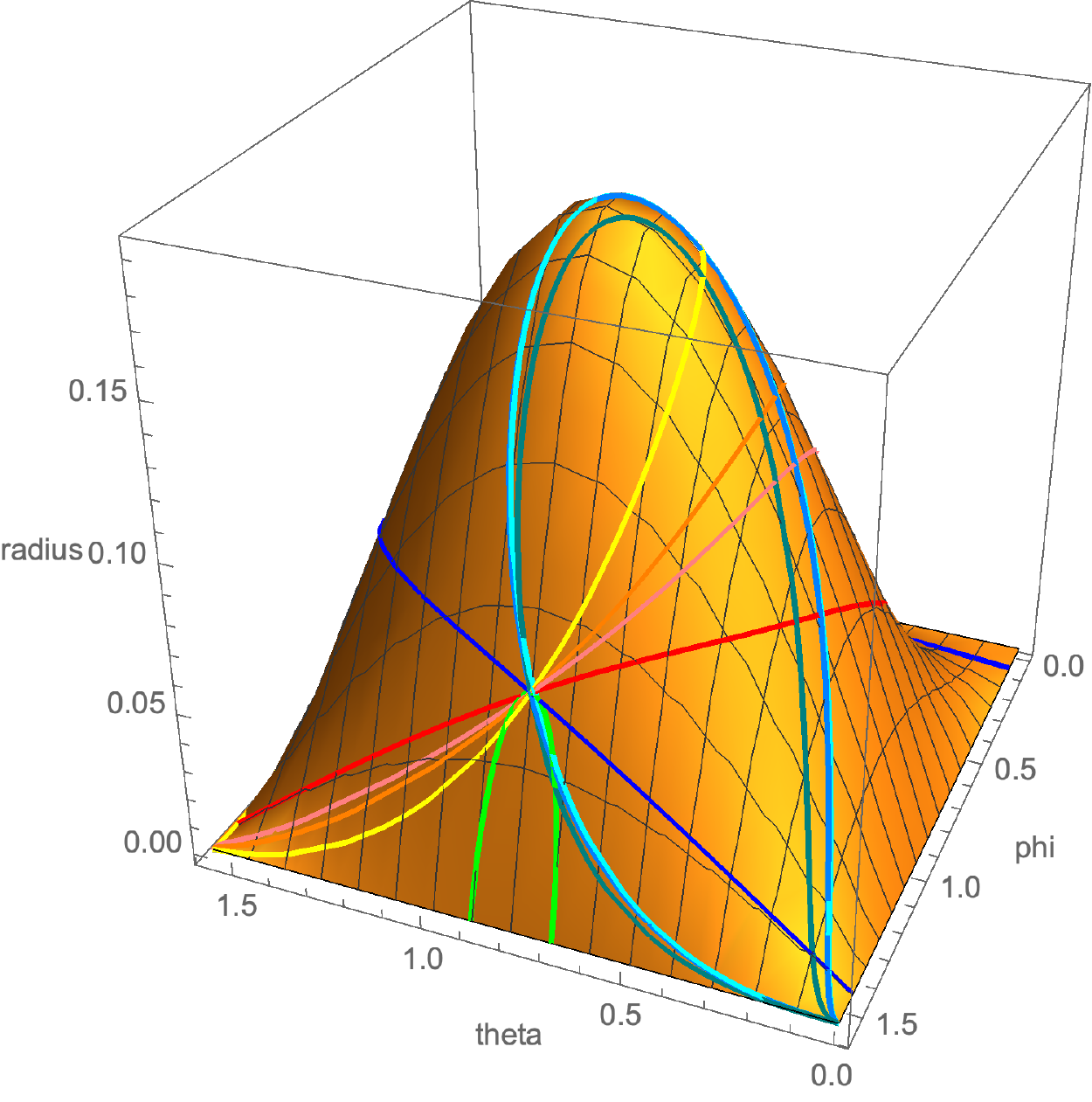}
	\caption{Forward and backward earthquake deformations about curves in $\curves$ on $\comp$ in inverted trace coordinates given starting points corresponding to trace coordinates $(3,3,3)$ (left), $\left(2\sqrt{2},2\sqrt{2},4\right)$ (centre), and $\left(10, 10, -10(-5 + \sqrt{23})\right)$ (right).}
	\label{fig:coordInv}
\end{figure}

\begin{figure}[htb]
	\centering	
	\includegraphics[width=0.29\textwidth]{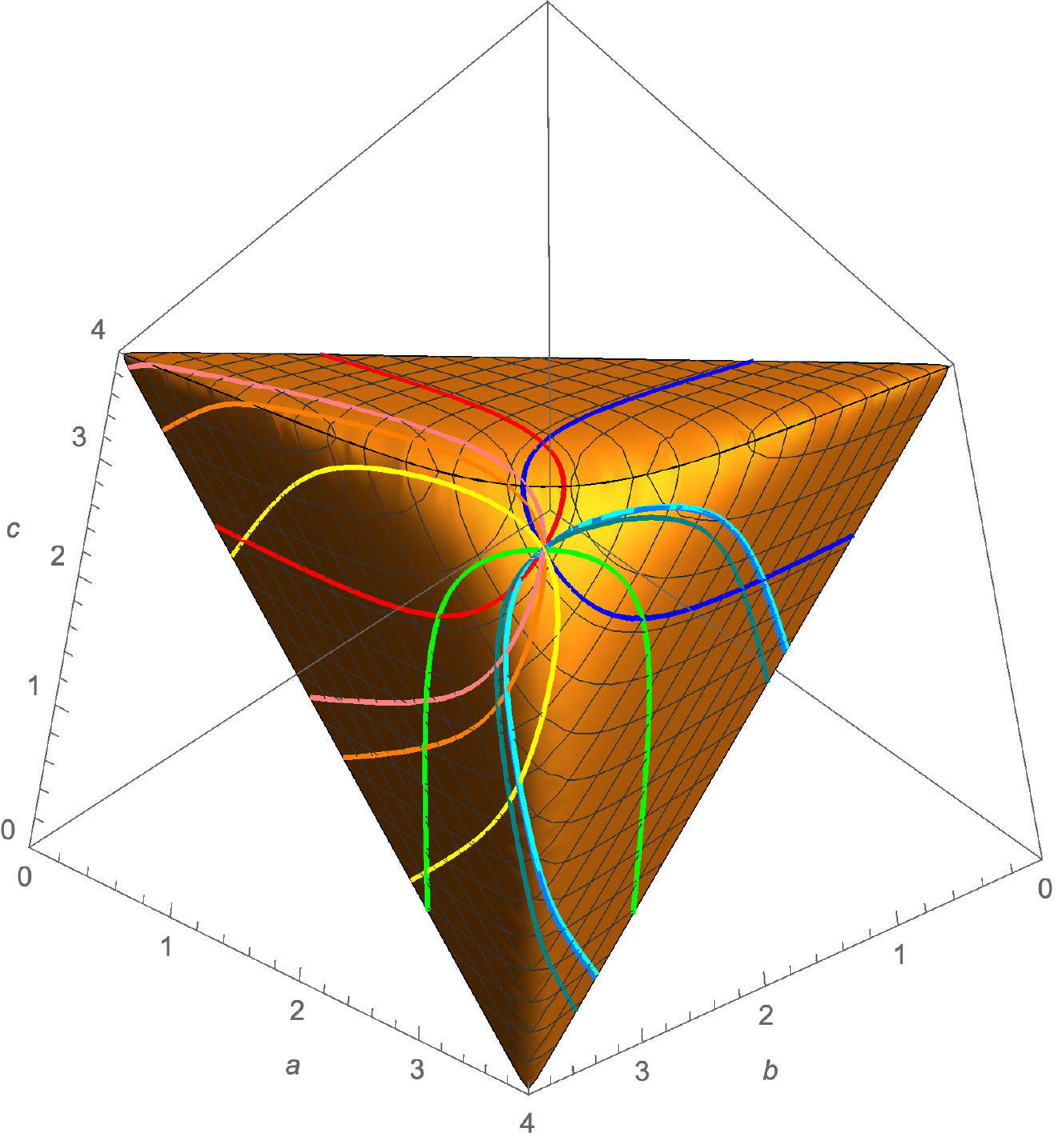}
	\includegraphics[width=0.29\textwidth]{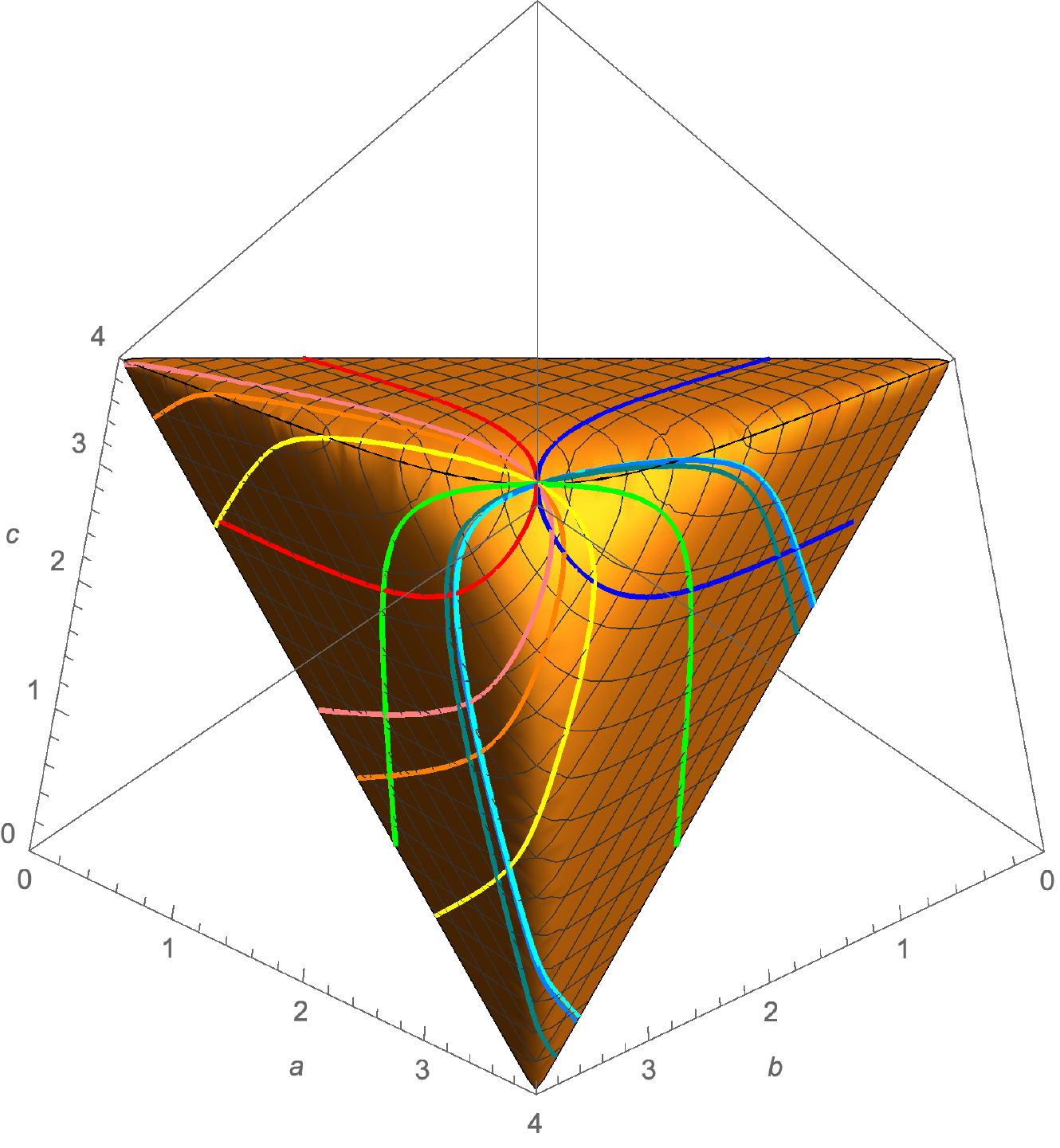}
	\includegraphics[width=0.29\textwidth]{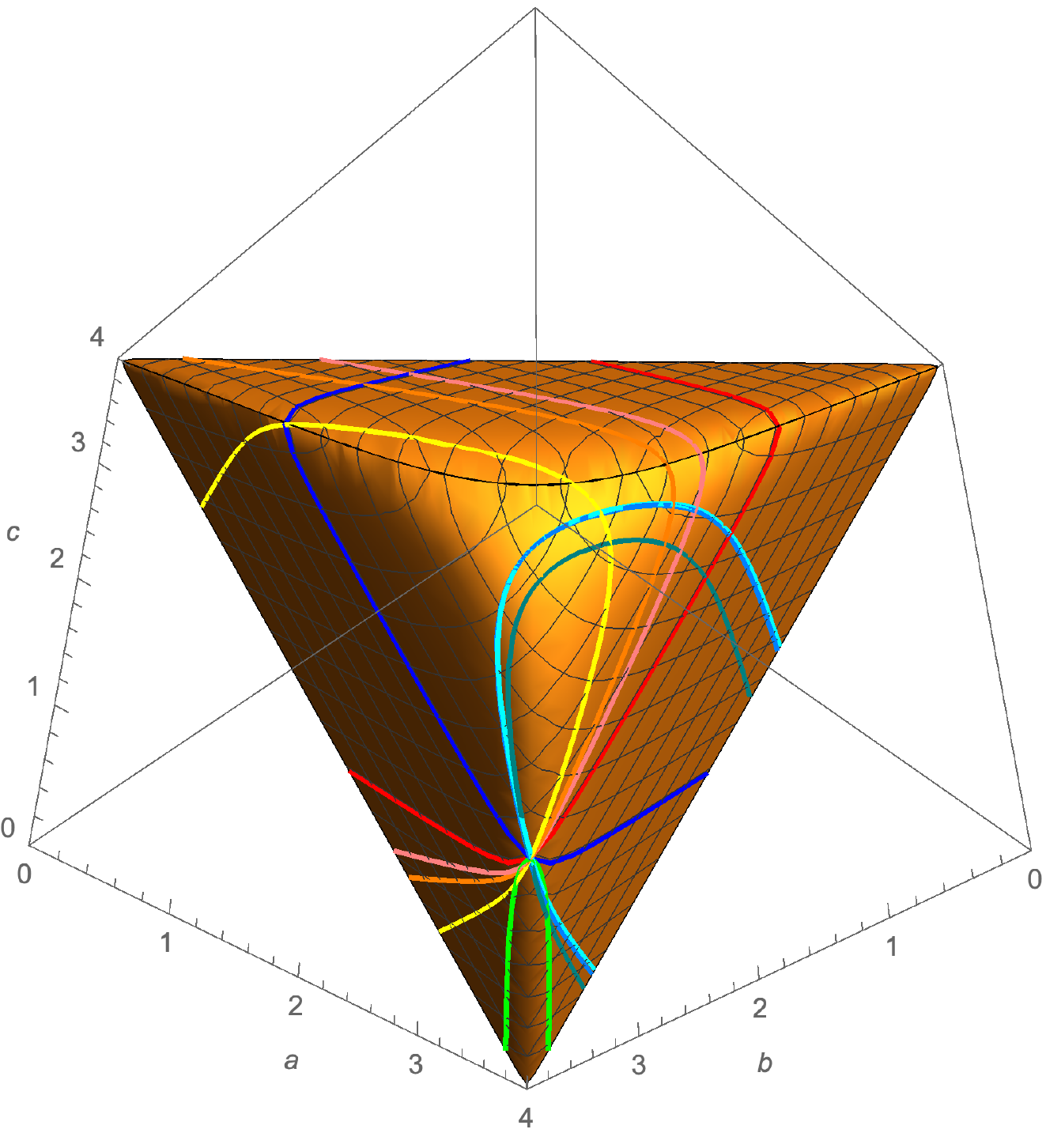}
	\caption{Forward and backward earthquake deformations about curves in $\curves$ on $\RPlusSpace$ in triangle length coordinates given starting points corresponding to trace coordinates $(3,3,3)$ (left), $\left(2\sqrt{2},2\sqrt{2},4\right)$ (centre), and $\left(10, 10, -10(-5 + \sqrt{23})\right)$ (right).}
	\label{fig:coordCP}
\end{figure}

\begin{figure}[htb]
	\centering
	\includegraphics[width=0.28\textwidth]{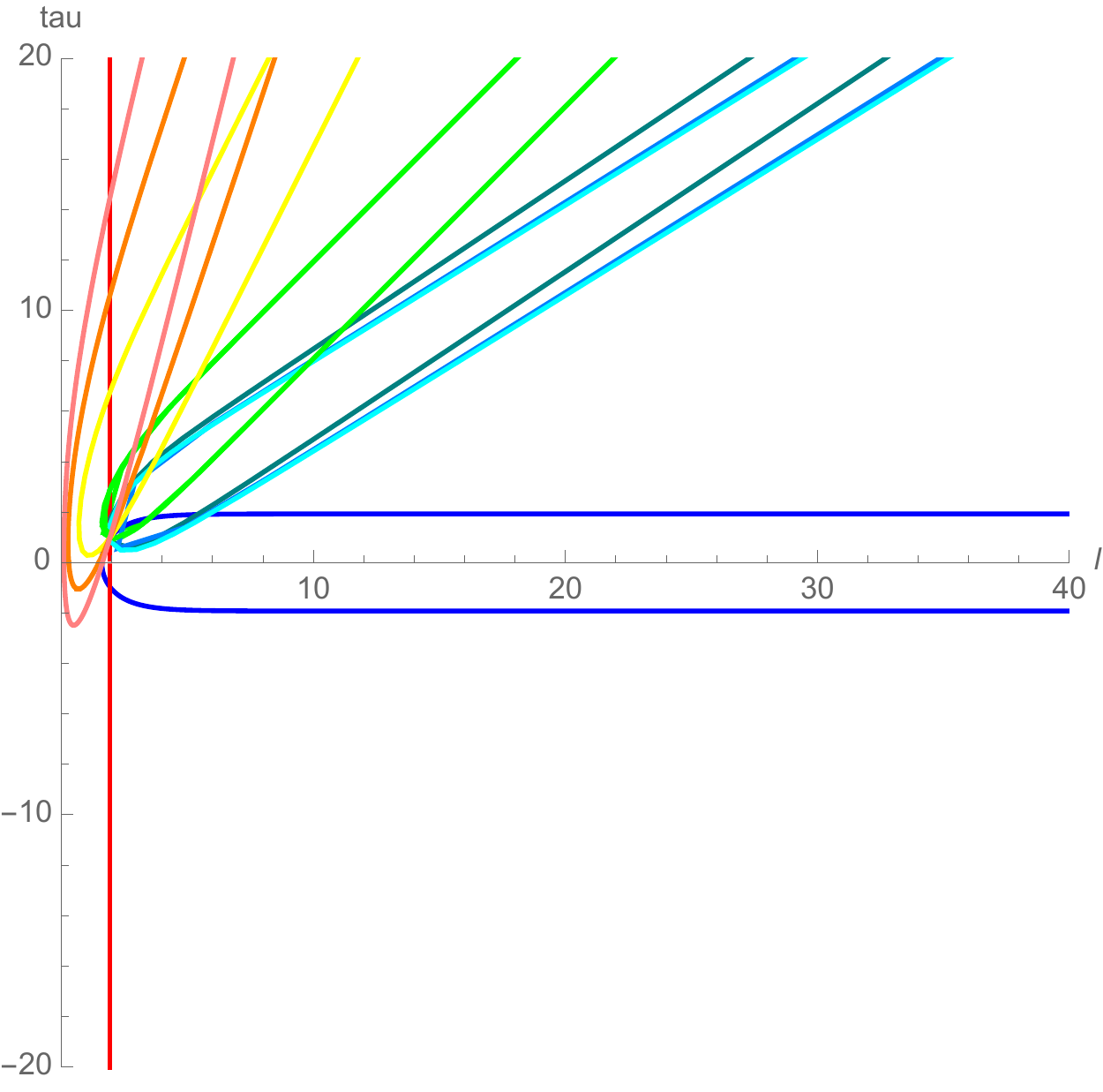}
	\includegraphics[width=0.28\textwidth]{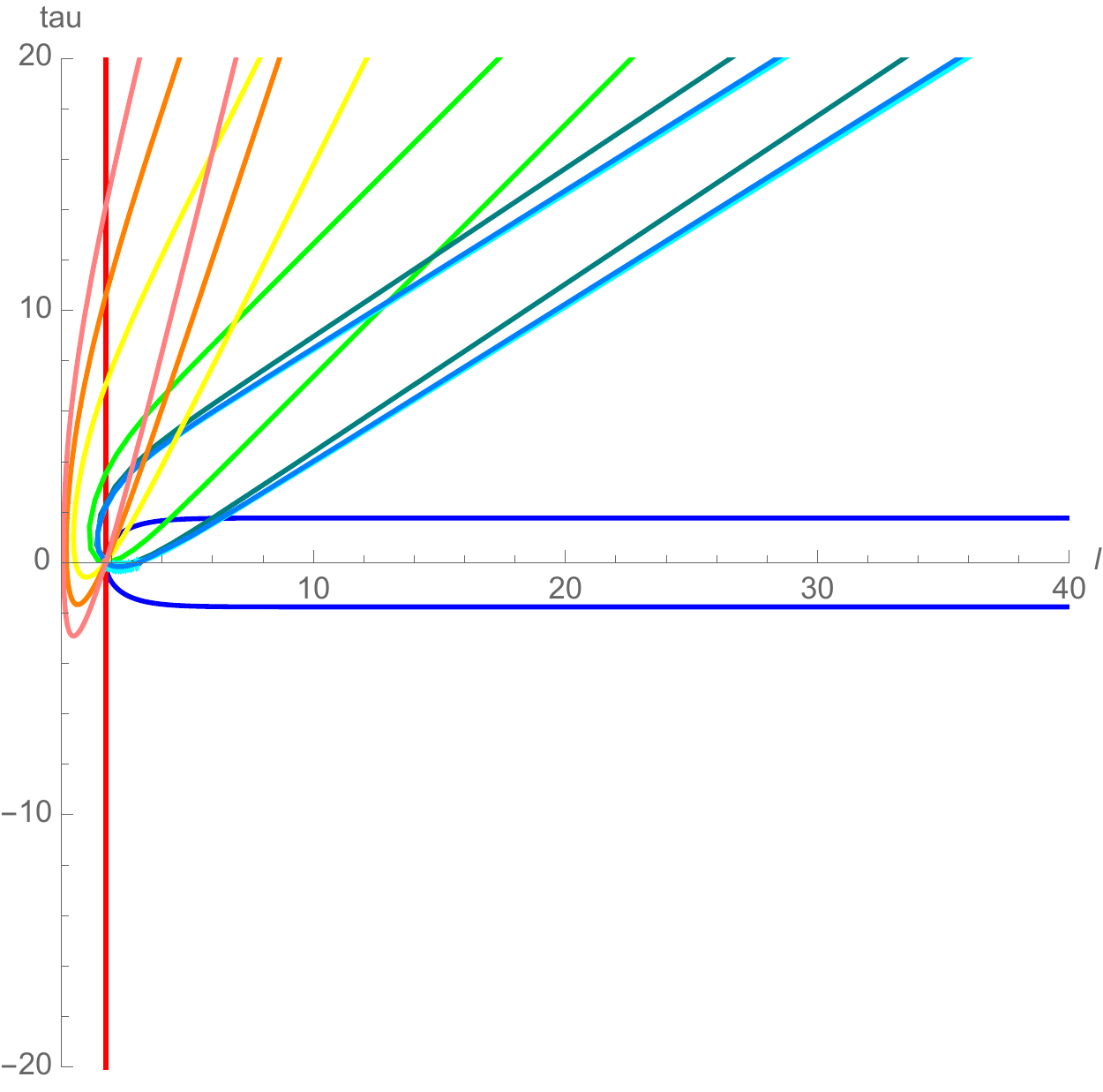}
	\includegraphics[width=0.28\textwidth]{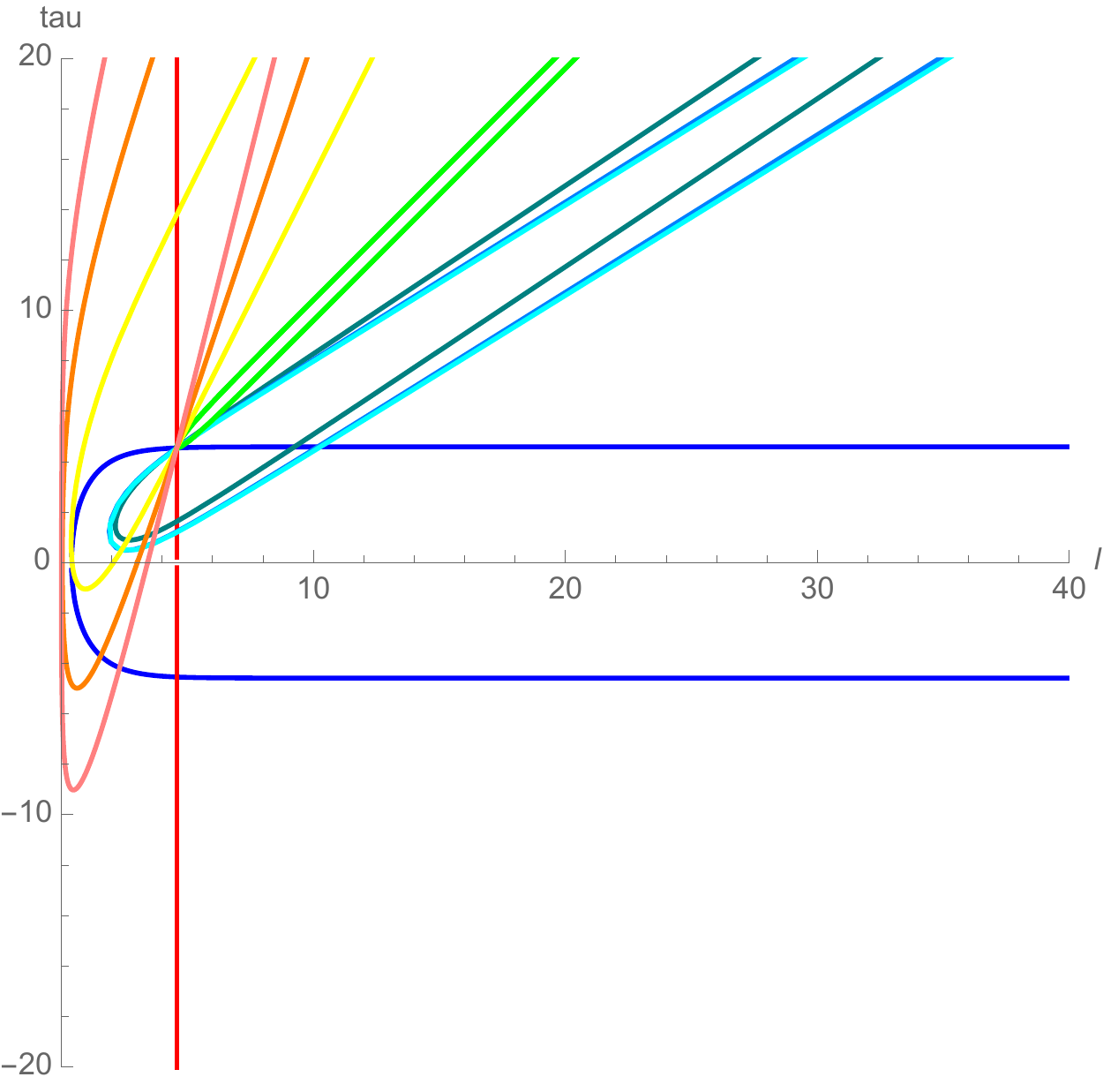}
	\caption{Forward and backward earthquake deformations about curves in $\curves$ in Fenchel-Nielsen coordinates given starting points corresponding to trace coordinates $(3,3,3)$ (left), $\left(2\sqrt{2},2\sqrt{2},4\right)$ (centre), and $\left(10, 10, -10(-5 + \sqrt{23})\right)$ (right).}
	\label{fig:coordFN}
\end{figure}

\begin{figure}[H]
	\centering
	\includegraphics[width=0.28\textwidth]{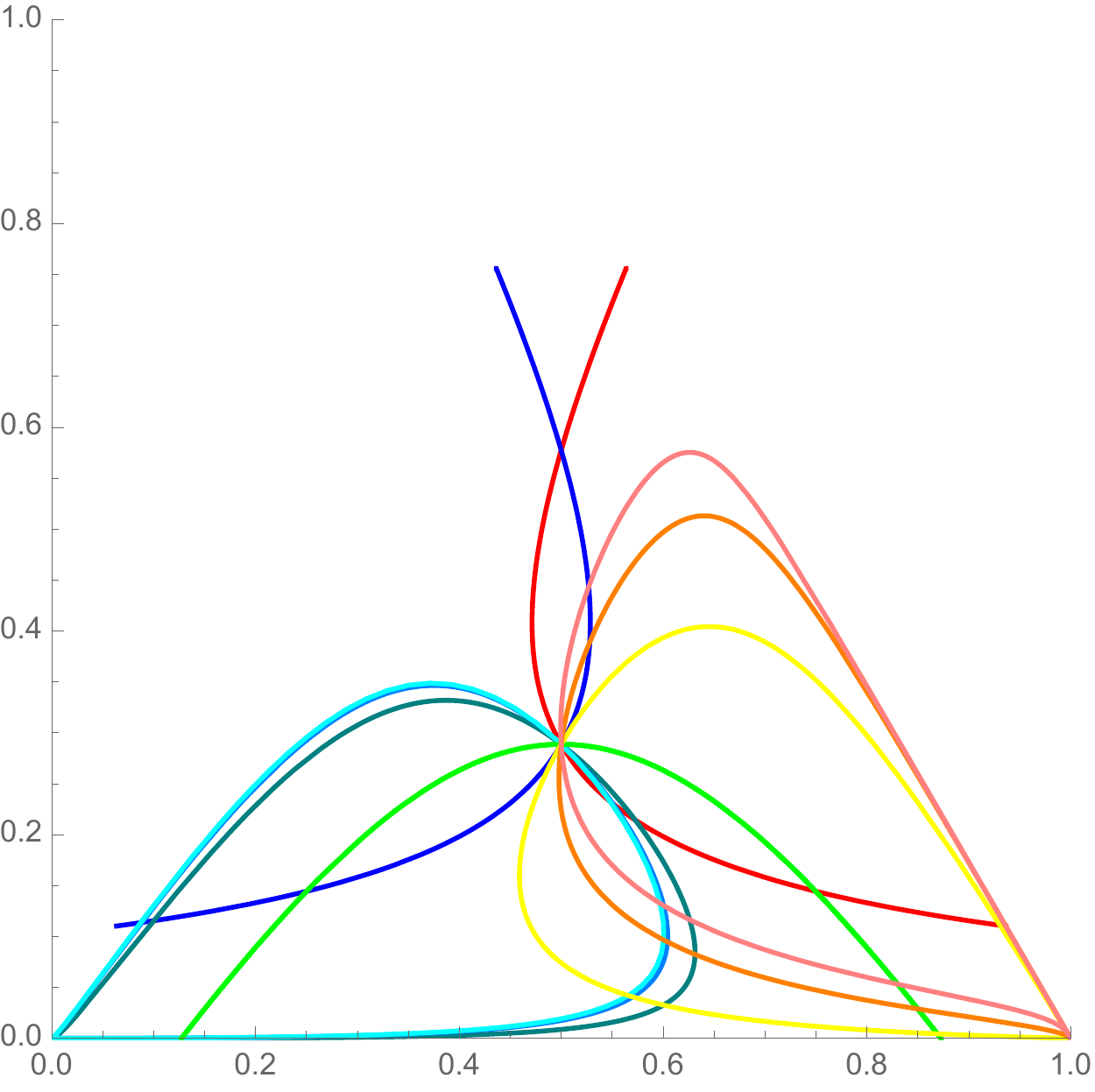}
	\includegraphics[width=0.28\textwidth]{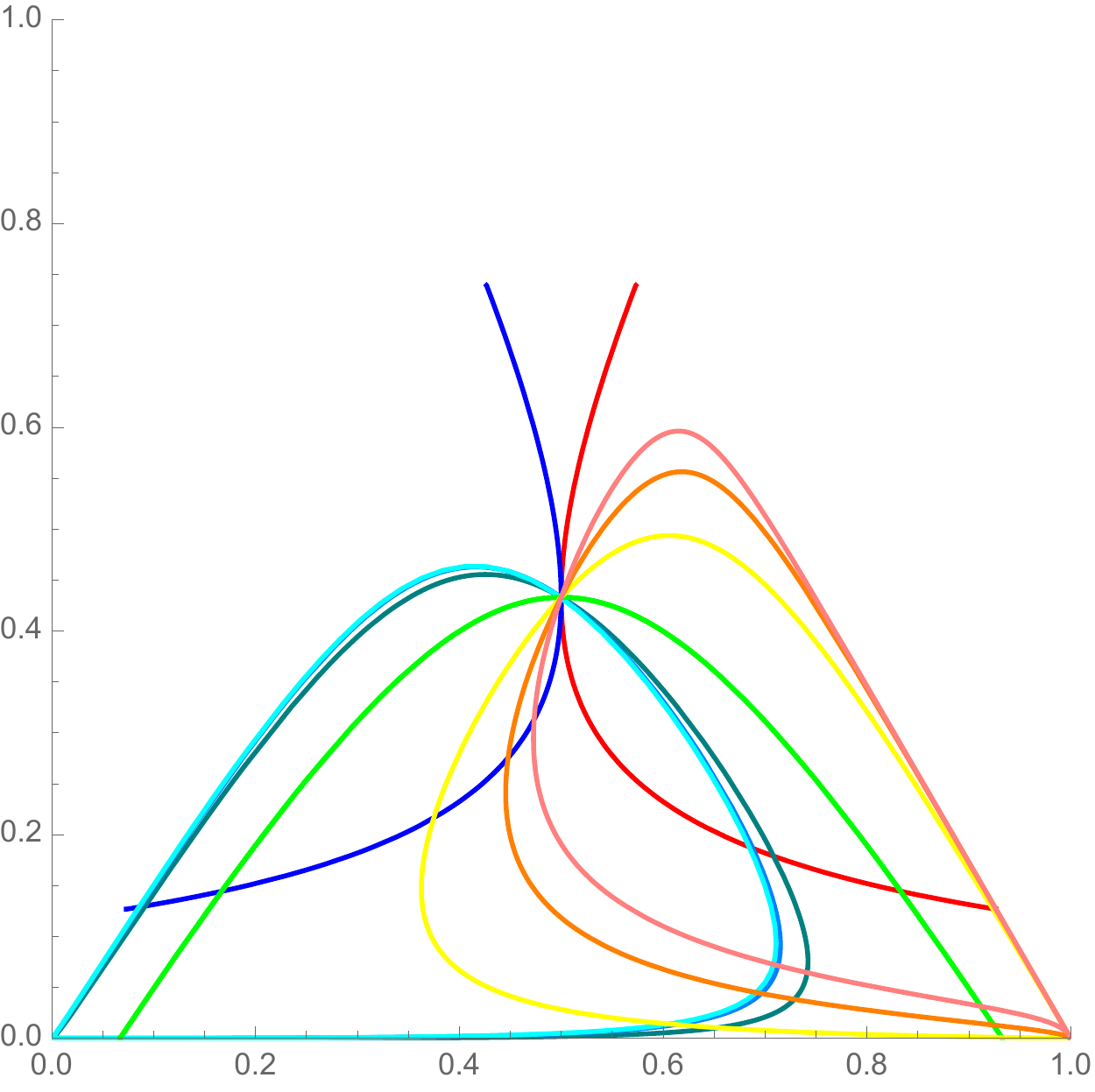}
	\includegraphics[width=0.28\textwidth]{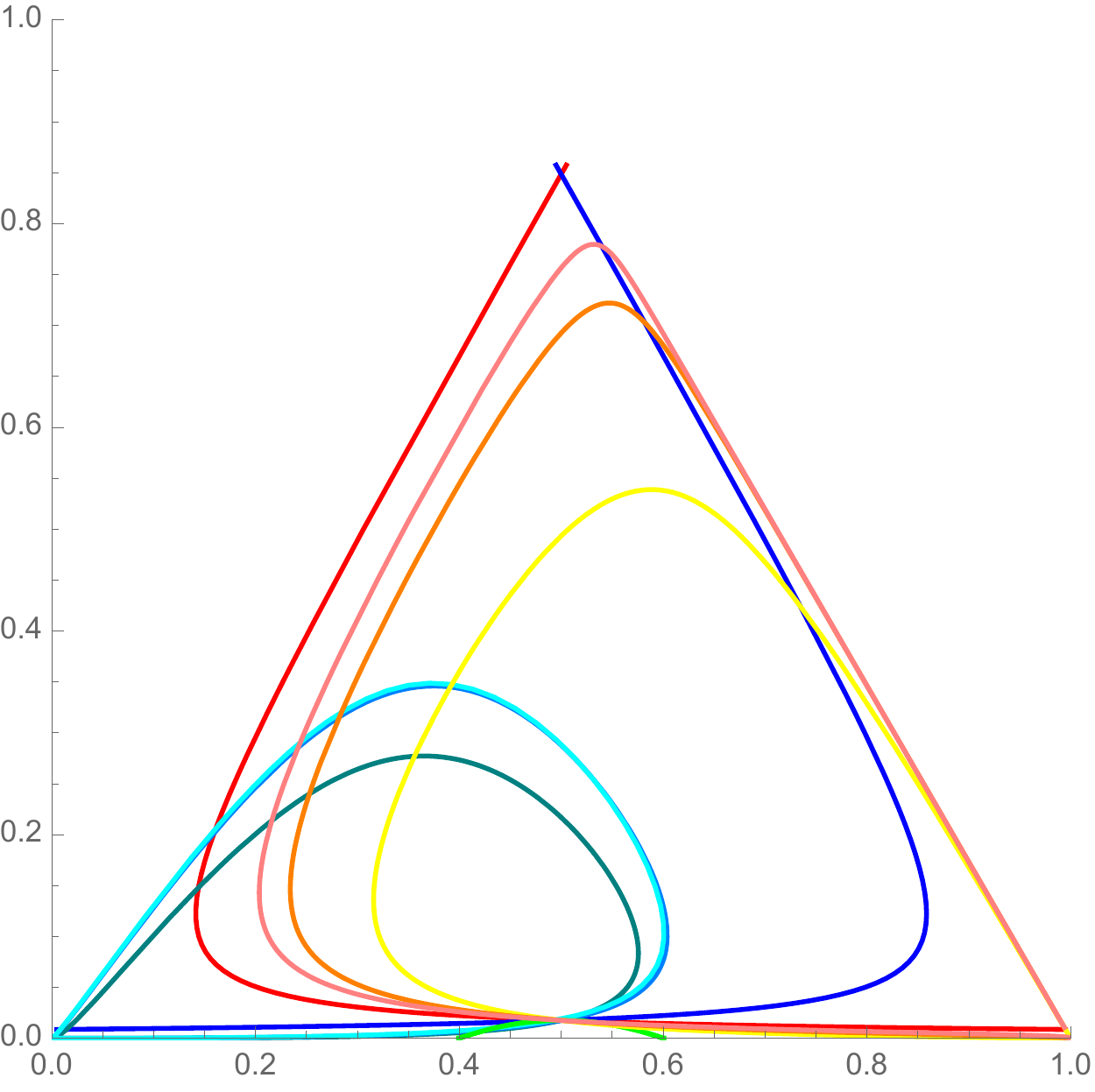}
	\caption{Forward and backward earthquake deformations about curves in $\curves$ on the simplex given starting points corresponding to trace coordinates $(3,3,3)$ (left), $\left(2\sqrt{2},2\sqrt{2},4\right)$ (centre), and $\left(10, 10, -10(-5 + \sqrt{23})\right)$ (right).}
	\label{fig:coordSim}
\end{figure}
\FloatBarrier

\emph{Inverted trace coordinates} are given by inverting the radius $r$ in the spherical trace coordinates, 
\[
	\begin{pmatrix} \theta \\ \phi \\ r \end{pmatrix} = \begin{pmatrix} \tan^{-1}\left(\frac{y}{x}\right) \\ \tan^{-1}\left(\frac{\sqrt{x^2+y^2}}{z}\right) \\ \frac{1}{\sqrt{x^2+y^2}} \end{pmatrix}
\]
with $1-\frac{1}{r}\cos(\theta)\sin(\theta)\cos(\phi)\sin(\phi)^2=0$. Pictures of earthquakes about the curves in $\curves$ in inverted trace coordinates are given in Figure \ref{fig:coordInv}. 

Recall triangle lengths as defined in Equation \ref{eqn:lengthCoord}. Pictures of earthquakes about the curves in $\curves$ in triangle length coordinates are given in Figure \ref{fig:coordCP}.

Also recall Fenchel-Nielsen coordinates using $\alpha$ are related to trace coordinates and triangle lengths by Equation \ref{eqn:FNdefn}. Pictures of earthquakes about the curves in $\curves$ in Fenchel-Nielsen coordinates are given in Figure \ref{fig:coordFN}.

Coordinates for the simplex can be defined by the transformation,
\[
	\begin{pmatrix} p \\ q \\ r\end{pmatrix} =  \begin{pmatrix} \frac{x}{yz} \\ \frac{y}{xz} \\ \frac{z}{xy} \end{pmatrix}
\]
with $p+q+r=1, p,q,r>0$. The simplex coordinates can be transformed to the plane $\left(p',q'\right)$. Pictures of earthquakes about the curves in $\curves$ on the simplex are given in Figure \ref{fig:coordSim}.

\begin{figure}[hbt]
	\centering
	\includegraphics[width=0.44\textwidth]{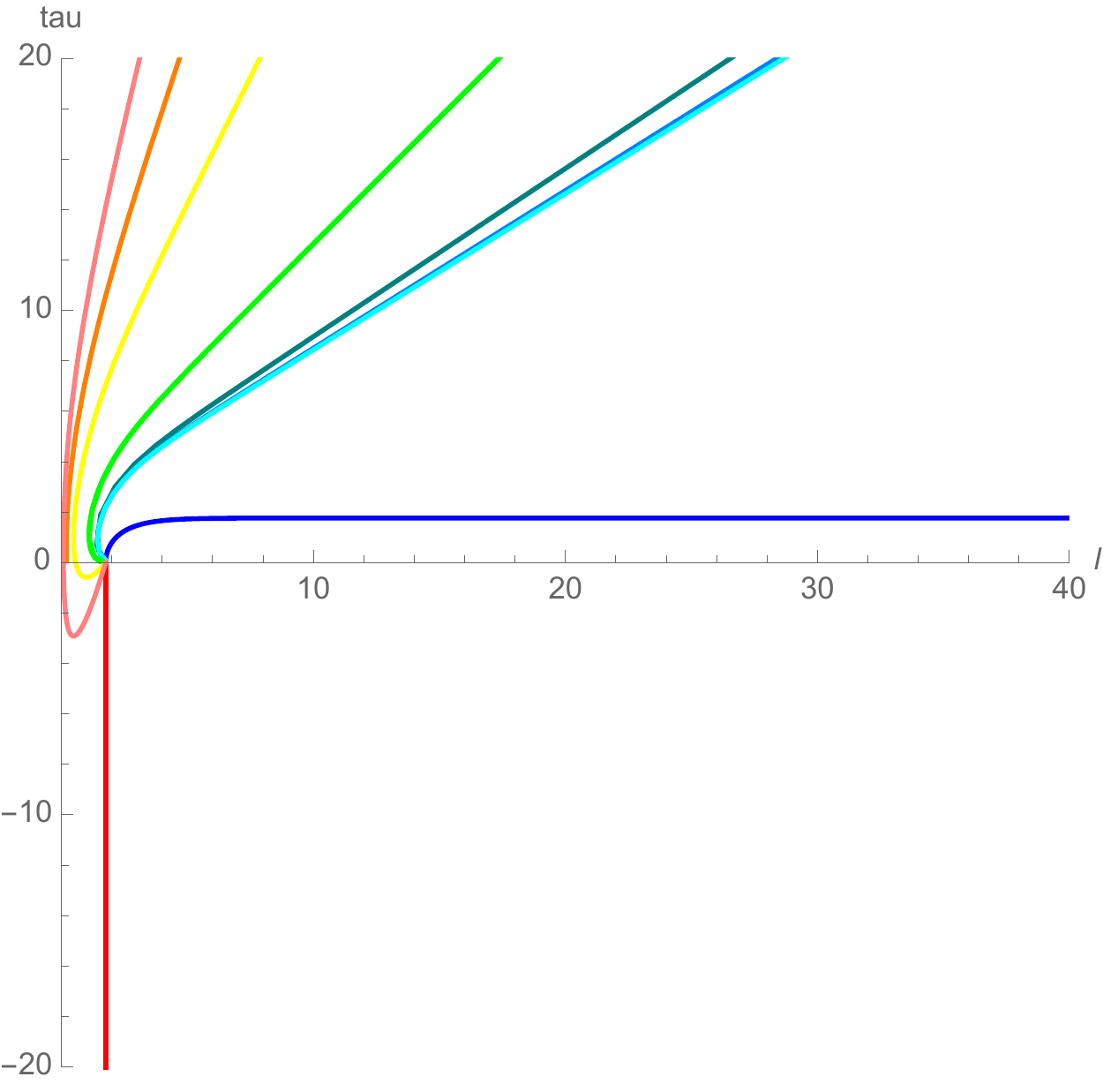}
	\includegraphics[width=0.44\textwidth]{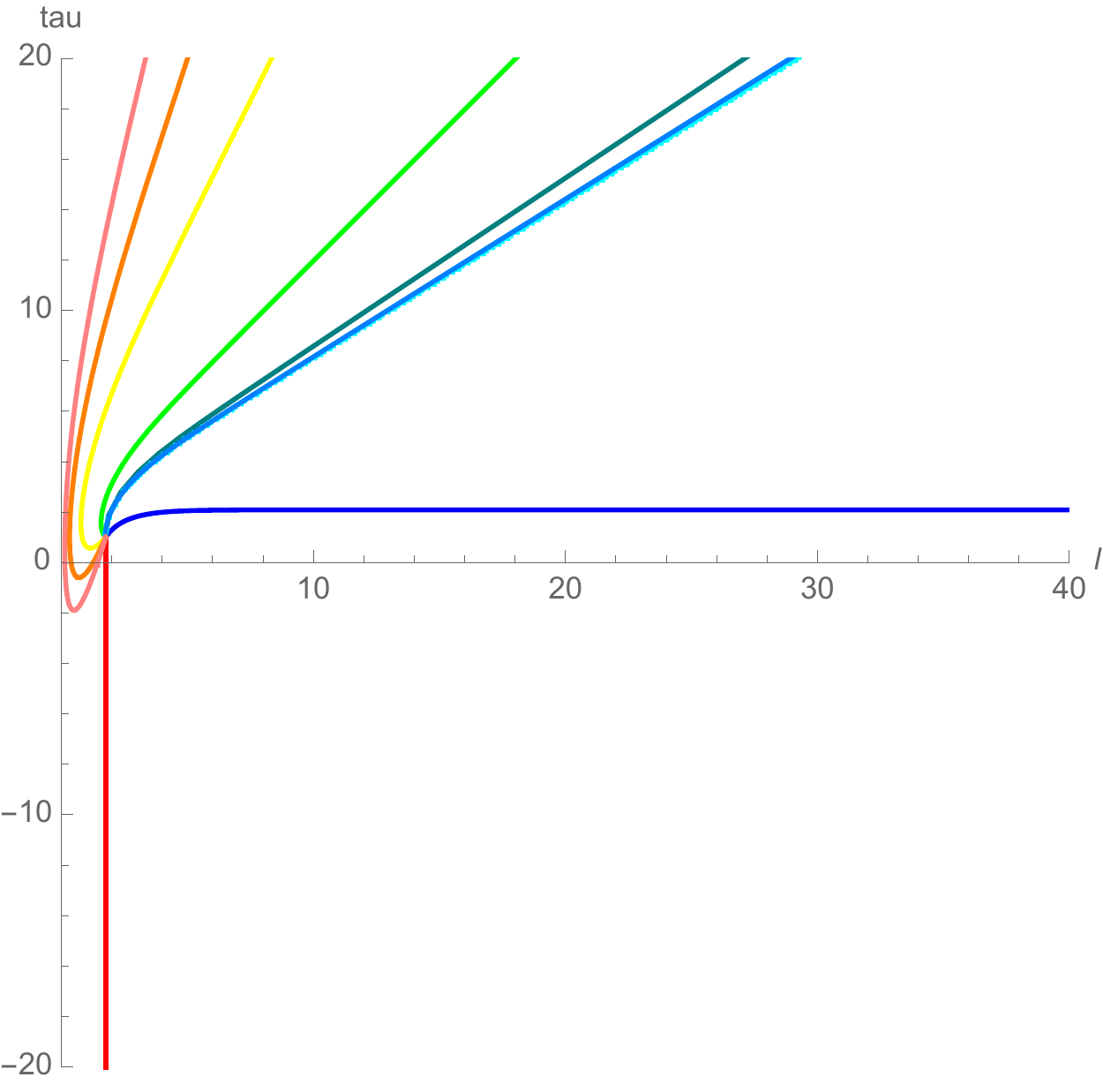}
	\caption{Forward earthquake deformations about curves in $\curves$ in Fenchel-Nielsen coordinates given starting points $\left(\cosh^{-1}(3),0\right)$ (left), $\left(\cosh^{-1}(3),1\right)$ (right).}
	\label{fig:WatWol}
\end{figure}

We can compare these pictures to previous work on earthquakes. In particular, Waterman and Wolpert \cite{WatWol83} used Fenchel-Nielsen coordinates to plot forward earthquake deformations about simple closed geodesics starting at points $\left(\cosh^{-1}\left(3\right),0\right), \left(\cosh^{-1}\left(3\right),1\right)$. Plots using our formulas converted to Fenchel-Nielsen coordinates are shown in Figure \ref{fig:WatWol} and are consistent with Waterman and Wolpert's results.

\FloatBarrier

%% --- ACKNOWLEDGEMENTS ---------------------------------------------------------------------------------------------------------------------------------------------------

\section*{Acknowledgements}
Research of the author was supported under the Australian Research Council's Discovery funding scheme (project number DP190102259). The author thanks Daryl Cooper and Catherine Pfaff for sharing an early preprint of their paper, with extra thanks to Daryl for comments on this work, Francis Bonahon for useful discussions on measured geodesic laminations, and Thomas Le Fils for insight into asymptotics of earthquakes. The author also thanks their supervisor Stephan Tillmann.

%% --- DECLARATION OF INTEREST ---------------------------------------------------------------------------------------------------------------------------------------------------

\section*{Declaration of interest statement}
The author reports there are no competing interests to declare.

%% --- REFERENCES ---------------------------------------------------------------------------------------------------------------------------------------------------

\bibliographystyle{amsplain}

\end{document}